\documentstyle[12pt]{article}
\title{Limit theorems  for  maximum flows on a lattice
\footnotetext{AMS classification: 60K 35.}
\footnotetext{Key words and phrases: maximum flow and minimum cut, random surfaces, cluster boundary, and first passage percolation.}} 
\author{Yu Zhang }
\begin{document}
\baselineskip 0.22in
\maketitle
\begin{abstract}

We  independently assign a non-negative value, as a capacity for the quantity of flows per unit time,  
with a distribution $F$ to each 
edge on the ${\bf Z}^d$ lattice. We consider the maximum flows through the edges 
from a source to a sink, in a large cube.
In this paper, we show that the ratio of the maximum flow and the size of source is asymptotic to a constant.
This constant is denoted by the flow constant.

\end{abstract}

\section {Introduction of the model and results.}
We consider the ${\bf Z}^d$ lattice, $d \geq 2,$  with integer vertices and edges between
${\bf u}=(u_1,\cdots, u_d)$ and ${\bf v}=(v_1,\cdots, v_d)$ 
when
$$\sum_{i=1}^d |u_i-v_i|=1.$$
Two vertices ${\bf u}$ and ${\bf v}$ with an edge connecting them are said to be ${\bf Z}^d$-adjacent or ${\bf Z}^d$-connected.
The edge is identified as a ${\bf Z}^d$-edge $e=({\bf u},{\bf v})$, or simply, an edge,  
with the open line segment in ${\bf R}^d$ from ${\bf u}$ to ${\bf v}$.
Two vertices ${\bf u}$ and ${\bf v}$ are said to be ${\bf L}^d$-adjacent or ${\bf L}^d$-connected if
$$\max_{1\leq i\leq d} |u_i-v_i|=1.$$
Clearly, if ${\bf u}$ and ${\bf v}$ are ${\bf Z}^d$-connected, then they are also ${\bf L}^d$-connected.

For two vertices ${\bf u}$ and ${\bf v}$, let {\bf dist}$({\bf u},{\bf v})$ be the
Euclidian distance between the two vertices.
For any two vertex sets ${\bf A}, {\bf B}\subset {\bf Z}^d$, the distance between ${\bf A}$ and ${\bf B}$ is also
defined by
$${\bf dist}({\bf A},{\bf B})=\min\{{\bf dist}({\bf u,v}):{\bf u}\in {\bf A}\mbox{ and } {\bf v}\in {\bf B}\}.\eqno{}$$
Now we assign independently to each ${\bf Z}^d$-edge $e$ a non-negative value $\tau(e)$ with a distribution $F$.
More formally, we consider the following probability space. As
sample space, we take $\Omega=\prod_{e\in {\bf Z}^d} [0, \infty),$ points of which
are represented as   {\em configurations}. If we want to emphasize a particular configuration $\omega$, 
for a random set ${\bf A}$ or  a random variable  $N$, we may write 
${\bf A}({\omega})$ or $N(\omega)$ for each, respectively.
Let ${\bf P}$ be the corresponding product measure on $\Omega$. The
expectation  with respect to ${\bf P}$ is denoted by ${\bf E}(\cdot)$.
For simplicity, we assume that $\tau(e)$ has a short tail
$$ {\bf E}\exp(\eta \tau(e))=\int_0^\infty e^{\eta x} dF(x) < \infty \eqno{(1.1)}$$
for some $\eta >0$. For each finite graph ${\bf B}$ with vertices and edges, 
we may think of $\tau(e)$ as the non-negative {\em capacity} for the quantity of fluid that may flow 
along $e\in {\bf B}$ in unit time, where an edge in a set means that the two vertices of the edge belong to
the set. Let ${\bf S}$ and ${\bf T}$ be two disjoint sets in ${\bf B}$, called the {\em source}
and the {\em sink}.
A {\em flow } (see Kesten (1987); Grimmett (1999)), 
from a vertex set ${\bf S}$ to another vertex set ${\bf T}$ in ${\bf B}$, 
is an assignment of a non-negative number $f(e)$ and an 
orientation to each edge $e=({\bf v}, {\bf w})$ of ${\bf B}$ such that
$$I({\bf v})=\sum_{{\bf w}\in {\bf B}: {\bf v}\rightarrow {\bf w}}f(({\bf v}, {\bf w}))-\sum_{{\bf w}\in {\bf B}: {\bf w}\rightarrow {\bf v}} f(({\bf v}, {\bf w}))$$
satisfies $I({\bf v})=0$ for all vertices 
${\bf v}\not\in {\bf S}\cup {\bf T}$, 
where the first summation (with respect to the second summation) 
is  calculated over all neighbors 
${\bf w}$ of ${\bf v}$, which $e({\bf v},{\bf w})$ is oriented away from (respectively toward) 
${\bf v}$. Thus fluid is conserved at all vertices except,  possibly,  at sources and sinks. 
In other words, the current flowing into a vertex ${\bf v}\not\in {\bf S}\cup {\bf T}$ must equal to the current flowing out. This basic assumption is called Kirchhoff's law in physics.
A flow is admissible if 
$$f(e) \leq \tau(e) \mbox{ for all edges }e,$$
and the value of such a flow is defined to be $\sum_{{\bf v}\in {\bf S}} I({\bf v})$, 
the aggregate amount of fluid entering ${\bf B}$ at source vertices. 
The {\em maximum flow } is the largest value of all admissible flows. One of the fundamental questions of this physics topic concerns understanding how the maximum flow depends on the source and sink.
 It is believed (see Kesten (1987); Grimmett (1999)) that the maximum flow  approximately equals  
the  ``size" of $\min\{|{\bf S}|, |{\bf T}|\}$ with a certain ratio for a convex set ${\bf B}$, where $|{\bf A}|$ is the number of vertices in ${\bf A}$. We write the ratio
as the {\em flow constant}, which only depends on  $F$. 
The main purpose of this paper is to demonstrate the existence of the
flow constant.

To study the maximum flow, we  need to understand cutsets.
To define a cutset from ${\bf S}$ to ${\bf T}$ on ${\bf B}$, we may first define a path on ${\bf Z}^d$ as follows.
For any two vertices ${\bf u}$ and ${\bf v}$ of ${\bf Z}^d$, 
a ${\bf Z}^d$ path, or simply a path, $\gamma$  from ${\bf u}$ to ${\bf v}$ is an alternating sequence 
$({\bf v}_0, e_1, {\bf v}_1,...,{\bf v}_{n-1}, e_n, {\bf v}_n)$ of vertices ${\bf v}_i$ and 
edges $e_i=({\bf v}_{i-1}, {\bf v}_i)$ in ${\bf Z}^d$, with ${\bf v}_0={\bf u}$ and $ {\bf v}_n={\bf v}$. 
${\bf u}$ and ${\bf v}$ are called ${\bf Z}^d$-connected. A ${\bf Z}^d$-connected vertex set is called {\em cluster}.
An edge  set ${\bf X}$ of  ${\bf B}$ is called an ${\bf S}$-${\bf T}$ {\em cutset } if all paths on ${\bf B}$ from ${\bf S}$ to ${\bf T}$ use at least one edge of ${\bf X}$. For convenience, we also add all the vertices of  edges in 
${\bf X}$ to have an edge and a vertex set. We still denote the set by ${\bf X}$.

A cutset ${\bf X}$ is said to be {\em self-avoiding} if ${\bf X}$ is a cutset, and 
${\bf X}\setminus \{e\}$ will no longer  be a cutset  for every $e\in {\bf X}$.
Note that there might  be many self-avoiding cutsets. 
For any edge set ${\bf E}$, we denote the {\em passage time } of ${\bf E}$ by
$$\tau({\bf E})=\sum_{e\in {\bf E}} \tau(e).$$

One of the fundamental problems in percolation is to study the cutset. 
Cutsets are also related to the boundary of clusters. For each edge $e$, it is said
to be {\em open} or {\em closed} if $\tau(e)>0$ or $\tau(e) =0$. Clearly, for each $e$,
$${\bf P}[e\mbox{ is closed }]=F(0), \mbox{ and } {\bf P}[e\mbox{ is open} ]=1-F(0)=p.$$
Note that if $e$ is closed, its passage time is zero, so we also sometimes denote it as a zero-edge.
Let ${\bf C}(x)$ be an open cluster containing $x$ and let
$$\theta(p)={\bf P}[|{\bf C}({\bf 0})|=\infty], \mbox{ and } p_c=p_c(d)=\sup\{p: \theta(p)=0\}.$$
If $F(0) < 1-p_c$, there exists an infinite open cluster from the origin with a positive probability.
If $|{\bf C}({\bf 0})|$ is finite, there exists a closed cutset that cuts the origin from $\infty$.
An edge $e$ is called the boundary edge of ${\bf C}({\bf 0})$ if $e\not\in {\bf C}({\bf 0})$, but 
$e$ is ${\bf Z}^d$-adjacent to ${\bf C}({\bf 0})$. $\Delta {\bf C}({\bf 0})$ is defined as all the  boundary edges of 
${\bf C}({\bf 0})$. If $|{\bf C}({\bf 0})|$ is finite, 
then $\Delta {\bf C}({\bf 0})$ is a 
finite closed cutset.  Here we define the more general boundary edges of  open clusters, starting
at a large set. Let 
$${\bf B}({\bf k},m)=\prod_{i=1}^{d-1}[0,k_i]\times [0,m]\mbox{ for }{\bf k}=(k_1,\cdots, k_{d-1}).$$
We may also assume without  loss of generality that
$$0\leq k_{1} \leq k_{2}\leq \cdots \leq k_{d-1}. \eqno{(1.2)}$$
When $k_1=k_2=\cdots=k_{d-1}=m=0$, ${\bf B}({\bf k},m)$ is  the origin.
We also denote  the volume of $[0,k_1]\times \cdots [0,k_{d-1}]$ by
$$\|{\bf k}\|_v=k_1\times k_2\times \cdots \times k_{d-1}.$$
As we defined above, a set is said to be  a cutset that cuts  ${\bf B}({\bf k},m)$ from $\infty$ if any path from 
${\bf B}({\bf k},m)$ to $\infty$  uses at least one edge of the set. 
We select,  from 
these cutsets, a cutset  ${ \bf X}({\bf k}, m)$ with the minimum passage time among  all the
cutsets. We also denote by $\chi({\bf k}, m)$ the passage time of 
${ \bf X}({\bf k}, m)$:
$$\tau({\bf X}({\bf k}, m))=\chi ({\bf k}, m).$$
There might be many such cutsets.
 If so, we select the one with the  minimum number of edges among all such cutsets by using a unique method. 
We still denote it by ${\bf X}({\bf k}, m)$.
With this selection, ${\bf X}({\bf k},m)$ must be self-avoiding. 
In this paper, the unique method of selecting cutsets always involves using  the same selection rule for each configuration.

Furthermore, a set ${\bf X}({\bf k},m)$ is said to be  a zero-cutset (or a closed cutset) that cuts  ${\bf B}({\bf k},m)$ from $\infty$
 if 
$$\tau({\bf X}({\bf k},m))=0.$$
In this case,  any path from 
${\bf B}({\bf k},m)$ to $\infty$ must use at least one closed edge of the set. 
In other words, there is no open path from ${\bf B}({\bf k},m)$ to $\infty$.
Let $N({\bf k},m)$ be the number of edges in ${\bf X}({\bf k}, m)$. 
We have the following  fundamental geometric theorem to show that $N({\bf k}, m)$ cannot be
 much larger than $\|{\bf k}\|_v$ when $F(0) < 1-p_c$.\\

{\bf Theorem 1.} {\em If $F(0) < 1-p_c$,
then there exist constants $\beta=\beta(F, d)$ and $C_i=C_i(F, \beta, d)$
for $i=1,2$ such that for all $n \geq \beta  \|{\bf k}\|_v$ and $m \leq \min_{1\leq i\leq d-1}k_i$},
$${\bf P}\left[ \infty > N({\bf k},m) \geq n\right ]\leq C_1\exp(-C_2 n).$$

In this paper, we always denote by $C$ or $C_i$  a large or a small positive constant that will be used for some upper or lower bound in inequalities.
$C$ and $C_i$ do not depend on $k_1,\cdots, k_{d-1}$, $m$,
$l$, $n$, $w_1,\cdots, w_{d-1}$. 
In addition, both values of $C$ and $C_i$ may  change from  appearance to  appearance. 
For a finite open cluster ${\bf C}({\bf 0})$, 
its exterior boundary edges are the edges in   $\Delta {\bf C}({\bf 0})$ such that
there is a ${\bf Z}^d$-path from the vertices  of the edges to $\infty$ without using any edges of ${\bf C}({\bf 0})$.
We denote by $\Delta_e {\bf C}({\bf 0})$ the exterior boundary of ${\bf C}({\bf 0})$.
Kesten and Zhang (1990) showed that there exists a constant $\sigma$ such that
$$\lim_{n\rightarrow \infty} -n^{-1} \log {\bf P}\left[ |\Delta_e {\bf C}({\bf 0})|=n\right] =\sigma(F(0)).$$
By (6.18) in Grimmett (1999), we know that 
$$\sigma (1-p_c)=0.\eqno{(1.3)}$$
It also follows from (6.13) in Grimmett (1999) that if $F(0) >1-p_c$, then
$$\sigma(F(0)) >0.\eqno{(1.4)}$$
It is natural to ask whether (1.4)  holds when $F(0) < 1-p_c$. 
Note that $\Delta_e {\bf C}({\bf 0})$ is a closed cutset for ${\bf C}({\bf 0})$, so on the event that there is no
infinite open cluster,
$$\mbox{ the number of edges in }\Delta_e {\bf C}({\bf 0})\geq N({\bf 0}, 0).$$
However, the size of $N({\bf 0}, 0)$ may be much less than the size of $\Delta_e {\bf C}({\bf 0})$
when $ {\bf C}({\bf 0})$ is finite. Thus, we still do not know whether $\sigma (F(0)) >0$ when $F(0) < 1-p_c$.

Now we focus on a specific cutset. In fact,
one of most interesting questions (see Kesten (1987); Grimmett (1999)) 
is to understand the
behavior of the cutsets on ${\bf B}({\bf k}, m)$ that cut its bottom face from its top face.
We denote by
$${\bf F}_0={\bf F}_0({\bf k}, m)=\{(x_1,\dots, x_d)\in {\bf B}({\bf k}, m): x_d=0\}$$
and 
$${\bf F}_m={\bf F}_m({\bf k}, m)=\{(x_1,\cdots, x_d)\in {\bf B}({\bf k}, m): x_d=m\}$$
 the bottom and the top faces of the box, respectively.
We select ${\bf W}({\bf k}, m)$ as a  cutset,  
cutting the bottom face of ${\bf B}({\bf k}, m)$ from its top face,
with  minimal passage time. Similarly, if there is more than one such cutsets,
we use  a unique method to select one
with the minimum number of  edges among all such cutsets.
We still denote it by ${\bf W}({\bf k}, m)$.
Let $\bar{N}({\bf k}, m)$ be the number of vertices in  cutset ${\bf W}({\bf k}, m)$. 
We now show the fundamental geometric theorem for this cutset.\\

{\bf Theorem 2.} {\em If $F(0) < 1-p_c$, and if $m=m(k_1,\cdots, k_{d-1})\rightarrow \infty $ as $k_{1}, k_{2}, \cdots, k_{d-1}\rightarrow \infty$
in such a way that
$$\log m\leq \|{\bf k}\|_v, \eqno{(1.5)}$$
then there exist constants $\beta=\beta(F, d)\geq 1$ and $C_i=C_i(F, \beta,  d)$
for $i=1,2$ such that for all $n \geq \beta \|{\bf k}\|_v$},
$${\bf P}\left[ \bar{N}({\bf k}, m) \geq n\right ]\leq C_1\exp(-C_2 n).\eqno{(1.6)}$$

{\bf Remark 1.} In the proof of Theorem 2, we can use a weak condition that 
$\log m\leq C\|{\bf k}\|_v$ to replace (1.5).\\

{\bf Remark 2.} In Theorems 1 and  2, we consider a cutset that cuts ${\bf B}({\bf k}, m)$ from $\infty$ or from
${\bf F}_0$ to ${\bf F}_m$. The same proof
can be shown for a general set rather than ${\bf B}({\bf k}, m)$. \\

With Theorem 2,
the number of vertices for each cutset is proportional to the size of  ${\bf F}_0$.
We call the results in Theorem 2  the {\em linearity}.
When $F(0) > 1-p_c$, it is known (see chapter 6 in Grimmett (1999)) that  Theorems 1 and  2 hold. 
On the other hand, it is known by Aizenman et al. (1983) 
that if $F(0) < 1-p_c$ and $m$ satisfies (1.5), then
$${\bf P}[ \exists \,\,\, {\bf W}({\bf k}, m) \mbox{ with } \tau({\bf W}({\bf k}, m))=0] \leq \exp(-C \|{\bf k}\|_v).
\eqno{(1.7)}$$
(1.7) is called the {\em area law}. Clearly,  $\bar{N}({\bf k}, m)$ is always larger than $\|{\bf k}\|_v$.  In fact, we may view the cutset (see Aizenman et al. (1983)) as a surface
between ${\bf F}_m$ and  ${\bf F}_0$. (1.6) tells us that it costs probability $\exp(-Ct)$ whenever the surface
increases $t$ units.
We call (1.6) the  {\em surface law}. The surface law has proved
to hold  (Kesten (1986) and (1988)) when $d=2$ and $F(0) < 1-p_c(2)$,  and when $d=3$ and $F(0)< 1/27$. 
As the main conjecture, Kesten believed that the surface law should hold for
all $d$ and all $F(0) < 1-p_c(d)$. In  Theorem 2, we answer Kesten's conjecture
affirmatively.
When $F(0)=1-p_c$, the closed cutsets are very chaotic. For example, we believe 
that $\bar{N}({\bf k}, m)$ should be  much larger than $\|{\bf k}\|_v$.

Now we focus on the maximum flow problem to discuss the existence of the flow constant. 
Without loss of generality (see Kesten (1987); Grimmett (1999)), we discuss the maximum flow on ${\bf B}({\bf k}, m)$ from ${\bf F}_0$ to
${\bf F}_m$.
The max-flow  min-cut theorem characterizes the maximum flow through the network in terms of the sizes of cutsets.  
The size 
of the $({\bf F}_0, {\bf F}_m)$-cutset ${\bf W}({\bf k}, m) $ is defined to be the sum of the capacities of edges in ${\bf W}({\bf k}, m)$. 
As we mentioned, one of fundamental  questions (see Kesten (1987); Grimmett (1999)) is how to understand the limit behavior of  the flow from ${\bf F}_0$ in ${\bf B }({\bf k}, m)$. 
 Let $\phi_{\max}({\bf k},m)$ denote the maximum flow through the edges of ${\bf Z}^d$ in ${\bf B}({\bf k},m)$
from ${\bf F}_0$ to ${\bf F}_m$.  Let 
$$\tau_{\min} ({\bf k}, m)=\tau ({\bf W}({\bf k}, m)).$$

By the  max-flow min-cut theorem, we have
$$\tau_{\min}({\bf k}, m)= \phi_{\max}({\bf k}, m).$$
In particular, if $\tau(e)$ only takes 0 or 1, the maximal flow $\phi_{\max} ({\bf k}, m)$ is the number of disjoint
open paths from ${\bf F}_0$ to ${\bf F}_m$ in ${\bf B}({\bf k}, m)$.

With these definitions, let us introduce the developments in this field.
When $F(0)=1-p_c$, the so-called  critical case, it has been proved (see Zhang (2000)) that
$$\lim_{k_1, \cdots, k_{d-1}, m\rightarrow \infty} {\tau_{\min}({\bf k}, m)\over \|{\bf k}\|_v}=0\mbox{ a.s. and in }L_1.\eqno{(1.8)}$$
When $F(0) > 1-p_c$, the so-called  supercritical case, we also have
$$\lim_{k_1,\cdots, k_{d-1}, m\rightarrow \infty} {\tau_{\min}({\bf k}, m)\over \|{\bf k}\|_v}=0\mbox{ a.s. and in }L_1.\eqno{(1.9)}$$
In fact, as we mentioned before (see chapter 6 in Grimmett (1999)),  with a large probability,
$$\tau_{\min}({\bf k}, m)=0 \mbox{ when }m \geq k_{d-1}^{\delta} \mbox{ for } \delta>0.\eqno{}$$
In other words, the flow constant will vanish in the supercritical and critical cases.

The most interesting case is  understanding the limit behaviors
when $F(0) < 1-p_c$, the  subcritical case. 
With the moment assumption in (1.1), we have
$$\limsup _{{\bf k}\rightarrow \infty} {{\bf E}\tau_{\min}({\bf k}, m)\over \|{\bf {\bf k}\|_v}} < \infty. \eqno{(1.10)}$$
In fact, by a standard large deviation estimate, we can show  that
$$P(\tau_{\min} ({\bf k}, m)\geq C \|{\bf k}\|_v)\leq C_1\exp(-C_2 \|{\bf k}\|_v).\eqno{(1.11)}$$
On the other hand, it can be  shown (see Chayes and Chayes (1986)) that when $F(0) < 1-p_c$,
$$0< \liminf _{{\bf k}\rightarrow \infty} {{\bf E}\tau_{\min}({\bf k}, m)\over \|{\bf k} \|_v}.\eqno{(1.12)}$$
With (1.10) and (1.12), it is natural to ask what  the limit behavior is. If the limit exists, then the flow
constant exists.
When $d=2$,
 Grimmett and Kesten (1984) showed that 
$$\lim_{k_1,m\rightarrow \infty}(k_1)^{-1} \tau_{\min} (k_1, m) =\nu(F)  \mbox{ a.s. and in }L_1\eqno{(1.13)}$$
when $k_1\rightarrow \infty, m\rightarrow \infty$ such that
$$\log m /k \rightarrow 0.\eqno{(1.14)}$$
In fact, when $d=2$, the min-cutset is just a dual path from the left to the right in ${\bf B}(k_1, m)$.
The techniques to handle paths have been well developed since Hammersley and Welsh created the first passage percolation model in 1965.

When $d=3$,
Kesten (1987) used  a surface consisting of a {\em plaquette} (see Aizenman et al. (1983); Kesten (1987))
to work on the limit behavior of $\tau_{\min}({\bf k}, m)$. 
He showed, in an extensive  proof, that if the surface law  holds, then
$$\lim_{k_1,  k_2, m\rightarrow \infty}(k_1\times k_2)^{-1} \tau_{min} \left((k_1, k_2), m\right) =\nu(F) \mbox{ a.s. and in }L_1\eqno{(1.15)}$$
when $k_1, k_2\rightarrow \infty,$ $ m(k_1, k_2)\rightarrow \infty$ as $k_1 \leq k_2$ in such a way for some $\delta >0$
such that
$$ \log m \leq k_1^{1-\delta}.\eqno{(1.16)}$$
Furthermore, he showed that the surface law holds  when   $d=3$ and  $F(0) < 1/27$. 
Therefore, his result implies that the flow constant exists when $d=3$ and $F(0) < 1/27$.
Kesten conjectured that the surface law should hold for all $F(0)< 1-p_c$. 
In Theorem 2, we show the surface law. Thus,  the flow constant exists for all $F(0) \leq 1-p_c$ when $d=3$.
In addition, Kesten also conjectured that
the flow constant should exist for all $d \geq 2$. In this paper, we answer the conjectures affirmatively to  show the existence of the flow constant for all $F$.
\\

{\bf Theorem  3.} {\em If  (1.1) holds, and if $m=m(k_1,\cdots, k_{d-1})\rightarrow \infty $ as $k_{1},k_{2}, \cdots, k_{d-1}\rightarrow \infty$
in such a way that
$$\log m/\|{\bf k}\|_v \rightarrow 0, \eqno{(1.17)}$$
then there exists a flow constant $\nu(F)$ such that}
$$\lim_{k_1, k_2,\cdots , k_{d-1}, m\rightarrow \infty} {\tau_{\min} ({\bf k}, m)\over \|{\bf k}\|_v}=
\lim_{k_1, k_2,\cdots , k_{d-1}, m\rightarrow \infty} {\phi_{\max} ({\bf k}, m)\over \|{\bf k}\|_v}=\nu(F)\mbox{ a.s. and in }
L_1.$$

{\bf Remark 3.}  Kesten believes (open problem 2.24 in Kesten (1987)) that
the second moment condition for $\tau$ should imply Theorem 3. Since we need to use a concentration inequality (see (9.1) below), we need the exponential tail assumption for
$\tau$ in (1.1). On the other hand, the condition of (1.17) is optimal. \\

{\bf Remark 4.} As we discussed before, if $\nu(F)$ exists in Theorem 3,  $\nu(F)$ is positive if and
only if $F(0) < 1- p_c$.
By using the concentration inequality (9.22) below, we might have an upper estimate for the variance of $\tau_{\mbox{min}}({\bf k}, m)$. However, an asymptotical estimate of the variance should be much more difficult to achieve. We can also estimate convergence rates for the limit in Theorem 3 by using the concentration inequality (9.22),
 With this convergence rate, we can show the continuity of the flow constant  $\nu(F)$ in $F$ when $F(0) < 1-p_c$. Note that the continuity of $\nu(F)$, when $F(0)=1-p_c$, is proved by Zhang (2000).\\

{\bf Remark 5.}  Theorem 3 can be generalized to any periodic lattice (see the definition in Kesten (1982)) with the $d$ coordinate axes.\\

{\bf Remark 6.} Cerf and Theret (2011) and Raphael and Theret (2010) showed the large deviations for the maximum flows.

\section{A construction for a linear cutset.}
In this section, we will construct  a special zero-, or closed, cutset
 about  the linear size of
$\|{\bf k}\|_v$. 
Since we only consider the cutset surrounding ${\bf B}({\bf k}, m)$, for convenience, 
we will assume that all edges inside ${\bf B}({\bf k},m)$ are open edges in this section.
Now our probability
measure is on the edges in ${\bf Z}^d\setminus {\bf B}({\bf k},m)$. 

For a finite ${\bf Z}^d$-connected set ${\bf A}$, 
$\partial {\bf A}$ is a vertex set, called the {\em boundary} of ${\bf A}$, that is ${\bf L}^d$-adjacent to ${\bf A}$ but is not in ${\bf A}$. 
We also denote by $\partial_i {\bf A}$  the vertex set, called the {\em interior boundary} of ${\bf A}$,
that is in ${\bf A}$ and is ${\bf L}^d$-adjacent to  $\partial {\bf A}$.
Furthermore,
we name  $\partial_e {\bf A}$ as its {\em exterior boundary} if its vertex $v\in \partial {\bf A}$
and there is a ${\bf Z}^d$-connected path from $v$ to $\infty$ that does not use vertices of ${\bf A}$ (see Fig. 1).
Note that 
$$\partial_e {\bf A}\subset \partial {\bf A}.$$
Recall from section 1 that $\Delta {\bf A}$  and $\Delta_e {\bf A}$ are defined as the ${\bf Z}^d$-edges from 
$\partial {\bf A}$ and from $\partial_e {\bf A}$, respectively, to $\partial_i {\bf A}$. 
They are  called {\em boundary edges} and 
{\em exterior boundary edges}. 
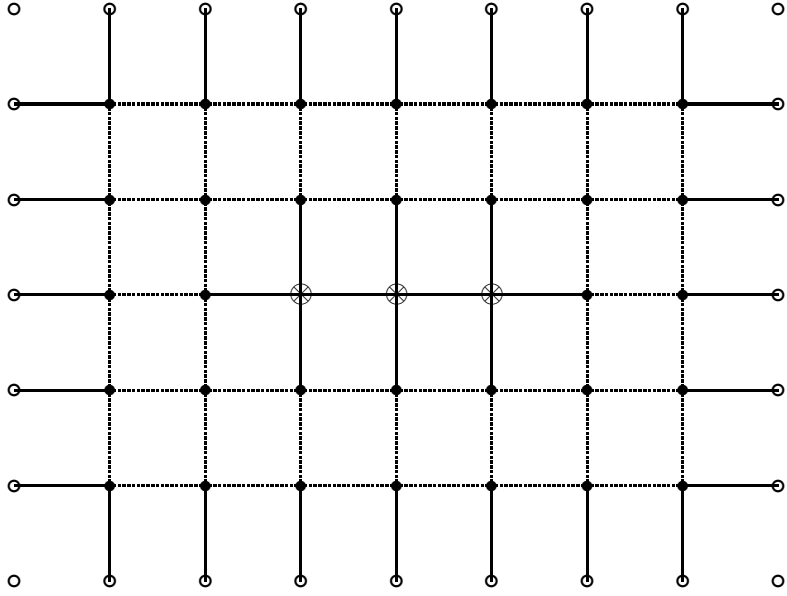
\begin{figure}\label{F:alphabeta}
\begin{center}
\setlength{\unitlength}{0.0125in}%
\begin{picture}(200,150)(67,800)
\thicklines
\put(40,760){\dashbox(40,40)[br]}
\put(80,760){\dashbox(40,40)[br]}
\put(120,760){\dashbox(40,40)[br]}
\put(120,760){\dashbox(40,40)[br]}
\put(160,760){\dashbox(40,40)[br]}
\put(200,760){\dashbox(40,40)[br]}
\put(240,760){\dashbox(40,40)[br]}

\put(40,880){\dashbox(40,40)[br]}
\put(80,880){\dashbox(40,40)[br]}
\put(120,880){\dashbox(40,40)[br]}
\put(120,880){\dashbox(40,40)[br]}
\put(160,880){\dashbox(40,40)[br]}
\put(200,880){\dashbox(40,40)[br]}
\put(240,880){\dashbox(40,40)[br]}

\put(40,800){\dashbox(40,40)[br]}
\put(40,840){\dashbox(40,40)[br]}

\put(240,800){\dashbox(40,40)[br]}
\put(240,840){\dashbox(40,40)[br]}

\put(0,720){\circle{4}}
\put(0,760){\circle{4}}
\put(0,800){\circle{4}}
\put(0,840){\circle{4}}
\put(0,880){\circle{4}}
\put(0,920){\circle{4}}
\put(0,960){\circle{4}}

\put(320,720){\circle{4}}
\put(320,760){\circle{4}}
\put(320,800){\circle{4}}
\put(320,840){\circle{4}}
\put(320,880){\circle{4}}
\put(320,920){\circle{4}}
\put(320,960){\circle{4}}

\put(40,720){\circle{4}}
\put(80,720){\circle{4}}
\put(120,720){\circle{4}}
\put(160,720){\circle{4}}
\put(200,720){\circle{4}}
\put(240,720){\circle{4}}
\put(280,720){\circle{4}}

\put(40,960){\circle{4}}
\put(80,960){\circle{4}}
\put(120,960){\circle{4}}
\put(160,960){\circle{4}}
\put(200,960){\circle{4}}
\put(240,960){\circle{4}}
\put(280,960){\circle{4}}

\put(40,760){\circle*{4}}
\put(40,800){\circle*{4}}
\put(40,840){\circle*{4}}
\put(40,880){\circle*{4}}
\put(40,920){\circle*{4}}

\put(80,760){\circle*{4}}
\put(80,800){\circle*{4}}
\put(80,840){\circle*{4}}
\put(80,880){\circle*{4}}
\put(80,920){\circle*{4}}

\put(280,760){\circle*{4}}
\put(280,800){\circle*{4}}
\put(280,840){\circle*{4}}
\put(280,880){\circle*{4}}
\put(280,920){\circle*{4}}

\put(240,760){\circle*{4}}
\put(240,800){\circle*{4}}
\put(240,840){\circle*{4}}
\put(240,880){\circle*{4}}
\put(240,920){\circle*{4}}

\put(200,760){\circle*{4}}
\put(200,800){\circle*{4}}
\put(195,837){$\otimes$}
\put(200,880){\circle*{4}}
\put(200,920){\circle*{4}}
\put(155,837){$\otimes$}
\put(115,837){$\otimes$}

\put(120,760){\circle*{4}}
\put(120,800){\circle*{4}}
\put(120,880){\circle*{4}}
\put(120,920){\circle*{4}}

\put(160,760){\circle*{4}}
\put(160,800){\circle*{4}}
\put(160,880){\circle*{4}}
\put(160,920){\circle*{4}}

\put(80,840){\line(1,0){160}}
\put(120,800){\line(0,1){80}}
\put(160,800){\line(0,1){80}}
\put(200,800){\line(0,1){80}}

\put(0,760){\line(1,0){40}}
\put(0,800){\line(1,0){40}}
\put(0,840){\line(1,0){40}}
\put(0,880){\line(1,0){40}}
\put(0,920){\line(1,0){40}}

\put(280,760){\line(1,0){40}}
\put(280,800){\line(1,0){40}}
\put(280,840){\line(1,0){40}}
\put(280,880){\line(1,0){40}}
\put(280,920){\line(1,0){40}}

\put(40,720){\line(0,1){40}}
\put(80,720){\line(0,1){40}}
\put(120,720){\line(0,1){40}}
\put(160,720){\line(0,1){40}}
\put(200,720){\line(0,1){40}}
\put(240,720){\line(0,1){40}}
\put(280,720){\line(0,1){40}}

\put(40,920){\line(0,1){40}}
\put(80,920){\line(0,1){40}}
\put(120,920){\line(0,1){40}}
\put(160,920){\line(0,1){40}}
\put(200,920){\line(0,1){40}}
\put(240,920){\line(0,1){40}}
\put(280,920){\line(0,1){40}}

\end{picture}
\end{center}
\vskip 1in
\caption{ \em ${\bf A}$ is a cluster containing  edges  (dashed lines) and  vertices (solid circles). 
The solid circles are also interior boundary $\partial_i ({\bf A})$.
$\partial {\bf A}$ is a vertex set (circles and otimes). $\partial_e {\bf A}$ is a subset of  $\partial {\bf A}$ (only circles). 
The four corners in $\partial_e {\bf A}$ are ${\bf L}^d$-adjacent to ${\bf A}$ 
$\Delta {\bf A}$ is an edge set (all solid lines), and $\Delta_e {\bf A}$ is a subset of $\Delta {\bf A}$ (the edges
are not adjacent to $\otimes$).}
\end{figure}
By the definition, 
$$\Delta_e {\bf A}\subset \Delta {\bf A}.$$
In addition, we denote the {\em interior vertex set} of ${\bf A}$ by 
$$\mbox{int}({\bf A})= {\bf A}\setminus \partial_i {\bf A}.$$ 
With these definitions, the following lemma is well known (see Lemma 2.23 Kesten (1986)).\\

{\bf Lemma 1} (Kesten). {\em If ${\bf A}$ is a ${\bf Z}^d$-connected finite set, then $\partial_e {\bf A}$ is
 a ${\bf Z}^d$-connected graph. }\\

Let 
$${\bf C}({\bf k},m)=\{v\in {\bf Z}^d: v\mbox{ is ${\bf Z}^d$-connected by an open path to }{\bf B}({\bf k},m)\}.$$
Note that ${\bf B}({\bf k},m)$ is a ${\bf Z}^d$-connected open set as we defined, so
$${\bf C}({\bf k},m) \mbox{ is a  ${\bf Z}^d$-connected open cluster,}\mbox{ and } {\bf B}({\bf k}, m)\subset {\bf C}({\bf k}, m).\eqno{(2.0)}$$ 
As we have defined, if there is a cutset that cuts  ${\bf B}({\bf k},m)$
from $\infty$, then any path from ${\bf B}({\bf k},m)$ to $\infty$ must use at least an edge of the cutset. 
Furthermore,
if there is a zero-cutset, then a path from 
${\bf B}({\bf k},m)$ to $\infty$ must not only use at least an edge of the cutset, but also a
zero-edge of the cutset.
If there is a zero-cutset,  then there is no open path from  ${\bf B}({\bf k}, m)$
to $\infty$, so
$$|{\bf C}({\bf k},m)| < \infty.\eqno{(2.1)}$$
On the other hand, if (2.1) holds, then there exists such a zero-cutset.
Let  ${\cal G}({\bf k},m)$ be the event that (2.1) occurs.
In this section, we will always discuss  particular fixed  configurations
in ${\cal G}({\bf k},m)$.

For each configuration in ${\cal G}({\bf k}, m)$, 
it follows from the definitions that the boundary edges of $\Delta {\bf C}({\bf k},m)$  are all closed
and they cut ${\bf B}({\bf k}, m)$ from $\infty$. However, $\Delta {\bf C}({\bf k},m)$ may contain too many
extra edges (see Fig. 1), so we would like to focus on $\Delta_e {\bf C}({\bf k},m)$. 
Since ${\bf C}({\bf k},m)$ is uniquely determined for each configuration in ${\cal G}({\bf k}, m)$, 
$\partial_e { \bf C}({\bf k},m)$ (see Fig. 1) is also uniquely determined. 
With these definitions, we have the following lemma.\\

{\bf Lemma 2.} {\em For all configurations in ${\cal G}({\bf k},m)$,
 $\Delta_e {\bf C}({\bf k},m)$ is a finite 
zero-cutset cutting ${\bf B}({\bf k},m)$ from $\infty$.}\\

{\bf Proof.} As we mentioned above, $\Delta_e {\bf C}({\bf k},m)$ is a zero-edge set.
 Since each vertex of $\partial_e {\bf C}({\bf k}, m)$  is ${\bf L}^d$-connected to 
$ {\bf C}({\bf k},m)$, $\Delta_e {\bf C}({\bf k},m)$ is finite. It remains to show that
$\Delta_e {\bf C}({\bf k},m)$ is a cutset. Since  ${\bf C}({\bf k},m)$ is finite,
for any ${\bf Z}^d$-path $\gamma$ from ${\bf B}({\bf k},m)$ to $\infty$, 
 some part of $\gamma$ must be outside of 
${\bf C}({\bf k},m)$. Let ${\bf u}$ be the last vertex in ${\bf C}({\bf k},m)$ such that after ${\bf u}$, 
the remaining piece
of $\gamma$ never uses another edge of ${\bf C}({\bf k},m)$. Suppose that after ${\bf u}$, $\gamma$ uses the edge $e$.
Thus, $e$ will be a zero-edge, otherwise, $e\in {\bf C}({\bf k},m)$. By the definition, 
$e\in \Delta {\bf C}({\bf k},m)$. On the other hand, the remaining  piece of $\gamma$ from $e$
will not return to ${\bf C}({\bf k},m)$ again as we defined. Thus, $e \in \Delta_e {\bf C}({\bf k},m)$.
Since any path must use an edge of $\Delta_e{\bf C}({\bf k},m)$, $\Delta_e{\bf C}({\bf k},m)$
will be a  cutset cutting ${\bf B}({\bf k},m)$ from $\infty$, so Lemma 2 follows. $\Box$\\

Furthermore, by Lemma 1 (see Fig. 1), we  know that
$$\partial_e {\bf C}({\bf k},m) \mbox{ is ${\bf Z}^d$-connected.}\eqno{(2.2)}$$

By Lemma 2, we know that $\Delta_e {\bf C}({\bf k},m)$ is a zero-cutset. However, we cannot use this cutset
to show Theorem 1, since this cutset might be too tangled. 
We need to eliminate
the tangled parts of  $\Delta_e {\bf C}({\bf k},m)$
to construct another zero-cutset.
To construct such a zero-cutset, we  use the idea of renormalization in  Kesten and Zhang (1990).
We define, for integer $t\geq 1$ and ${\bf u}=(u_1,\cdots, u_d)\in {\bf Z}^d$, the cube
$$B_t({\bf u})=\prod_{i=1}^d[tu_i, tu_i+t].$$
Here we need to take $t$ large, but much smaller than $m$ and $k_1, \cdots, k_{d-1}$.
Also, without loss of generality, we assume that
$k_i/t$ for $i=1,2,...,d-1$ and $m/t$ are integers; otherwise we may use $\lfloor k_i/t\rfloor$ and $\lfloor m/t\rfloor$. 
Usually, we consider  the ${\bf Z}^d$-vertices 
 in $B_t({\bf u})$. In addition, we can also consider the edges in $B_t({\bf u})$ if their two vertices
belong to $B_t({\bf u})$. For a cube $B_t({\bf u})$, we denote by $\bar {B}_t({\bf u})$ the cube $B_t({\bf u})$ and
its ${\bf L}^d$-adjacent neighbor cubes. 
 We call $B_t({\bf u})$  a {\em $t$-cube} and call $\bar {B}_t({\bf u})$  a {\em $3t$-cube}.
Through a simple computation, $\bar{\bf B}({\bf k}, m)$ contains 9 or 27 $t$-cubes when $d=2$ or 3. In general,
$$\mbox{ the number of $t$-cubes in }\bar{B}_t({\bf u})\mbox{}\leq 2^{2d}.\eqno{(2.3)}$$
If ${\bf A}=\cup_i B_t({\bf u}_i)$, then we denote
$\bar{\bf A}=\{ \bar{B}_t({\bf u}_i)\}$.
Also, a ${\bf Z}^d$-connected neighbor of $B_t({\bf u})$ and $B_t({\bf v})$ have  common vertices. We simply name these vertices   the {\em surface } of $B_t({\bf u})$. A cube has $2d$ surfaces.
In particular, two surfaces of $B_t({\bf u})$ with a distance $t$ 
are called {\em opposite surfaces}.

Two cubes, $B_t({\bf u})$ and $B_t({\bf v})$, for ${\bf u}\neq {\bf v}$, 
are said to be ${\bf Z}^d$-{\em adjacent cubically} or ${\bf L}^d$-{\em adjacent cubically} if ${\bf u}$ and ${\bf v}$ are ${\bf Z}^d$- or
${\bf L}^d$-adjacent. If $({\bf v}_0, e_1, {\bf v}_1,...,{\bf v}_{n-1}, e_n, {\bf v}_n)$ is a ${\bf Z}^d$-path, then
$B_t({\bf v}_0), B_t({\bf v}_1),...,B_t({\bf v}_n)$ is a {\em cubic} ${\bf Z}^d$-{\em path}.
With these cubic paths, we can define a {\em cubic} ${\bf Z}^d$-{\em cluster}.
Similarly, if we replace the ${\bf Z}^d$-path by a ${\bf L}^d$-path, we can define a {\em cubic} ${\bf L}^d$-{\em cluster}.
Let  
$$ {\bf C}_t({\bf k},m)=\{B_t({\bf u}): B_t({\bf u})\cap {\bf C}({\bf k},m)\neq \emptyset\}.$$
Note that ${\bf C}({\bf k},m)$ is ${\bf Z}^d$-connected and that our cubes contain their inside boundaries,
so ${\bf C}_t({\bf k},m)$ is ${\bf L}^d$-connected in the sense of the connection of cubes.
A boundary cube of ${\bf C}_t({\bf k},m)$ is also defined as (see Fig. 2)
\begin{eqnarray*}
\partial_t {\bf C}({\bf k},m)=&&\{B_t({\bf u}):  B_t({\bf u})\mbox{ contains a vertex of }\partial {\bf C}({\bf k},m)\cup
\partial_i {\bf C}({\bf k},m)\}. 
\end{eqnarray*}
Note that ${\bf C}_t({\bf k}, m)$ and $\partial_{t} {\bf C}({\bf k}, m)$ may
have a common cube. 
For convenience, we account $B_t({\bf u})$ in $\partial_t {\bf C}({\bf k}, m)$ if
$B_t({\bf u})$ is a common cube in both ${\bf C}_t({\bf k}, m)$ and $\partial_t {\bf C}({\bf k}, m)$. Since all the edges in ${\bf B}({\bf k}, m)$ are open and $k_i/t$ and $m/t$
are integers for $i=1,\cdots, d-1$, $\partial_t {\bf C}({\bf k}, m)$ and ${\bf B}({\bf k}, m)$ have no common edges, but possible common vertices. We account the cubes in ${\bf B}({\bf k}, m)$ if  they have common vertices.
Note also that $\partial_e{\bf C}({\bf k},m)$ contains all boundary edges of ${\bf C}({\bf k},m)$.

As we proved, on ${\cal G}({\bf k},m)$, $\Delta_e {\bf C}({\bf k},m)$  is a zero-cutset.
Therefore,  
$$\mbox{ there is a zero-cutset in } \partial_t {\bf C}({\bf k},m).\eqno{(2.4)}$$
In addition, the {\em exterior  cube-boundary} of ${\bf C}_t({\bf k},m)$ is defined as (see Fig. 2)
\begin{eqnarray*}
\partial_{t, e}{\bf C}({\bf k},m)=&&\{B_t({\bf u})\in  \partial_t {\bf C}({\bf k},m): B_t({\bf u}) \mbox{ is connected by}\\
&& \mbox{a cubic ${\bf Z}^d$-path  to $\infty$ outside $ \partial_t {\bf C}({\bf k},m)$}\}.
\end{eqnarray*}

Note that each cube in $\partial_{t, e}{\bf C}({\bf k},m)$ contains at least one vertex of  $\partial_e{\bf C}({\bf k},m)$.
Thus, by (2.2),
$$\mbox{ $\partial_{t,e} {\bf C}({\bf k},m)$ is ${\bf L}^d$-connected.}\eqno{(2.5)}$$
If we consider $\partial_{t, e} {\bf C}({\bf k},m)$ as an edge set, it  follows from the same proof of Lemma 2 to show that
$$\partial_{t, e} {\bf C}({\bf k},m)\mbox{  is a finite cutset that cuts
${\bf B}({\bf k},m)$ from $\infty$}.\eqno{(2.6)}$$
By (2.4), we know that $\partial_{t, e} {\bf C}({\bf k},m)$ is a cutset. On the other hand,
as the main task, Kesten and Zhang (1990) showed that $\partial_{t, e} {\bf C}({\bf k},m)$ is not very tangled.
However, unlike $\partial_e {\bf C}({\bf k},m)$, $\partial_{t, e} {\bf C}({\bf k},m)$ may not contain a zero-cutset  (see Fig. 2). 
The main task in this section is to combine $\partial_{t, e} {\bf C}({\bf k},m)$ with additional  edges
to construct a zero-cutset.  This construction is much easier to understand through 
Fig. 2, than through   rigorous written descriptions. We suggest that readers refer to Fig. 2 
while reading the following definitions.

\begin{figure}\label{F:alphabeta}
\begin{center}
\setlength{\unitlength}{0.0125in}%
\begin{picture}(200,130)(67,740)
\thicklines
\put(100, 670){\line(0,1){50}}
\put(100, 670){\line (1,0){70}}
\put(100, 720){\line (1,0){70}}
\put(170, 670){\line (0,1){50}}
\put(-20, 900){\dashbox(20,20)[br]}
\put(0, 900){\dashbox(20,20)[br]}
\put(20, 900){\dashbox(20,20)[br]}
\put(40, 900){\dashbox(20,20)[br]}
\put(60, 900){\dashbox(20,20)[br]}
\put(80, 900){\dashbox(20,20)[br]}
\put(100, 900){\dashbox(20,20)[br]}
\put(120, 900){\dashbox(20,20)[br]}
\put(140, 900){\dashbox(20,20)[br]}
\put(160, 900){\dashbox(20,20)[br]}

\put(70, 895){\line(0,1){24}}
\put(71, 895){\line(0,1){24}}
\put(48, 925){$\mbox{ not accounted}$}
\put(70, 920){\circle{5}}


\put(60, 880){\dashbox(20,20)[br]}

\put(-25,910){\circle*{4}}

\put(-20,910){\circle*{4}}
\put(-15,910){\circle*{4}}
\put(-10,910){\circle*{4}}
\put(-5,910){\circle*{4}}
\put(0,910){\circle*{4}}
\put(5,910){\circle*{4}}
\put(10,910){\circle*{4}}
\put(15,910){\circle*{4}}
\put(20,910){\circle*{4}}
\put(25,910){\circle*{4}}
\put(30,910){\circle*{4}}
\put(35,910){\circle*{4}}
\put(40,910){\circle*{4}}
\put(45,910){\circle*{4}}
\put(50,910){\circle*{4}}
\put(55,910){\circle*{4}}
\put(60,910){\circle*{4}}

\put(63,905){\circle*{4}}
\put(63,900){\circle*{4}}
\put(63,895){\circle*{4}}
\put(68,895){\circle*{4}}
\put(72,895){\circle*{4}}
\put(76,895){\circle*{4}}






\put(78,905){\circle*{4}}
\put(78,900){\circle*{4}}
\put(78,895){\circle*{4}}


\put(80,910){\circle*{4}}
\put(85,910){\circle*{4}}
\put(90,910){\circle*{4}}
\put(95,910){\circle*{4}}
\put(100,910){\circle*{4}}
\put(105,910){\circle*{4}}
\put(110,910){\circle*{4}}
\put(115,910){\circle*{4}}
\put(120,910){\circle*{4}}
\put(125,910){\circle*{4}}
\put(130,910){\circle*{4}}
\put(135,910){\circle*{4}}
\put(140,910){\circle*{4}}
\put(145,910){\circle*{4}}
\put(150,910){\circle*{4}}
\put(155,910){\circle*{4}}
\put(160,910){\circle*{4}}
\put(165,910){\circle*{4}}
\put(170,910){\circle*{4}}

\put(170,905){\circle*{4}}
\put(170,900){\circle*{4}}
\put(170,895){\circle*{4}}
\put(170,890){\circle*{4}}
\put(170,885){\circle*{4}}

\put(175,885){\circle*{4}}
\put(180,885){\circle*{4}}
\put(185,885){\circle*{4}}
\put(190,885){\circle*{4}}
\put(195,885){\circle*{4}}
\put(200,885){\circle*{4}}
\put(205,885){\circle*{4}}
\put(210,885){\circle*{4}}
\put(215,885){\circle*{4}}
\put(220,885){\circle*{4}}

\put(223,880){\circle*{4}}
\put(223,875){\circle*{4}}
\put(223,870){\circle*{4}}
\put(223,865){\circle*{4}}
\put(223,860){\circle*{4}}
\put(223,855){\circle*{4}}
\put(223,850){\circle*{4}}
\put(223,845){\circle*{4}}
\put(223,840){\circle*{4}}
\put(223,835){\circle*{4}}

\put(220,835){\circle*{4}}
\put(215,835){\circle*{4}}
\put(210,835){\circle*{4}}
\put(205,835){\circle*{4}}
\put(200,835){\circle*{4}}
\put(195,835){\circle*{4}}
\put(190,835){\circle*{4}}
\put(185,835){\circle*{4}}
\put(180,835){\circle*{4}}
\put(180,830){\circle*{4}}
\put(180,825){\circle*{4}}
\put(180,820){\circle*{4}}
\put(180,815){\circle*{4}}
\put(185,815){\circle*{4}}
\put(190,815){\circle*{4}}
\put(195,815){\circle*{4}}
\put(200,815){\circle*{4}}
\put(205,815){\circle*{4}}
\put(210,815){\circle*{4}}
\put(215,815){\circle*{4}}
\put(220,815){\circle*{4}}
\put(225,815){\circle*{4}}
\put(230,815){\circle*{4}}
\put(235,815){\circle*{4}}
\put(240,815){\circle*{4}}
\put(245,815){\circle*{4}}
\put(250,815){\circle*{4}}
\put(255,815){\circle*{4}}
\put(255,810){\circle*{4}}
\put(255,805){\circle*{4}}
\put(255,800){\circle*{4}}
\put(255,795){\circle*{4}}
\put(255,790){\circle*{4}}
\put(250,790){\circle*{4}}
\put(245,790){\circle*{4}}
\put(240,790){\circle*{4}}
\put(235,790){\circle*{4}}
\put(230,790){\circle*{4}}
\put(225,790){\circle*{4}}
\put(220,790){\circle*{4}}
\put(215,790){\circle*{4}}
\put(215,785){\circle*{4}}
\put(215,780){\circle*{4}}
\put(220,780){\circle*{4}}
\put(255,780){\circle*{4}}
\put(250,780){\circle*{4}}
\put(245,780){\circle*{4}}
\put(240,780){\circle*{4}}
\put(235,780){\circle*{4}}
\put(230,780){\circle*{4}}
\put(225,780){\circle*{4}}
\put(255,775){\circle*{4}}
\put(255,770){\circle*{4}}
\put(250,770){\circle*{4}}
\put(245,770){\circle*{4}}
\put(240,770){\circle*{4}}
\put(235,770){\circle*{4}}
\put(230,770){\circle*{4}}
\put(225,770){\circle*{4}}
\put(220,770){\circle*{4}}
\put(215,770){\circle*{4}}
\put(210,770){\circle*{4}}
\put(205,770){\circle*{4}}
\put(200,770){\circle*{4}}
\put(195,770){\circle*{4}}
\put(190,770){\circle*{4}}
\put(185,770){\circle*{4}}
\put(180,770){\circle*{4}}
\put(175,770){\circle*{4}}
\put(170,770){\circle*{4}}
\put(165,770){\circle*{4}}
\put(160,770){\circle*{4}}
\put(155,770){\circle*{4}}
\put(150,770){\circle*{4}}
\put(145,770){\circle*{4}}
\put(140,770){\circle*{4}}
\put(140,775){\circle*{4}}
\put(140,780){\circle*{4}}
\put(140,785){\circle*{4}}
\put(140,790){\circle*{4}}
\put(140,795){\circle*{4}}
\put(140,800){\circle*{4}}
\put(140,805){\circle*{4}}
\put(145,805){\circle*{4}}
\put(150,805){\circle*{4}}
\put(155,805){\circle*{4}}
\put(160,805){\circle*{4}}
\put(165,805){\circle*{4}}
\put(165,810){\circle*{4}}
\put(165,815){\circle*{4}}
\put(165,820){\circle*{4}}
\put(165,825){\circle*{4}}
\put(165,830){\circle*{4}}
\put(160,830){\circle*{4}}
\put(155,830){\circle*{4}}
\put(150,830){\circle*{4}}

\put(145,830){\circle*{4}}
\put(140,830){\circle*{4}}
\put(135,830){\circle*{4}}
\put(130,830){\circle*{4}}
\put(125,830){\circle*{4}}
\put(120,830){\circle*{4}}

\put(115,830){\circle*{4}}
\put(110,830){\circle*{4}}
\put(105,830){\circle*{4}}
\put(100,830){\circle*{4}}

\put(95,830){\circle*{4}}
\put(90,830){\circle*{4}}
\put(90,825){\circle*{4}}
\put(90,820){\circle*{4}}
\put(90,815){\circle*{4}}
\put(90,810){\circle*{4}}
\put(90,805){\circle*{4}}
\put(90,800){\circle*{4}}
\put(90,795){\circle*{4}}
\put(90,790){\circle*{4}}
\put(90,785){\circle*{4}}
\put(90,780){\circle*{4}}
\put(90,775){\circle*{4}}
\put(90,770){\circle*{4}}
\put(90,765){\circle*{4}}
\put(90,760){\circle*{4}}
\put(90,755){\circle*{4}}
\put(95,755){\circle*{4}}
\put(100,755){\circle*{4}}
\put(105,755){\circle*{4}}
\put(110,755){\circle*{4}}
\put(115,755){\circle*{4}}
\put(120,755){\circle*{4}}
\put(125,755){\circle*{4}}
\put(130,755){\circle*{4}}
\put(135,755){\circle*{4}}
\put(140,755){\circle*{4}}
\put(145,755){\circle*{4}}
\put(150,755){\circle*{4}}
\put(155,755){\circle*{4}}
\put(160,755){\circle*{4}}
\put(165,755){\circle*{4}}

\put(170,755){\circle*{4}}
\put(175,755){\circle*{4}}
\put(180,755){\circle*{4}}
\put(185,755){\circle*{4}}
\put(190,755){\circle*{4}}
\put(195,755){\circle*{4}}
\put(200,755){\circle*{4}}
\put(205,755){\circle*{4}}
\put(210,755){\circle*{4}}
\put(215,755){\circle*{4}}
\put(220,755){\circle*{4}}
\put(225,755){\circle*{4}}
\put(230,755){\circle*{4}}
\put(235,755){\circle*{4}}
\put(240,755){\circle*{4}}
\put(245,755){\circle*{4}}
\put(250,755){\circle*{4}}
\put(255,755){\circle*{4}}
\put(260,755){\circle*{4}}
\put(265,755){\circle*{4}}
\put(270,755){\circle*{4}}
\put(275,755){\circle*{4}}
\put(280,755){\circle*{4}}
\put(285,755){\circle*{4}}
\put(290,755){\circle*{4}}
\put(295,755){\circle*{4}}
\put(300,755){\circle*{4}}
\put(305,755){\circle*{4}}

\put(305,760){\circle*{4}}
\put(305,765){\circle*{4}}
\put(305,770){\circle*{4}}
\put(305,775){\circle*{4}}
\put(305,780){\circle*{4}}
\put(305,785){\circle*{4}}
\put(305,790){\circle*{4}}
\put(305,795){\circle*{4}}
\put(305,800){\circle*{4}}
\put(305,805){\circle*{4}}
\put(305,810){\circle*{4}}
\put(305,815){\circle*{4}}
\put(305,820){\circle*{4}}
\put(305,825){\circle*{4}}
\put(305,830){\circle*{4}}
\put(305,835){\circle*{4}}
\put(305,840){\circle*{4}}
\put(305,845){\circle*{4}}
\put(305,850){\circle*{4}}
\put(305,855){\circle*{4}}

\put(300,855){\circle*{4}}
\put(295,855){\circle*{4}}
\put(290,855){\circle*{4}}
\put(285,855){\circle*{4}}
\put(280,855){\circle*{4}}
\put(275,855){\circle*{4}}
\put(275,855){\circle*{4}}
\put(270,855){\circle*{4}}
\put(265,855){\circle*{4}}
\put(260,855){\circle*{4}}
\put(255,855){\circle*{4}}
\put(250,855){\circle*{4}}
\put(245,855){\circle*{4}}

\put(238,880){\circle*{4}}
\put(238,875){\circle*{4}}
\put(238,870){\circle*{4}}
\put(238,865){\circle*{4}}
\put(238,860){\circle*{4}}
\put(238,855){\circle*{4}}

\put(240,885){\circle*{4}}
\put(245,885){\circle*{4}}
\put(250,885){\circle*{4}}
\put(255,885){\circle*{4}}
\put(260,885){\circle*{4}}
\put(265,885){\circle*{4}}
\put(270,885){\circle*{4}}
\put(275,885){\circle*{4}}
\put(280,885){\circle*{4}}
\put(285,885){\circle*{4}}
\put(290,885){\circle*{4}}
\put(295,885){\circle*{4}}
\put(300,885){\circle*{4}}
\put(305,885){\circle*{4}}
\put(310,885){\circle*{4}}
\put(315,885){\circle*{4}}
\put(320,885){\circle*{4}}
\put(325,885){\circle*{4}}
\put(330,885){\circle*{4}}
\put(335,885){\circle*{4}}
\put(340,885){\circle*{4}}
\put(345,885){\circle*{4}}
\put(350,885){\circle*{4}}

\put(160, 880){\dashbox(20,20)[br]}
\put(180, 880){\dashbox(20,20)[br]}
\put(220, 880){\dashbox(20,20)[br]}
\put(240, 880){\dashbox(20,20)[br]}
\put(260, 880){\dashbox(20,20)[br]}
\put(280, 880){\dashbox(20,20)[br]}
\put(300, 880){\dashbox(20,20)[br]}
\put(320, 880){\dashbox(20,20)[br]}

\put(210, 910){${\bf D}(1)$}
\put(230, 900){\line(0,-1){80}}
\put(231, 900){\line(0,-1){80}}
\put(230, 846){\line(1,0){45}}
\put(230, 847){\line(1,0){45}}
\put(275, 846){\line(0,-1){84}}
\put(274, 846){\line(0,-1){84}}
\put(275, 810){\line(1,0){30}}
\put(275, 811){\line(1,0){30}}
\put(274, 763){\line(-1,0){160}}
\put(274, 764){\line(-1,0){160}}
\put(114, 763){\line(0,1){58}}
\put(115, 763){\line(0,1){58}}
\put(115, 817){\line(1,0){38}}
\put(115, 816){\line(1,0){38}}

\put(200, 700){\footnotesize $\partial_t {\bf C}({\bf k},m)$}
\put(220, 710){\vector (0,1){30}}
\put(255, 705){\vector (1,0){25}}
\put(220, 695){\vector (0,-1){30}}

\put(200, 880){\dashbox(20,20)[br]}
\put(220, 860){\dashbox(20,20)[br]}
\put(220, 840){\dashbox(20,20)[br]}
\put(220, 820){\dashbox(20,20)[br]}

\put(200, 820){\dashbox(20,20)[br]}
\put(180, 820){\dashbox(20,20)[br]}
\put(160, 820){\dashbox(20,20)[br]}

\put(160, 800){\dashbox(20,20)[br]}
\put(120, 800){\dashbox(20,20)[br]}
\put(120, 780){\dashbox(20,20)[br]}
\put(180, 760){\dashbox(20,20)[br]}
\put(200, 760){\dashbox(20,20)[br]}
\put(220, 760){\dashbox(20,20)[br]}
\put(240, 760){\dashbox(20,20)[br]}
\put(240, 780){\dashbox(20,20)[br]}
\put(240, 800){\dashbox(20,20)[br]}
\put(200, 780){\dashbox(20,20)[br]}

\put(143, 785){{\footnotesize$\partial_{t, e}{\bf D}(1, 1,1)$}}
\put(203, 790){\vector (1,0){10}}
\put(0, 815){{$\partial_{t, e} {\bf D}(1,1)$}}
\put(173, 840){\circle{5}}
\put(30, 800){\vector (4,-1){100}}
\put(110, 865){{$\partial_{t, e} {\bf D}(1)$}}
\put(230, 900){\circle{5}}
\put(130, 855){\vector (0,-1){15}}

\put(173, 857){\line (0,-1){63}}
\put(174, 857){\line (0,-1){63}}
\put(174, 805){\line (1,0){70}}
\put(174, 806){\line (1,0){70}}
\put(275, 785){\line (-1,0){57}}
\put(275, 786){\line (-1,0){57}}
\put(260, 786){\circle{5}}
\put(173, 857){\line (-1,0){110}}
\put(173, 858){\line (-1,0){110}}
\put(63, 858){\line (0,-1){170}}
\put(62, 858){\line (0,-1){170}}
\put(62, 688){\line (1,0){37}}
\put(62, 687){\line (1,0){37}}
\put(400, 810){\vector (-1,0){70}}
\put(370, 800){{\footnotesize $\mbox{\bf (block property)}$}}
\put(400, 810){\vector (-1,1){23}}

\put(360, 820){\dashbox(20,20)[br]}
\put(320, 800){\dashbox(20,20)[br]}

\put(180, 740){\dashbox(20,20)[br]}
\put(200, 740){\dashbox(20,20)[br]}
\put(220, 740){\dashbox(20,20)[br]}
\put(240, 740){\dashbox(20,20)[br]}
\put(260, 740){\dashbox(20,20)[br]}
\put(300, 740){\dashbox(20,20)[br]}
\put(280, 740){\dashbox(20,20)[br]}

\put(300, 760){\dashbox(20,20)[br]}
\put(300, 780){\dashbox(20,20)[br]}
\put(300, 800){\dashbox(20,20)[br]}
\put(300, 820){\dashbox(20,20)[br]}
\put(300, 840){\dashbox(20,20)[br]}

\put(280, 840){\dashbox(20,20)[br]}
\put(260, 840){\dashbox(20,20)[br]}
\put(240, 840){\dashbox(20,20)[br]}

\put(140, 820){\dashbox(20,20)[br]}
\put(120, 820){\dashbox(20,20)[br]}
\put(100, 820){\dashbox(20,20)[br]}
\put(80, 820){\dashbox(20,20)[br]}
\put(80, 800){\dashbox(20,20)[br]}
\put(80, 780){\dashbox(20,20)[br]}
\put(80, 760){\dashbox(20,20)[br]}
\put(80, 740){\dashbox(20,20)[br]}
\put(100, 740){\dashbox(20,20)[br]}
\put(120, 740){\dashbox(20,20)[br]}
\put(140, 740){\dashbox(20,20)[br]}
\put(160, 740){\dashbox(20,20)[br]}
\put(160, 800){\dashbox(20,20)[br]}
\put(140, 800){\dashbox(20,20)[br]}
\put(120, 760){\dashbox(20,20)[br]}
\put(140, 760){\dashbox(20,20)[br]}
\put(160, 760){\dashbox(20,20)[br]}
\put(220, 780){\dashbox(20,20)[br]}



\put(120,700){\mbox{${\bf B}({\bf k}, m)$}}

\put(-25,910){\circle*{4}}
\put(-25,905){\circle*{4}}
\put(-25,900){\circle*{4}}
\put(-25,895){\circle*{4}}
\put(-25,890){\circle*{4}}
\put(-25,885){\circle*{4}}
\put(-25,880){\circle*{4}}
\put(-25,875){\circle*{4}}
\put(-25,870){\circle*{4}}
\put(-25,865){\circle*{4}}
\put(-25,860){\circle*{4}}
\put(-25,855){\circle*{4}}
\put(-25,850){\circle*{4}}
\put(-25,845){\circle*{4}}
\put(-25,840){\circle*{4}}
\put(-25,835){\circle*{4}}
\put(-25,830){\circle*{4}}
\put(-25,825){\circle*{4}}
\put(-25,820){\circle*{4}}
\put(-25,815){\circle*{4}}
\put(-25,810){\circle*{4}}
\put(-25,805){\circle*{4}}
\put(-25,800){\circle*{4}}
\put(-25,795){\circle*{4}}
\put(-25,790){\circle*{4}}
\put(-25,785){\circle*{4}}
\put(-25,780){\circle*{4}}
\put(-25,775){\circle*{4}}
\put(-25,770){\circle*{4}}
\put(-25,765){\circle*{4}}
\put(-25,760){\circle*{4}}
\put(-25,755){\circle*{4}}

\put(-20,755){\circle*{4}}
\put(-15,755){\circle*{4}}
\put(-10,755){\circle*{4}}
\put(-05,755){\circle*{4}}
\put(-05,750){\circle*{4}}
\put(-05,745){\circle*{4}}
\put(-05,740){\circle*{4}}
\put(-05,735){\circle*{4}}
\put(-05,730){\circle*{4}}
\put(-05,725){\circle*{4}}
\put(-05,720){\circle*{4}}
\put(-05,715){\circle*{4}}
\put(-05,710){\circle*{4}}
\put(-05,705){\circle*{4}}
\put(-05,700){\circle*{4}}
\put(-05,695){\circle*{4}}
\put(-5,690){\circle*{4}}
\put(-5,685){\circle*{4}}
\put(-5,680){\circle*{4}}
\put(-5,675){\circle*{4}}
\put(-5,670){\circle*{4}}
\put(-5,665){\circle*{4}}
\put(-5,660){\circle*{4}}
\put(-5,655){\circle*{4}}
\put(-5,650){\circle*{4}}

\put(0,650){\circle*{4}}
\put(5,650){\circle*{4}}

\put(10,650){\circle*{4}}
\put(15,650){\circle*{4}}
\put(20,650){\circle*{4}}
\put(25,650){\circle*{4}}
\put(30,650){\circle*{4}}
\put(35,650){\circle*{4}}
\put(40,650){\circle*{4}}
\put(45,650){\circle*{4}}
\put(50,650){\circle*{4}}
\put(55,650){\circle*{4}}
\put(60,650){\circle*{4}}
\put(65,650){\circle*{4}}
\put(70,650){\circle*{4}}
\put(75,650){\circle*{4}}
\put(80,650){\circle*{4}}
\put(85,650){\circle*{4}}
\put(90,650){\circle*{4}}
\put(95,650){\circle*{4}}
\put(100,650){\circle*{4}}
\put(105,650){\circle*{4}}
\put(110,650){\circle*{4}}
\put(115,650){\circle*{4}}
\put(120,650){\circle*{4}}
\put(125,650){\circle*{4}}
\put(130,650){\circle*{4}}
\put(135,650){\circle*{4}}
\put(140,650){\circle*{4}}
\put(145,650){\circle*{4}}
\put(150,650){\circle*{4}}
\put(155,650){\circle*{4}}
\put(160,650){\circle*{4}}
\put(165,650){\circle*{4}}
\put(170,650){\circle*{4}}
\put(175,650){\circle*{4}}
\put(180,650){\circle*{4}}
\put(185,650){\circle*{4}}
\put(190,650){\circle*{4}}
\put(195,650){\circle*{4}}
\put(200,650){\circle*{4}}
\put(205,650){\circle*{4}}
\put(210,650){\circle*{4}}
\put(215,650){\circle*{4}}
\put(220,650){\circle*{4}}
\put(225,650){\circle*{4}}
\put(230,650){\circle*{4}}
\put(235,650){\circle*{4}}
\put(240,650){\circle*{4}}
\put(265,650){\circle*{4}}
\put(270,650){\circle*{4}}

\put(245,650){\circle*{4}}
\put(245,650){\circle*{4}}
\put(245,650){\circle*{4}}
\put(250,650){\circle*{4}}
\put(255,650){\circle*{4}}
\put(260,650){\circle*{4}}
\put(275,650){\circle*{4}}
\put(280,650){\circle*{4}}
\put(285,650){\circle*{4}}
\put(290,650){\circle*{4}}
\put(295,650){\circle*{4}}
\put(300,650){\circle*{4}}
\put(305,650){\circle*{4}}
\put(310,650){\circle*{4}}
\put(315,650){\circle*{4}}
\put(320,650){\circle*{4}}
\put(325,650){\circle*{4}}
\put(330,650){\circle*{4}}
\put(335,650){\circle*{4}}
\put(340,650){\circle*{4}}
\put(345,650){\circle*{4}}
\put(350,650){\circle*{4}}
\put(350,655){\circle*{4}}
\put(350,665){\circle*{4}}

\put(350,665){\circle*{4}}
\put(350,670){\circle*{4}}
\put(350,675){\circle*{4}}
\put(350,660){\circle*{4}}
\put(350,665){\circle*{4}}
\put(350,670){\circle*{4}}
\put(350,675){\circle*{4}}
\put(350,680){\circle*{4}}
\put(350,685){\circle*{4}}
\put(350,690){\circle*{4}}
\put(350,695){\circle*{4}}
\put(350,700){\circle*{4}}

\put(350,705){\circle*{4}}
\put(345,705){\circle*{4}}
\put(340,705){\circle*{4}}
\put(335,705){\circle*{4}}
\put(330,705){\circle*{4}}
\put(325,705){\circle*{4}}
\put(320,705){\circle*{4}}
\put(315,705){\circle*{4}}
\put(310,705){\circle*{4}}
\put(305,705){\circle*{4}}
\put(300,705){\circle*{4}}
\put(295,705){\circle*{4}}
\put(295, 711){\line(1,0){65}}
\put(280,670){\mbox{$\partial_{t, e}{\bf D}(2)$}}
\put(290, 680){\vector (0,1){20}}
\put(370,710){\mbox{${\bf D}(2)$}}
\put(360, 710){\circle{5}}
\put(-110,900){\mbox{$\partial_{t, e}{\bf C}({\bf k},m)$}}
\put(-80, 890){\vector (3,-1){40}}
\put(-110,850){\mbox{$\partial_e{\bf C}({\bf k},m)$}}
\put(-90, 840){\vector (2,-1){60}}

\put(295, 710){\line(1,0){65}}

\put(350,715){\circle*{4}}
\put(345,715){\circle*{4}}
\put(340,715){\circle*{4}}
\put(335,715){\circle*{4}}
\put(330,715){\circle*{4}}
\put(325,715){\circle*{4}}
\put(320,715){\circle*{4}}
\put(315,715){\circle*{4}}
\put(310,715){\circle*{4}}
\put(305,715){\circle*{4}}
\put(300,715){\circle*{4}}
\put(295,715){\circle*{4}}
\put(295,710){\circle*{4}}

\put(350,715){\circle*{4}}
\put(350,720){\circle*{4}}
\put(350,725){\circle*{4}}
\put(350,730){\circle*{4}}
\put(350,735){\circle*{4}}

\put(350,740){\circle*{4}}
\put(350,745){\circle*{4}}
\put(350,750){\circle*{4}}
\put(350,755){\circle*{4}}
\put(350,760){\circle*{4}}
\put(350,765){\circle*{4}}
\put(350,770){\circle*{4}}
\put(350,775){\circle*{4}}
\put(350,780){\circle*{4}}
\put(350,785){\circle*{4}}
\put(350,790){\circle*{4}}
\put(350,795){\circle*{4}}
\put(350,800){\circle*{4}}
\put(350,805){\circle*{4}}
\put(350,810){\circle*{4}}
\put(350,815){\circle*{4}}
\put(350,820){\circle*{4}}
\put(350,825){\circle*{4}}
\put(350,830){\circle*{4}}
\put(350,835){\circle*{4}}
\put(350,840){\circle*{4}}
\put(350,845){\circle*{4}}
\put(350,850){\circle*{4}}
\put(350,855){\circle*{4}}
\put(350,860){\circle*{4}}
\put(350,865){\circle*{4}}
\put(350,870){\circle*{4}}
\put(350,875){\circle*{4}}
\put(350,880){\circle*{4}}
\put(-40, 900){\dashbox(20,20)[br]}
\put(-40, 880){\dashbox(20,20)[br]}
\put(-40, 860){\dashbox(20,20)[br]}
\put(-40, 840){\dashbox(20,20)[br]}
\put(-40, 820){\dashbox(20,20)[br]}
\put(-40, 800){\dashbox(20,20)[br]}
\put(-40, 780){\dashbox(20,20)[br]}
\put(-40, 760){\dashbox(20,20)[br]}
\put(-40, 740){\dashbox(20,20)[br]}
\put(-20, 740){\dashbox(20,20)[br]}
\put(-20, 720){\dashbox(20,20)[br]}
\put(-20, 700){\dashbox(20,20)[br]}
\put(-20, 680){\dashbox(20,20)[br]}
\put(-20, 660){\dashbox(20,20)[br]}
\put(-20, 640){\dashbox(20,20)[br]}
\put(0, 640){\dashbox(20,20)[br]}
\put(20, 640){\dashbox(20,20)[br]}
\put(40, 640){\dashbox(20,20)[br]}
\put(-20, 640){\dashbox(20,20)[br]}
\put(0, 640){\dashbox(20,20)[br]}
\put(20, 640){\dashbox(20,20)[br]}
\put(40, 640){\dashbox(20,20)[br]}
\put(60, 640){\dashbox(20,20)[br]}
\put(80, 640){\dashbox(20,20)[br]}
\put(100, 640){\dashbox(20,20)[br]}
\put(120, 640){\dashbox(20,20)[br]}
\put(140, 640){\dashbox(20,20)[br]}
\put(160, 640){\dashbox(20,20)[br]}
\put(180, 640){\dashbox(20,20)[br]}
\put(200, 640){\dashbox(20,20)[br]}
\put(220, 640){\dashbox(20,20)[br]}
\put(240, 640){\dashbox(20,20)[br]}
\put(260, 640){\dashbox(20,20)[br]}
\put(280, 640){\dashbox(20,20)[br]}
\put(300, 640){\dashbox(20,20)[br]}
\put(320, 640){\dashbox(20,20)[br]}
\put(340, 640){\dashbox(20,20)[br]}
\put(340, 660){\dashbox(20,20)[br]}
\put(340, 680){\dashbox(20,20)[br]}
\put(340, 700){\dashbox(20,20)[br]}
\put(340, 720){\dashbox(20,20)[br]}
\put(340, 740){\dashbox(20,20)[br]}

\put(340, 760){\dashbox(20,20)[br]}
\put(340, 780){\dashbox(20,20)[br]}
\put(340, 800){\dashbox(20,20)[br]}
\put(340, 820){\dashbox(20,20)[br]}
\put(340, 840){\dashbox(20,20)[br]}
\put(340, 860){\dashbox(20,20)[br]}
\put(340, 880){\dashbox(20,20)[br]}

\put(320, 700){\dashbox(20,20)[br]}
\put(300, 700){\dashbox(20,20)[br]}
\put(280, 700){\dashbox(20,20)[br]}

\put(180, 800){\dashbox(20,20)[br]}
\put(200, 800){\dashbox(20,20)[br]}
\put(220, 800){\dashbox(20,20)[br]}

\end{picture}
\end{center}
\vskip 1in
\caption{ \em This graph shows how to construct  linear cubes that contain a zero-, or closed, cutset.
The large dotted line, $\partial_e {\bf C}({\bf k}, m)$,  is a zero-boundary for open cluster  ${\bf B}({\bf k}, m)$. 
The cubes 
$\{B_t({\bf u})\}$ that contain $\partial_e {\bf C}({\bf k}, m)$ are $\partial_t {\bf C}({\bf k}, m)$. 
Part of the cubes from $\partial_t {\bf C}({\bf k}, m)$ enclose a cubic circuit $\partial_{t,e} {\bf C}({\bf k}, m)$.
The surface of the circuit, ${\bf U}(\partial_{t,e} {\bf C}({\bf k}, m)),$ divides ${\bf Z}^d$ into two parts:
the inside part including ${\bf B}({\bf k}, m)$ and the outside part.
On the  surface, there are three exits  such that  open clusters (boldfaced lines)  
 penetrate $\partial_{t,e} {\bf C}({\bf k}, m)$ from the outside. The first one will not be accounted
 since it stays in $\partial_{t, e}\bar{\bf C}({\bf k}, m)$. The other two are ${\bf D}(i_1)$ for $i_1=1,2$.  For  ${\bf D}(1)$,  its exterior-cube-boundary consists of $\partial_{t, e}{\bf D}(1)$.  There is another open cluster,  ${\bf D}(1, 1)$ from ${\bf U}(\partial_{t, e}{\bf D}(1))$ penetrating     $\partial_{t, e}{\bf D}(1)$.  Its exterior-cube-boundary consists of $\partial_{t, e}{\bf D}(1,1)$.  Furthermore, there is another open cluster ${\bf D}(1, 1,1)$ from ${\bf U}(\partial_{t, e}{\bf D}(1,1))$ penetrating     $\partial_{t, e}{\bf D}(1,1)$.  For  ${\bf D}(1,1,1)$,  its exterior-cube-boundary consists of $\partial_{t, e}{\bf D}(1,1,1)$.  Note that ${\bf D}(1,1)$ is connected to ${\bf B}({\bf k}, m)$, so it is a joint open cluster. The others are disjoint.
 ${\bf \Gamma}({\bf k}, t)$ consists of all the cubes in $\partial_{t,e} {\bf C}({\bf k}, m)$,  
$\partial_{t, e}{\bf D}(1)$,   $\partial_{t, e}{\bf D}(1,1)$, $\partial_{t, e}{\bf D}(1,1,1)$
and $\partial_{t, e}{\bf D}(2)$.
  For each cube  ${B}_t({\bf u})$ in   $\partial_{t, e}{\bf C}({\bf k}, m)$ or in $\partial_{t, e}{\bf D}(i_1, \cdots, i_j)$, there is a cube adjacent  to  ${B}_t'({\bf u})$ cubically with a block property indicted in the graph.}
\end{figure}
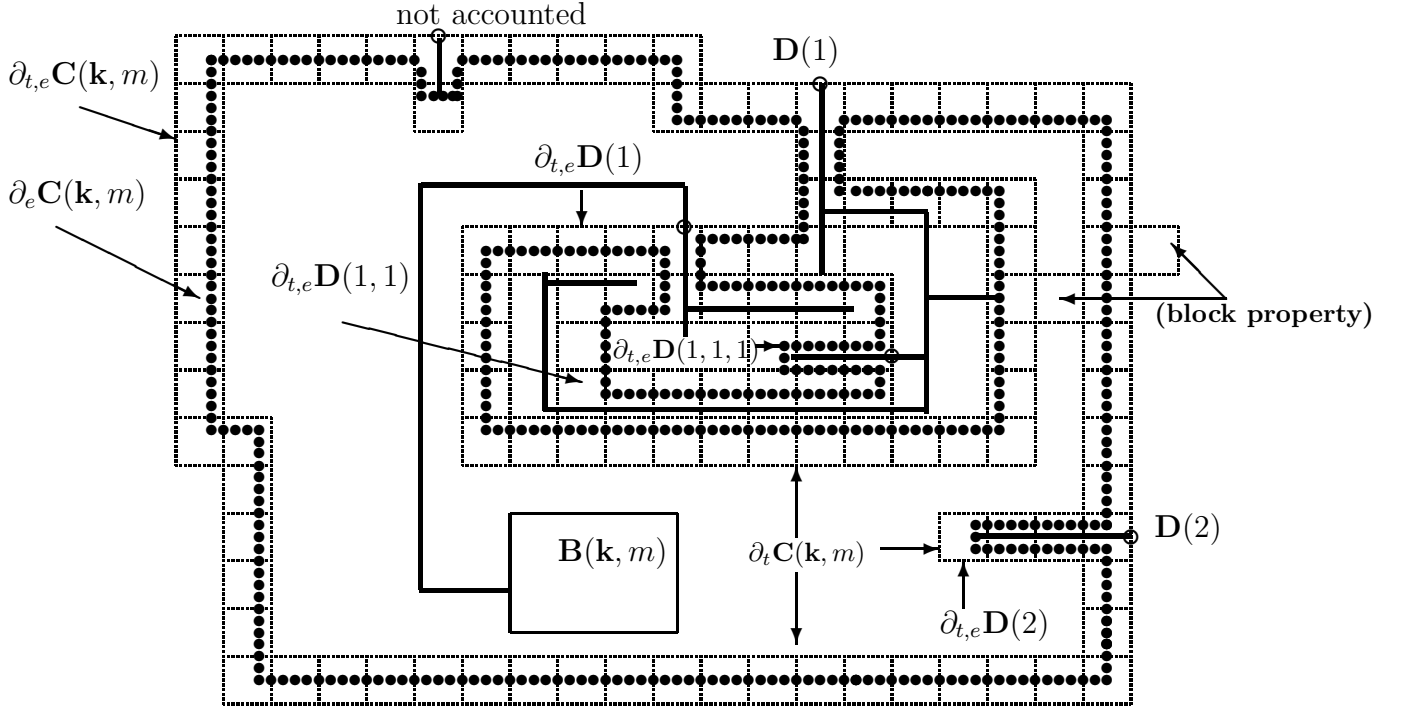

For configurations in ${\cal G}({\bf k},m)$, $\partial_{t, e} {\bf C}({\bf k},m)$ is a cutset.  Each path from 
$\infty$ to ${\bf B}({\bf k},m)$ must meet a vertex of $\partial_{t, e} {\bf C}({\bf k},m)$ and then go to
$\partial {\bf B}({\bf k}, m)$ from the vertex. 
We name these vertices
the {\em surface} of $\partial_{t, e} {\bf C}({\bf k},m)$ and denote them by ${\bf U}(\partial_{t, e} {\bf C}({\bf k},m))$. 
More precisely,
\begin{eqnarray*}
{\bf U}(\partial_{t, e} {\bf C}({\bf k},m))=&&\{{\bf v}\in \partial_{t, e} {\bf C}({\bf k},m): \exists 
\mbox{ a ${\bf Z}^d$-path from
${\bf v}$ to $\infty$ }\\
&&\mbox{ without using vertices, except ${\bf v}$, in $\partial_{t, e} {\bf C}({\bf k},m)$}\}.
\end{eqnarray*}
By (2.6) and the definition of the surface, 
$${\bf U}(\partial_{t, e} {\bf C}({\bf k},m))\cap {\bf C}({\bf k}, m)=\emptyset.\eqno{(2.7)}$$
Since $\partial_{t, e} {\bf C}({\bf k},m)$  cuts ${\bf B}({\bf k},m)$ from $\infty$, it divides 
the vertices of ${\bf Z}^d$ into two parts: the inside and outside parts, where the inside part, containing ${\bf B}({\bf k}, m)$,
is enclosed by ${\bf U}(\partial_{t, e} {\bf C}({\bf k},m))$. We denote the vertices in the inside part 
 by
${\bf L}_t({\bf k}, m)$.

We consider  the  open clusters in ${\bf L}_t({\bf k}, m)$ (see Fig. 2)  starting from
the  boundary of  ${\bf U}(\partial_{t, e} {\bf C}({\bf k}, m))$. Many open clusters might 
be inside $\partial_{t, e}\bar{\bf C}({\bf k}, m)$ and we ignore them,
where 
$$\partial_{t, e}\bar{\bf C}({\bf k}, m)=\{\bar{B}({\bf u}): B({\bf u})\in \partial_{t, e}{\bf C}({\bf k}, m)\}.$$
We denote all the remaining open clusters by
${\bf D}(1), \cdots, {\bf D}(i_1), \cdots,  {\bf D}(l_1)$.  It follows from our definition that
$${\bf D}(i)\cap {\bf D}(j)=\emptyset \mbox{ for } 1\leq i\neq j \leq l_1.\eqno{(2.8)}$$
We also denote by
$$\{{\bf u}_1(i_1),\cdots,  {\bf u}_{q_1}(i_1)\}={\bf D}(i_1)\cap {\bf U}(\partial_{t, e} {\bf C}({\bf k}, m))\mbox{ for each }i_1.$$
They are called {\em exits}.
We consider the exterior boundary of ${\bf D}(i_1)$ for $1\leq i_1\leq l_1$. We denote them by
$\partial_{e}  {\bf D}(i_1)$. With the exterior boundary,  let 
$$\partial_{t,e}  {\bf D}(i_1)= \{ B_t({\bf u}): B_t({\bf u})\cap \partial_{e}  {\bf D}(i_1)\neq \emptyset\}.$$
Note that it is possible that $\partial_{t,e}  {\bf D}(i_1)$ and $\partial_{t, e}\bar{\bf C}({\bf k}, m)$
have many common cubes.
 For each $B_t({\bf u})\in \partial_{t,e}  {\bf D}(i_1)$, but not in $\partial_{t, e}\bar{\bf C}({\bf k}, m)$,  it follows from the definition of the exterior boundary of 
  $\partial_{t, e}{\bf D}(i_1)$ that (see Fig. 2)
  \begin{eqnarray*}
  &&\exists \mbox{ an open path in ${\bf D}(i_1)$ from }B_t({\bf u})\mbox{ to } \partial \bar{B}_t({\bf u})\mbox{ and } 
  B_t({\bf u})  \mbox{ has  an ${\bf L}^d$-neighbor  $B_t'({\bf u})$ } \\
  && \mbox{ with } B_t'({\bf u})\cap {\bf D}(i_1) =\emptyset.\hskip 10cm {(2.9)}
  \end{eqnarray*}
  Similarly,  for $i_1$, we may consider the surface of $\partial_{t, e}{\bf D}(i_1)$ as
  \begin{eqnarray*}
{\bf U}(\partial_{t, e} {\bf D}(i_1))=&&\{{\bf v}\in \partial_{t, e} {\bf D}(i_1): \exists 
\mbox{ a ${\bf Z}^d$-path from
${\bf v}$ to $\infty$ }\\
&&\mbox{ without using vertices, except ${\bf v}$, in $\partial_{t, e} {\bf D}(i_1)$}\}.
\end{eqnarray*}
${\bf U}(\partial_{t, e} {\bf D}(i_1))$ also divides ${\bf Z}^d$ into two parts: inside and
outside parts. We continue  to find the open clusters  from  ${\bf U}(\partial_{t, e} {\bf D}(i_1))$ in the inside part.  Note that ${\bf D}(i_1)$ is one of them. Similarly, we ignore the open clusters   in $\partial_{t, e} \bar{\bf D}(i_1)$  and  ${\bf D}(i_1)$.
where
$$\partial_{t, e}\bar{\bf D}(i_1)=\{\bar{B}({\bf u}): B({\bf u})\in \partial_{t, e}{\bf D}(i_1)\}.$$
 We denote the remaining open clusters by ${\bf D}(i_1, i_2)$ for $1\leq i_2\leq l_2$ (see Fig. 2).
 We also denote the exits of these open clusters by
 $$\{{\bf u}_1( i_2),\cdots, u_{q_2}( i_2)\}= D(i_1, i_2)\cap {\bf U}(\partial_{t, e} {\bf D}(i_1)).$$
  Note that ${\bf L}_t({\bf k}, m)$ is a finite vertex set, so we can continue this way to find all the open clusters ${\bf D}(i_1, \cdots i_j)$ for $1\leq i_j\leq l_j$ and for $1\leq j\leq s$,
 their exterior cubic boundaries   $\partial_{t, e} {\bf D}(i_1, \cdots i_j)$, and their exits
 $\{{\bf u}_1(i_j),\cdots, u_{q_{j}}(i_j)\}$.
  Note that some of open clusters ${\bf D}(i_1, \cdots i_j)$ will never connect to
  ${\bf B}({\bf k}, m)$ (see Fig. 2), but some of them do. We call them {\em disjoint} open clusters or {\em joint} open clusters.  By Lemma 2, ${\bf D}(i_1)$ is always disjoint.
  Similarly to (2.8) and (2.9), we have for $x\neq y$,
$${\bf D}(i_1, \cdots i_{j}, x)\cap {\bf D}(i_1, \cdots i_{j}, y)\neq \emptyset\eqno{(2.10)}$$
and  for each $B_t({\bf u}) \in \partial_{t, e}{\bf D}(i_1, \cdots i_j)$ for $j\leq s$, 
 \begin{eqnarray*}
  &&\exists \mbox{ an open path in ${\bf D}(i_1, \cdots, i_j)$ from }B_t({\bf u})\mbox{ to } \partial \bar{B}_t({\bf u})\mbox{ and } 
  B_t({\bf u})  \mbox{ has  an ${\bf L}^d$-neighbor}   \\
  && B_t'({\bf u}) \mbox{ with } B_t'({\bf u})\cap {\bf D}(i_1, \cdots, i_j) =\emptyset.\hskip 8.1cm {(2.11)}
  \end{eqnarray*}
Recall that  $\partial_{e, t}{\bf D}(i_1, \cdots i_j)$  is the exterior cube-boundaries of ${\bf D}(i_1, \cdots i_j)$. 
We collect all these $t$-cubes  of the open clusters together with 
$\partial_{t, e}{\bf C}({\bf k}, m)$ to have the following cubic set
\begin{eqnarray*}
{\bf \Gamma}_t({\bf k}, m)&= &\partial_{t, e}{\bf C}({\bf k}, m)\bigcup_{i_1}\partial {\bf D}_{t, e} (i_1)\cdots \bigcup_{i_j} \partial_{t, e}{\bf D}(i_1, \cdots, i_j)\cdots\bigcup_{i_s} \partial_{t, e}{\bf D}(i_1, \cdots, i_s)\\
&=&\partial_{t, e}{\bf C}({\bf k}, m)\bigcup_{j=1}^s \bigcup_{i_j} \partial_{t, e}{\bf D}(i_1, \cdots, i_j)
\end{eqnarray*}
and
$$ \bar{ \bf{\Gamma}}_t({\bf k}, m)= \{\bar{B}({\bf u}): {B}({\bf u})\in  \bf {\Gamma}_t({\bf k}, m)\}.$$
Note that, by Lemma 1, $\partial_{t,e}  {\bf D}(i_1)$ is ${\bf L}^d$-connected. By the definition,   $\partial_{t,e}  {\bf D}(i_1)$ is also connected to
$\partial_{t, e}{\bf C}({\bf k}, m)$ cubically.  Therefore,
$$\partial_{t, e}{\bf C}({\bf k}, m)\bigcup_{i_1} \partial_{t,e}  {\bf D}(i_1)\mbox{ is ${\bf L}^d$-
connected cubically}.\eqno{(2.12)}$$
Similarly, $\partial_{t,e}  {\bf D}(i_1, i_2)$ is ${\bf L}^d$-connected  and also connected to $\partial_{t,e}  {\bf D}(i_1)$ cubically.  Thus,
$$\partial_{t, e}{\bf C}({\bf k}, m)\bigcup_{i_1} \partial_{t,e}  {\bf D}(i_1)\bigcup_{ i_2} \partial_{t,e}  {\bf D}(i_1, i_2)\mbox{ is ${\bf L}^d$-
connected cubically}.\eqno{(2.13)}$$
By (2.12) -(2.13), and a simple induction, we have the following Lemma.\\ 

{\bf Lemma 3.} {\em For all configurations in ${\cal G}({\bf k},m)$,
$\bar{{\bf \Gamma}}_t({\bf k},m)$ is an ${\bf L^d}$-connected $t$-cube set. }\\

$\bar{\bf\Gamma}_t(\bf{k}, m)$ is defined to be a vertex set. On the other side,
we may view it as as an edge set by considering all edges with vertices in 
$\bar{\bf\Gamma}_t(\bf{k}, m)$.
With these definitions and lemmas, we would like to show the following fundamental geometric lemma.\\

{\bf Lemma 4.} {\em  For all configurations in ${\cal G}({\bf k},m)$, $\bar{{\bf \Gamma}}_t({\bf k},m)$ contains 
a zero-cutset cutting ${\bf B}({\bf k}, m)$ from $\infty$.}\\

{\bf Proof.}  The proof is a pure topology argument. It is easy to be convinced by  Fig. 2.
To show Lemma 4, we need to show that there is a cutset in $\bar{{\bf \Gamma}}_t({\bf k},m)$ and that all its edges are closed.
If we collect all  closed ${\bf Z}^d$-edges in $\bar{{\bf \Gamma}}_t({\bf k},m)$,  for any ${\bf Z}^d$-path $\gamma$ from $\infty$ to ${\bf B}({\bf k},m)$, 
we only need to show that $\gamma$ must use one of these closed edges. 
We will now go along $\gamma$ from $\infty$ to ${\bf B}({\bf k},m)$. By (2.6) and
the definition of the surface,
$\gamma$ must first  reach the surface ${\bf U}(\partial_{t,e} {\bf C}({\bf k},m))$.  Let ${\bf v}_1$ be the last vertex
of $\gamma$ at the surface such that the remaining piece of $\gamma$ from ${\bf v}_1$ to ${\bf B}({\bf k},m)$, denoted by
$\gamma_1$,  will not have  common vertices with the surface. 

 We suppose that in the following case a: $\gamma_1$  will not follow an open path to  a vertex of ${\bf D}(i_1)$ for all $1\leq i_1\leq l_1$.
 Thus, 
$\gamma_1$ will  use a closed  boundary edge of the open cluster containing ${\bf v_1}$ in $\partial_{t,e} \bar{\bf C}({\bf k}, m)$,  so Lemma 4 follows
for case a.

Now we suppose the following case b: $\gamma_1$ follows from an open path to a vertex  of  ${\bf D}(i_1)$ for some $i_1$. 
By (2.7),  $  {\bf D}(i_1)\cap {\bf B}({\bf k}, m)=\emptyset$. Thus,  together with
the assumptions that $k_i/t$ and $m/t$ are integers for $i=1, \cdots , d-1$,
  $\partial_{t, e} {\bf D}(i_1)$ and ${\bf B}({\bf k}, m)$ have no common cubes. So they do not have a common edge.
  Note that 
$\partial_{t, e}{\bf D}(i_1, \cdots i_j, i_{j+1}) \subset \partial_{t, e}{\bf D}(i_1, \cdots i_j, i_{j}),$
so 
$$\partial_{t, e}{\bf D}(i_1, \cdots i_j )\mbox{ and } {\bf B}({\bf k}, m) \mbox{ have no common edges for all }j \geq 1.\eqno{(2.14)}$$
Thus, $\gamma_1$ will leave any $\partial_{t, e}{\bf D}(i_1, \cdots i_j )$ for $j \geq 1$
before meeting ${\bf B}({\bf k}, k)$. Since we  will show that $\gamma_1$  uses a closed edge
in $\cup_{2\leq j\leq s}\cup_{i_j} \partial_{t, e} {\bf D}(i_1, i_2, \cdots, i_j)$, we only need to focus on the piece of $\gamma_1$ in the domain enclosed by
${\bf U}(\partial_{t, e}{\bf D}(i_1 ))$.  For simplicity,  we still use $\gamma_1$ to denote this piece.

 Let  $ q$ be the number of  all open clusters, constructed above,  in
$ \{{\bf D}(i_1, i_2, \cdots , i_j)\}$ for $1\leq i_j\leq l_j$ for $j=2, \cdots, s$.
We use an induction to show that $\gamma_1$ has to use a closed edge in 
$\cup_{2\leq j\leq s}\cup_{ i_j} \partial_{t, e} {\bf D}(i_1, i_2, \cdots, i_j)$. We first focus on $q=0$.  
By (2.14),
$\gamma_1$ has to eventually leave from $\partial_{t, e}{\bf D}(i_1)$.
By using  $q=0$ and ${\bf D}(i_1)\cap {\bf B}({\bf k}, m)=\emptyset$,   all the open paths from ${\bf U}(\partial_{t, e} {\bf D}(i_1))$ 
cannot penetrate  $\partial_{t, e} \bar{\bf D}(i_1)$.   In other words,
any path from  ${\bf U}(\partial_{t,e} {\bf D}(i_1))$  to a vertex of ${\bf D}(i_1)$   has to use a closed edge in $\partial_{t, e} \bar{\bf D}(i_1)$.
This implies that $\gamma_1$ will use a closed edge
inside $\partial_{t, e} \bar{\bf D}(i_1)\subset \bar{{\bf \Gamma}}({\bf k}, m)$.
If we replace ${\bf D}(i_1)$ by ${\bf D}(i_1, \cdots, i_j)$ with ${\bf D}(i_1, \cdots, i_j)\cap {\bf B}({\bf k}, m)=\emptyset$, this argument still holds for $q=0$.

For $q\leq l$ and for $i_1, \cdots, i_j$,  we consider open cluster ${\bf D}(i_1, \cdots, i_j)$ with
$${\bf D}(i_1, \cdots, i_j)\cap {\bf B}({\bf k}, m)=\emptyset.$$
We suppose that  $\gamma_1$ will come to use a vertex of  $ {\bf D}(i_1, \cdots, i_j)$  and will leave   $\partial_{t,e}{\bf D}(i_1, \cdots, i_{j-1})$, but not from its exits.
We make the induction hypothesis  that $\gamma_1$ has to use a closed edge in 
$\cup_{2\leq j\leq s}\cup_{ i_j} \partial_{t, e} {\bf D}(i_1, i_2, \cdots, i_j)$. 
Let us focus on  $q=l+1$.  There are two sub-cases:
case b (1), $\gamma_1$ will  first reach exits $\{{\bf u}_1(i_j), \cdots, {\bf u}_{q_j}(i_j)\}$ and   will use   vertices in ${\bf D}(i_1, i_2, \cdots i_{j})$   with ${\bf D}(i_1, i_2, \cdots i_j)\cap {\bf B}({\bf k}, m)=\emptyset$; 
or 
case b (2),   $\gamma_1$ will not.  Here we remark that if $\gamma_1$  enters 
an open cluster, but leaves from its exits and never comes back, then we still account 
$\gamma_1$ never uses the open cluster.

  In case b (1),    $\gamma_1$ will first enter to ${\bf D}(i_1, i_2,\cdots, i_j )$ from $\{{\bf u}_1(i_j), \cdots, {\bf u}_{q_j}(i_j)\}$  for some $j\geq 3$, and then will leave $\partial_{t,e}{\bf D}(i_1, i_2,\cdots, i_{j-1} )$, but not from $\{{\bf u}_1(i_j), \cdots {\bf u}_{q_j}(i_j)\}$. Note that the number of open clusters in $\{ {\bf D}(i_1, i_2, i_3, \cdots, i_j )\}$ 
 for $i_j\leq l_j$ and for $j\geq 3$ is strictly less than $q+1$. 
 In addition, 
$$ {\bf D}(i_1, i_2, \cdots, i_j)\cap {\bf B}({\bf k}, m)=\emptyset.$$
 Thus, by the induction hypothesis,  $\gamma_1$ will use a closed edge in 
 $\cup_{2\leq j\leq s}\cup_{ i_j} \partial_{t, e} {\bf D}(i_1, i_2, \cdots, i_j)$ in case b (1).
 
 In case (b) (2),  let us focus on $\{{\bf D}(i_1, i_2)\}$ for all $i_2$ with exits
 $\{{\bf u}_1(i_2), \cdots {\bf u}_{q_2}(i_2)\}$. If $\gamma_1$   never uses vertices
 of $\{{\bf D}(i_1, i_2)\}$, note that open clusters of the vertices on $U(\partial_{t, e}{\bf D}(i_1))\setminus 
 \{{\bf u}_1(i_2), \cdots {\bf u}_{q_2}(i_2)\}$ cannot penetrate $\partial_{t, e}{\bf D}(i_1)$, 
 so by using the same proof of  $q=0$, $\gamma_1$ will use a closed edge in 
 $\partial_{t, e}{\bf D}(i_1)$. Now we may assume that $\gamma_1$ will first reach  to  ${\bf D}(i_1, i_2)$ and will use at least one of its edges
 for some $i_2$ after ${\bf v_1}$.  Thus,  $\gamma_1$ will meet  the surface $U(\partial_{t, e}{\bf D}(i_1, i_2))$ at ${\bf u}$. Since  $\gamma_1$ is a path in  the domain enclosed by $U(\partial_{t, e}{\bf D}(i_1))$, ${\bf u}\not\in \{{\bf u}_1(i_2), \cdots {\bf u}_{q_2}(i_2)\}$. We denote by ${\bf D}({\bf u})$ the open cluster at ${\bf u}$ in the domain enclosed by $U(\partial_{t, e}{\bf D}(i_1, i_2))$.
 Since ${\bf u}\not\in \{{\bf u}_1(i_2), \cdots {\bf u}_{q_2}(i_2)\}$ and ${\bf u} \in U(\partial_{t, e}{\bf D}(i_1, i_2))$,  
 $$ {\bf D}({\bf u})\cap {\bf D}(i_1, i_2) =\emptyset .\eqno{(2.15)}$$
 Note that by (2.7),  $\{{\bf u}_1(i_2), \cdots {\bf u}_{q_2}(i_2)\}$ are the only possible vertices on ${\bf B}({\bf k}, m)$,  so by (2.15), ${\bf D}({\bf u})\cap {\bf B}({\bf k}, m)=\emptyset$. 
 Therefore, ${\bf D}({\bf u})\subset \partial_{t, e}{\bf D}(i_1, i_2)$;
 otherwise it would be one of open clusters of 
 $\{ {\bf D}(i_1, i_2, i_3, \cdots, i_j )\}$  for $i_j\leq l_j$ and $2\leq j\leq s$ with
 ${\bf D}(i_1, i_2, i_3, \cdots, i_j )\cap {\bf B}({\bf k}, m)=\emptyset$.
This contradicts the assumption of  case b (2). 
So $\gamma_1$ will use a boundary edge (a closed edge)  of ${\bf D}({\bf u})$ in
$\partial_{t, e}{\bf D}(i_1, i_2)$. Lemma 4 for case b, (2) follows. Together with case a,
case b (1) and case b (2), Lemma 4 is proved.
$\Box$\\

Since $\bar{\bf \Gamma}_t({\bf k},m)$ contains a closed cutset for ${\bf B}({\bf k}, m)$, we select a closed
self-avoiding cutset inside  $\bar{\bf \Gamma}_t({\bf k},m)$ using a unique method and denote it by $\bar{\bf \Gamma}({\bf k},m)$
with ${\bf Z}^d$-edges.
Now we will show another geometric property for $\bar{\bf \Gamma}_t({\bf k},m)$.  
A cube $B_t({\bf u})$
is said to have a  {\em blocked property} if  there is
 an open path inside $\bar{B}_t({\bf u})$ from  $B_t({\bf u})$ to the boundary of $\bar{B}_t({\bf u})$, 
without
connecting in $\bar{B}_t({\bf u})$ 
to a surface of  $t$-cubes in $\bar{B}_t({\bf u})$. For an independent purpose, we require that 
the above open paths will only use the edges in int$(\bar{ B}({\bf u}))$.
Intuitively, open paths are blocked to reach certain surfaces.
Note that $B_t({\bf u})$ is a blocked cube that only depends on configurations of edges in int$(\bar{B} _t({\bf u}))$.
For a fixed cube $B_t({\bf u})$, we say it has  a {\em disjoint property}
if there exist two disjoint open paths in $\bar{B}_t({\bf u})$
from  cube $B_t({\bf u})$ to
$\partial \bar {\bf B}_t({\bf u})$. Similarly, for an independent purpose, we require that 
the above open paths will  use the  edges in int$(\bar{ B}_t({\bf u}))$.
With this definition, $B_t({\bf u})$ has a disjoint property
 depending only on the configurations of edges in int$(\bar{B} _t({\bf u}))$.
 
 If a cube $B_t({\bf u})\in {\bf \Gamma}_t({\bf k},m)$, we shall show that it has 
a blocked  property. To see this, note that if $B_t({\bf u})\in \partial_{t,e} {\bf C}({\bf k},m)$, there is an open path from  ${\bf B}({\bf k}, m)$ to $B_t({\bf u})$, but cannot go out of ${\bf U}(\partial_{t,e} {\bf C}({\bf k},m))$, then it has a blocked property.  If
$B_t({\bf u})\in \partial_{t, e}{\bf D}(i_1, \cdots i_j)$, by (2.11), $B_t({\bf u})$ has an ${\bf L}^d$-neighbor cube  $B_t'({\bf u})$ with
$$ B_t'({\bf u})\cap {\bf D}(i_1, \cdots i_j)=\emptyset.\eqno{(2.15)}$$
In addition, there is an open in ${\bf D}(i_1, \cdots i_j)$  at least from $B_t({\bf u})$ to
$\bar{B}_t({\bf u})$. Thus,
  $B_t({\bf u})$ has a  block property.
We summarize this geometric property as the following lemma.\\

{\bf Lemma 5.} {\em For all configurations in ${\cal G}({\bf k},m)$,  the cubes in ${\Gamma}_t({\bf k},m)$ have a blocked property.}

\section{ Probability estimates for the linear zero-cutset.} 
In  section 3, we will first estimate the probabilities of events on Lemma 6.\\

{\bf Lemma 6.} {\em If $F(0) < 1-p_c$, then there exist $C_i=C_i(F(0), d)$ for $i=1,2$ such that for each cube $B_t({\bf u})$,}
$${\bf P}[ B_t({\bf u}) \mbox{ has a disjoint property}] \leq C_1 \exp(-C_2 t).$$

{\bf Proof.} The proof of Lemma 6 follows from Lemma 7.89 in Grimmett (1999). $\Box$\\

{\bf Lemma 7.}
{\em If $F(0) < 1-p_c$, then there exist $C_i=C_i(F(0), d)$ for $i=1,2$
 such that for each cube $B_t({\bf u})$,}
$${\bf P}[ B_t({\bf u}) \mbox{ has a blocked property}] \leq C_1 \exp(-C_2 t).$$

{\bf Proof.} By Lemma 7.104 in Grimmett (1999), if ${\cal A}_t({\bf u})$ is the event that
any two surfaces in the cubes of $\bar{B}_t({\bf u})$ are connected by open paths in $\bar{B}_t({\bf u})$, then 
$${\bf P} \left[ {\cal A}_t({\bf u})\right] \geq  1-C_1 \exp(C_2 t).$$
Now we suppose that there is an open path from $B_t({\bf u})$ to $\bar{B}_t({\bf u})$ for some ${\bf u}$, 
but the path cannot be
further connected to one of the surfaces in the cubes of  $\bar{B}_t({\bf u})$.  
We denote this event by ${\cal B}_t(u)$.
By the above inequality and Lemma 6,
\begin{eqnarray*}
&&{\bf P}\left[ {\cal B}_t(u)\right]\leq {\bf P}\left[ {\cal B}_t(u)\cap {\cal A}_t({\bf u})\right] +C_1\exp(-C_2  t)\\
&\leq & {\bf P}[B_t({\bf u})\mbox{ has a disjoint property}]   +C_1 \exp(-C_2 t)\leq  C_3\exp(-C_4 t).
\end{eqnarray*}
Lemma 7 follows from the two inequalities above. $\Box$\\

For a configuration $\omega$, recall that  ${\bf X} ({\bf k}, m)$ is the selected cutset 
with passage time $\chi({\bf k}, m)$.
We also set the following edge set (see Fig. 3) as the {\em surface} edges of ${\bf B}({\bf k}, m)$:
$$\alpha({\bf k}, m)=\{e: e \mbox{ is a ${\bf Z}^d$-edge in }\Delta {\bf B}({\bf k}, m)\}.$$
Clearly, $\alpha({\bf k}, m)$ is a cutset that cuts  ${\bf B}({\bf k}, m)$ from $\infty$.
 Therefore, 
$$\chi( {\bf k}, m)\leq \tau(\alpha({\bf k}, m)).\eqno{(3.1)}$$
Note that if $m \leq \min_{1\leq i\leq d-1}k_{i}$, 
there are at most $2d\|{\bf k}\|_v$ ${\bf Z}^d$-edges in $\alpha({\bf k}, m)$, so
$${\bf E}\chi({\bf k}, m)\leq {\bf E}\tau(\alpha({\bf k}, m))\leq 2d \|{\bf k}\|_v {\bf E}\tau (e).\eqno{(3.2)}$$
Also, with our moment assumption in (1.1), by a standard large deviation result, 
there exist $C_i=C_i(F, d)$ for $i=1,2$ such that
for all $u \geq 4d {\bf E}\tau (e) \|{\bf k}\|_v$,
$${\bf P}[\tau(\alpha({\bf k}, m))\geq u]\leq C_1\exp(-C_2 u).\eqno{(3.3)}$$
With these observations, we have the following lemma.\\

{\bf Lemma 8.} {\em If the conditions in Theorem 1 hold, and  $u \geq 4d \|{\bf k}\|_v{\bf E}\tau(e)$, then}
$${\bf P}[\chi({\bf k}, m) \geq u]\leq {\bf P}[\tau(\alpha({\bf k}, m))\geq u]\leq C_1 \exp(-C_2 u).$$

\section{Connectedness of cutsets.}
In section 4, we need to show that each self-avoiding cutset is connected. 
Beforehand, we will show a lemma.\\

{\bf Lemma 9.} {\em If ${\bf Z}({\bf k}, m)$ is a self-avoiding cutset with ${\bf Z}^d$-edges
that cuts ${\bf B}({\bf k}, m)$
from $\infty$,
then for each edge $e\in {\bf Z}({\bf k}, m)$ with two vertices $v_1(e)$ and $v_2(e)$,
 there exist disjoint paths $\gamma_1$ and $\gamma_2$
from $v_1(e)$ to ${\bf B}({\bf k}, m)$ and from $v_2(e)$ to $\infty$ without using ${\bf Z}({\bf k}, m)$}.\\

{\bf Proof.} For each $e\in {\bf Z}({\bf k}, m)$, note that the cutset is self-avoiding, so
${\bf Z}({\bf k}, m)\setminus e$ is not a cutset. 
There exists a path $\gamma$ without using ${\bf Z}({\bf k}, m)\setminus e$ from 
$ {\bf B}({\bf k}, m)$ to $\infty$. 
If $\gamma$ does not pass through $e$, then $\gamma$, without using ${\bf Z}({\bf k}, m)$,
connects $ {\bf B}({\bf k}, m)$ to $\infty$. This  contradicts  the assumption that ${\bf Z}({\bf k}, m)$
is a cutset. Therefore, $e\subset \gamma$
and  $e$ is the only edge of ${\bf Z}({\bf k}, m)$ contained in $\gamma$. Let
$v_1(e)$ and $v_2(e)$ be the two vertices of $e$. By this observation, there exist 
paths $\gamma_1$ and $\gamma_2$ from $v_1(e)$ to ${\bf B}({\bf k}, m)$, and from $v_2(e)$ to $\infty$, respectively, such that
$$\gamma_1\cap {\bf Z}({\bf k}, m)=v_1(e)\mbox{ and }\gamma_2\cap {\bf Z}({\bf k}, m)=v_2(e).$$
Therefore,  Lemma 9 is proved. $\Box$\\

For ${\bf Z}({\bf k}, m)$ defined in Lemma 9, 
let $\hat{\bf Z}({\bf k}, m)$ be all the vertices that are connected by ${\bf Z}^d$-paths to
${\bf B}({\bf k}, m)$ without using ${\bf Z}({\bf k}, m)$. Note that ${\bf B}({\bf k}, m)$
is ${\bf Z}^d$-connected, and so is $\hat{\bf Z}({\bf k}, m)$. 
Recall that $\Delta _e \hat{\bf Z}({\bf k}, m)$ is denoted by the ${\bf Z}^d$-edges between
$\partial_e \hat{\bf Z}({\bf k}, m)$ and $\partial_i \hat{\bf Z}({\bf k}, m)$.
For each edge $e\in {\bf Z}({\bf k}, m)$, as we proved in Lemma 10, there exist
$\gamma_1$ and $\gamma_2$ from $v_1(e)$ to ${\bf B}({\bf k}, m)$, and from $v_2(e)$ to $\infty$, respectively, 
without using ${\bf Z}({\bf k}, m)$.
Note also that 
$$\gamma_2\cap \hat{\bf Z}({\bf k}, m)=\emptyset;$$
otherwise, ${\bf Z}({\bf k}, m)$ would not be a cutset.
Therefore, $v_1(e)\in \hat{\bf Z}({\bf k}, m)$, but $v_2(e)\not\in \hat{\bf Z}({\bf k}, m)$ and is connected by a path $\gamma_2$ 
without using
an edge of ${\bf Z}({\bf k}, m)$ from $v_2(e)$ to $\infty$. This implies that $e\in \Delta_e \hat{\bf Z}({\bf k}, m)$, so
$$ \mbox{ all edges of } {\bf Z}({\bf k}, m)\subset \Delta_e \hat{\bf Z}({\bf k}, m).\eqno{(4.1)}$$
For $e\in \Delta_e \hat{\bf Z}({\bf k}, m)$,
by the definition of $\hat{\bf Z}({\bf k}, m)$, 
$v_1(e)$ is connected to ${\bf B}({\bf k}, m)$ by $\gamma_1$ without using an edge of 
${\bf Z}({\bf k}, m)$, and $v_2(e)\not \in \hat{\bf Z}({\bf k}, m)$. This tells us that $e\in {\bf Z}({\bf k}, m)$,
since, otherwise, $\gamma_1\cup \{e\}$ would be a path that does not use ${\bf Z}({\bf k}, m)$ from $v_2(e)$
to ${\bf B}({\bf k}, m)$. So $v_2(e)\in \hat{\bf Z}({\bf k}, m)$.
Therefore, 
$$ \Delta_e \hat{\bf Z}({\bf k}, m)\subset \mbox{ all edges of } {\bf Z}({\bf k}, m). \eqno{(4.2)}$$
By (4.1) and (4.2), we have
$$ \Delta_e \hat{\bf Z}({\bf k}, m)= \mbox{ all edges of } {\bf Z}({\bf k}, m). \eqno{(4.3)}$$

By Lemma 1, $\partial_e \hat{\bf Z}({\bf k}, m)$ is ${\bf Z}^d$-connected. 
By (4.3), each vertex of $\partial_e \hat{\bf Z}({\bf k}, m)$ is either ${\bf Z}^d$-adjacent to $\hat{\bf Z}({\bf k}, m)$
by  a ${\bf Z}^d$-edge in ${\bf Z}({\bf k}, m)$ or  ${\bf L}^d$-adjacent to $\hat{\bf Z}({\bf k}, m)$.
Suppose that  ${\bf v}\in \partial_e \hat{\bf Z}({\bf k}, m)$ is only ${\bf L}^d$-adjacent to $\hat{\bf Z}({\bf k}, m)$, 
but is not ${\bf Z}^d$-adjacent. It is easy to verify (see Fig. 1) that one of its ${\bf L}^d$-neighbors is
${\bf Z}^d$-adjacent to $\hat{\bf Z}({\bf k}, m)$. In other words, one of its ${\bf L}^d$-neighbors is adjacent to 
$\hat{\bf Z}({\bf k}, m)$ by an edge in ${\bf Z}({\bf k}, m)$. Let us account for the number of 
 ${\bf L}^d$-neighbors for a vertex. Without loss of generality, we account for the origin.
 We assume that $(x_1, \cdots, x_d)$ is an ${\bf L}^d$-neighbor of 
the origin. Thus, $x_i$ can take either $\pm 1$ and zero. Hence, there are at most $3^{d}$ ${\bf L}^d$-neighbors
for the origin.  With this observation, for each ${\bf Z}^d$-edge $e$ in ${\bf Z}({\bf k}, m)$,
there are at most $3^{d+1}$ vertices in $\partial_e \hat{\bf Z}({\bf k}, m)$ that are ${\bf L}^d$-adjacent  to $e$.
With Lemma 1 and with these observations above, we have the following lemma to show the connectedness of cutsets.\\

{\bf Lemma 10.} {\em  If ${\bf Z}({\bf k}, m)$ 
is a self-avoiding cutset that cuts ${\bf B}({\bf k}, m)$
from $\infty$, then $\partial_e \hat{\bf Z}({\bf k}, m)$ is ${\bf Z}^d$-connected and }
$$ \mbox{ the number of ${\bf Z}^d$-edges in }{\bf Z}({\bf k}, m) \geq | \partial_e \hat{\bf Z}({\bf k}, m)|/ 3^{d+1}.\eqno{(4.4)}$$

Now we focus on the connectedness of the cutsets that cut ${\bf F}_0$ from ${\bf F}_m$. 
Let ${\bf V}({\bf k}, m)$ be a self-avoiding cutset that cuts ${\bf F}_0$ from ${\bf F}_m$.
Similarly, 
let  $\hat {\bf V}({\bf k}, m)$ be the all vertices in ${\bf B}({\bf k}, m)$ 
 that are connected by ${\bf Z}^d$-paths in ${\bf B}({\bf k}, m)$ to ${\bf F}_0$
without using  ${\bf V}({\bf k}, m)$.  In addition, let $\Delta \hat {\bf V}({\bf k}, m)$ be the all boundary 
edges of $\hat {\bf V}({\bf k}, m)$ in ${\bf B}({\bf k}, m)$.  For each edge $e\in \Delta \hat {\bf V}({\bf k}, m)$,
if there exists a path in ${\bf B}({\bf k}, m)$ from one of its vertices to ${\bf F}_m$ 
without using edges in ${\bf V}({\bf k}, m)$, then $e$ is an exterior boundary edge of 
$\hat {\bf V}({\bf k}, m)$. We denote  by $\Delta_e \hat {\bf V}({\bf k}, m)$ all the exterior boundary edges
of $\hat{\bf V}({\bf k}, m)$.
By the same proof as (4.3), we can show that
$$ \Delta_e \hat{\bf V}({\bf k}, m)= \mbox{ all edges of } {\bf V}({\bf k}, m). \eqno{(4.5)}$$
Kesten  (Lemma 3.17 in Kesten (1987)) showed 
that  for $d=3$, $\partial_e \hat{\bf V}({\bf k}, m)$ is ${\bf Z}^d$-connected; 
but his proof can be directly adapted to apply for all $d \geq 3$. On the other hand, it can also use the same proof
of Lemma 1 to show the ${\bf Z}^d$-connectedness of $\partial_e \hat{\bf V}({\bf k}, m)$.
By the same discussion
of (4.4), we can work on the number of ${\bf Z}^d$-edges of ${\bf V}({\bf k}, m)$. 
We summarize the above results as the
following lemma.\\

{\bf Lemma 11.} {\em If ${\bf V}({\bf k}, m)$ is a self-avoiding cutset that cuts ${\bf F}_0$ from ${\bf F}_m$
in ${\bf B}({\bf k}, m)$, then $\partial_e \hat{\bf V}({\bf k}, m)$ is ${\bf Z}^d$-connected and }
$$ \mbox{ the number of ${\bf Z}^d$-edges in }{\bf V}({\bf k}, m) \geq | \partial_e \hat{\bf V}({\bf k}, m)|/ 3^{d+1}.\eqno{(4.6)}$$

\section{Proof of Theorem 1.} 
If $\tau(e)$ only takes 0 or 1, Theorem 1 can be directly shown by using  the connectedness in Lemma 3, the blocked property in Lemma 5 and the exponential small probability estimate for the block property
in Lemma 7 together with  a Peierls argument in (5.26). If $\tau(e)$ can  take any non negative values, the proof of Theorem 1  is much more complicated. 
 In this section, we assume that $F(0) < 1-p_c$. For each $0<\epsilon <1$, 
$e$ is  said to be an $\epsilon^+$-edge or $\epsilon^-$-edge if  $\tau(e) > \epsilon$ or $0< \tau(e) \leq \epsilon$.
Let
$N^+({\bf k}, m)$ and $N^-({\bf k}, m)$ be  the numbers of $\epsilon^+$-edges  and $\epsilon^-$-edges
in ${\bf X}({\bf k}, m)$, respectively.
Note that
$$\epsilon N^+({\bf k}, m)\leq \chi({\bf k}, m),$$
so if we take $\beta_1= 4d{\bf E}\tau (e)$, by    Lemma 8 for $n \geq \epsilon^{-2}\beta_1 \|{\bf k}\|_v$
there exist $C_i=C_i(F, d, \epsilon)$ for $i=1,2$,
$${\bf P}\left[ N^+({\bf k}, m) \geq  \epsilon n\right]\leq {\bf P}\left[\chi({\bf k}, m)\geq \epsilon^2 n\right]\leq C_1\exp(-C_2 n).\eqno{(5.1)}$$ 

Now we take care of the $\epsilon^-$-edges in the cutset. By our definition,
$${\bf P}[ e \mbox{ is an $\epsilon^-$ edge}] \leq F(\epsilon)-F(0)=\delta_1=\delta_1 (\epsilon),\eqno{(5.2)}$$
where $\delta_1\rightarrow 0$ as $\epsilon \rightarrow 0$.
 We need to fix a vertex of ${\bf X}({\bf k}, m)$.
Since ${\bf X}({\bf k}, m)$ is a cutset, it must intersect the line ${\bf L}$:
$${\bf L}=\{(x_1, x_2,\cdots, x_d): x_i=0 \mbox{ for } i\geq 2\}.$$
We let ${\bf z}=(x_1,0,\cdots,0)$ be   the intersection vertex of
${\bf X}({\bf k}, m)$ and ${\bf L}$. 
If there are many intersections, we select one with the largest $x_1$-coordinate and still denote it by ${\bf z}$.
Note that if
$\mbox{{\bf dist}}({\bf 0}, {\bf z}) = l$ for some $l$,  then
$$\mbox{ the number of  edges } {\bf X}({\bf k},m) \mbox{ is larger than } l.\eqno{(5.3)}$$
To show (5.3), simply note that each layer between the hyperplanes $x_1=i$ and $x_1=(i+1)$ for $i\leq l$
contains at least one edge of ${\bf X}({\bf k},m)$.

Now we estimate the following probability for small $\delta_1$ defined in (5.2) and for a constant 
$D=D(d)$ selected later:
$${\bf P}\left[ N({\bf k}, m) \geq n, N^-({\bf k}, m)\geq -(D \log^{-1} (\delta_1))N({\bf k}, m)\right].$$
By (5.3), we fix ${\bf z}$ to have
\begin{eqnarray*}
&&{\bf P}\left[ N({\bf k}, m) \geq n, N^-(\epsilon, {\bf k}, m)\geq -(D \log^{-1} (\delta_1))N({\bf k}, m)\right]\\
&&=\sum_{j=n}^{\infty} \sum_{i=0}^j {\bf P}\left[ N({\bf k}, m) =j,  x_1=i, N^-({\bf k}, m)\geq - Dj\log^{-1} (\delta_1)\right].
\end{eqnarray*}
Recall $\hat{\bf X}({\bf k}, m)$ and $\partial_e \hat{\bf X}({\bf k}, m)$ defined above (see (4.2)).
As we defined before, ${\bf X}({\bf k}, m)$ is unique for each configuration, and
so is $\partial_e \hat{\bf X}({\bf k}, m)$. Thus, for two different  fixed
sets $\Gamma_1$ and $\Gamma_2$, we have
$$\{ \partial_e \hat{\bf X}({\bf k}, m)=\Gamma_1\}\mbox{ and } \{ \partial_e \hat{\bf X}({\bf k}, m)=\Gamma_2\} \mbox{ are disjoint}.$$
If $ \partial_e\hat{\bf X}({\bf k}, m)=\Gamma_1$, we say that $\partial_e\hat{\bf X}({\bf k}, m)$ has a choice $\Gamma_1$. 

If $x_1=i$,  then by Lemma 11, $\partial_e \hat{\bf X}({\bf k}, m)$ is ${\bf Z}^d$-connected and 
$$|\partial_e \hat{\bf X}({\bf k}, m)|\leq 3^{d+1} j.$$
Thus,  by  using  (4.24) in Grimmett (1999), 
there are at most $7^{d3^{d+1} j}$ choices of these  cutsets for  $\partial_e \hat{\bf X}({\bf k}, m)$
 when ${\bf X}({\bf k}, m)$ has $j$ edges.  After $\partial_e \hat{\bf X}({\bf k}, m)$ is a fixed vertex set,
we select the vertices in $\partial_e \hat{\bf X}({\bf k}, m)$ with ${\bf Z}^d$-edges of ${\bf X}({\bf k}, m)$.
We next select the ${\bf Z}^d$-edges of ${\bf X}({\bf k}, m)$ adjacent to these vertices.
Note that each vertex has at most $2d$ adjacent edges, so there are at most
$$\sum_{k=1}^{3^{d+1}j}{ 3^{d+1} j \choose k} (2d)^j\leq 2^{3^{d+1}j}(2d)^j$$
for the selections.
Thus, if $|\Gamma|_e$ is denoted by the number of edges in $\Gamma$,
\begin{eqnarray*}
&&{\bf P}\left[ N({\bf k}, m) \geq n, N^-( {\bf k}, m)\geq -(D\log^{-1} (\delta_1))N({\bf k}, m)\right ]\\
&&=\sum_{j=n}^{\infty} \sum_{i=0}^j {\bf P}\left [ N({\bf k}, m) =j,  x_1=i, N^-({\bf k}, m)\geq - Dj
\log^{-1} (\delta_1)\right]\\
&&\leq \sum_{j=n}^{\infty} j   7^{ d3^{d+1}j}2^{3^{d+1}j}(2d)^j\max_{\Gamma} {\bf P}\left[|\Gamma|_e=j, 
\mbox{ $\Gamma$ contains more than $-Dj \log^{-1} (\delta_1)$ $\epsilon^--$edges}\right]\hskip 0.1cm (5.4)
\end{eqnarray*}
where $\Gamma$  is a fixed cutset that cuts ${\bf B}({\bf k}, m)$ from $\infty$ with a number of edges
$|\Gamma|_e=j$, and the maximum is taking over all the
possible $\Gamma$. Let us estimate the probability in (5.4) for a fixed set $\Gamma$ 
with more than  $-Dj\log^{-1} (\delta_1)$ $\epsilon^-$-edges for a large number $D$.
$$
{\bf P}\left[ |\Gamma|_e=j, 
\mbox{ $\Gamma$ contains more than $- (Dj)\log^{-1} (\delta_1)$ $\epsilon^-$-edges}\right]
\leq \!\!\!\!\!\sum_{i\geq - Dj\log^{-1} (\delta_1)}^{j} {j \choose i} \delta_1^i. \eqno{(5.5)}
$$
By using Corollary 2.6.2 in Engle (1997), there are at most 
$${j \choose i} \leq \exp\left ( j H({i\over j})\right),\eqno{(5.6)}$$
where
$$H(x)=-x \log x -(1-x)\log (1-x).$$
By (5.6), we have
$$\sum_{i\geq  - Dj\log^{-1} (\delta_1)}^{j} {j \choose i} \delta_1^i\leq
\sum_{i\geq -Dj \log^{-1} (\delta_1)}^{j} \exp(j H(i/j))\delta_1^i.\eqno{(5.7)}$$
If $0<x<1$, we have
$$H(x) \leq 2x \log (1/x).\eqno{(5.8)}$$
By (5.7) and (5.8),
$$\sum_{i\geq  - Dj\log^{-1} (\delta_1)}^{j} {j \choose i} \delta_1^i\leq
\sum_{i\geq  - Dj\log^{-1} (\delta_1)}^{j} \exp\left[2i\log (j/i)+ i \log (\delta_1)\right].\eqno{(5.9)}$$
Note that $-Dj \log^{-1} (\delta_1)\leq i\leq {j}$ and $\log (j/i)$ is decreasing when $i$ is increasing until  $j$, 
so for small $\delta_1$ and $D\geq 1$, we have
$$\log {j\over i} \leq  \log {j\over - Dj\log^{-1} (\delta_1)}.\eqno{(5.10)}$$
Note also that if $\delta_1\rightarrow 0$, then
$${\log (-\log(\delta_1))\over \log \delta_1}\rightarrow 0\mbox{ from the left},$$
so we take a small $\delta_1$,
$$\log (-\log(\delta_1))\leq -(\log \delta_1)/4.$$
Hence, by (5.9)
 $$\sum_{i\geq - Dj\log^{-1} (\delta_1)}^{j} {j \choose i} \delta_1^i\leq 
\sum_{i\geq - Dj\log^{-1} (\delta_1)}^{j}\exp({i\over 2} \log (\delta_1)).\eqno{(5.11)}$$
We use (5.11) in (5.5) to produce
\begin{eqnarray*}
&&{\bf P}\left[ |\Gamma|_e=j, 
\mbox{ $\Gamma$ contains more than $- Dj\log^{-1} (\delta_1)$ $\epsilon^-$-edges}\right]\\
&&=\sum_{i\geq -Dj \log^{-1} (\delta_1)}^{j} {j \choose i} \delta_1^i \\
&&\leq \sum_{i\geq -Dj\log^{-1} (\delta_1)}^{\infty} \exp(i\log (\delta_1)/2)\\
&&\leq 2 \exp(-Dj/2).
\end{eqnarray*}
With this observation and (5.4),
\begin{eqnarray*}
&&{\bf P}\left[ N({\bf k}, m) \geq n, N^-( {\bf k}, m)\geq -(D\log^{-1} (\delta_1))N({\bf k}, m)\right ]\\
&&=\sum_{j=n}^{\infty} \sum_{i=0}^j {\bf P}\left [ N({\bf k}, m) =j, x_1=i, N^-({\bf k}, m)\geq -Dj 
\log^{-1} (\delta_1)\right]\\
&&\leq \sum_{j=n}^{\infty} j   7^{d3^{d+1} j}2^{3^{d+1}j}(2d)^j\max_{\Gamma} {\bf P}\left[ \Gamma, |\Gamma|_e=j, 
\mbox{ $\Gamma$ contains more than $-Dj \log^{-1} (\delta_1)$ $\epsilon^--$edges}\right]\\
&&\leq \sum_{j=n}^{\infty} j   7^{d3^{d+1}j}2^{3^{d+1}j}(2d)^j\exp(-Dj/2).
\end{eqnarray*}
By taking $D=D(d)$ and a small $\delta_1$, there are $C_i=C_i(F,d, \delta_1)$ for $i=1,2$ such that
$${\bf P}\left[ N({\bf k}, m) \geq n, N^-( {\bf k}, m)\geq -(D\log^{-1} (\delta_1))N({\bf k}, m)\right ]
\leq C_1 \exp(-C_2 n).\eqno{(5.12)}$$

For a configuration $\omega$,  we denote edges $\epsilon^\pm$ in the cutset by $e_1,\cdots, e_J$. 
Therefore,
$$\chi({\bf k}, m)=\sum_{i=1}^J \tau(e_i).\eqno{(5.13)}$$
By (5.1) and (5.12), for all small  $\epsilon$ and corresponding $\delta_1$ and $\beta_1$, and
for all  $n\geq \epsilon^{-2}\beta_1 \|{\bf k}\|_v$, there are $C_i=C_i(F, d, \epsilon)$ for $i=1,2$ such that
\begin{eqnarray*}
&&{\bf P}\left[ N({\bf k}, m) \geq n\right]\\
&\leq&  \sum_{j\geq n}{\bf P}\left[ N({\bf k}, m) =j, N^+({\bf k}, m) \leq  \epsilon j, N^-({\bf k}, m)\leq -(Dj)\log^{-1} (\delta_1) \right]+  C_1 \exp(-C_2 n).\hskip 0.2cm (5.14)
\end{eqnarray*}
On $\left\{N^+({\bf k}, m) \leq  \epsilon j, N^-({\bf k}, m)\leq -(Dj)\log^{-1} (\delta_1) \right\}$,
 we have
$$J \leq \left(\epsilon -D\log^{-1} (\delta_1)\right) j.\eqno{(5.15)}$$
For a large number $t$, by (5.15), we take $\epsilon$ small and then $n$ large such that
$$
{\bf P}\left[ N({\bf k}, m) \geq n\right]
\leq  \sum_{j\geq n}{\bf P}\left[ N({\bf k}, m) =j, J\leq j/ (4 (16t)^d)\right]+  C_1 \exp(-C_2 n)\eqno{(5.16)}
$$
for $C_i=C_i(F, d, \epsilon, t)$ with $i=1,2$.
Note that the edges in ${\bf X}({\bf k}, m)$ other than  $e_i$ for $i=1,\cdots, J$ are zero-edges, so
if we change these $\epsilon^\pm$-edges from positive to zero,  we will have
a closed cutset corresponding to another configuration $\omega'$. 
More precisely, for each configuration $\omega$, if we make the changes for these $\epsilon^\pm$-edges,
then $\omega$ will change to another configuration $\omega'$. Since ${\bf X}({\bf k}, m)(\omega)$ is uniquely selected,
$\omega'$ is determined uniquely for each $\omega$.
If there exists a closed cutset for $\omega'$,  by Lemma 4,
 there exists the self-avoiding zero-cutset $\bar{\bf \Gamma}({\bf k}, m)(\omega')$
inside $\bar{\bf \Gamma}_t({\bf k}, m)(\omega')$.  
By Lemma 3, $\bar{\bf \Gamma}_t({\bf k}, m)(\omega')$ is ${\bf L}^d$-connected cubically. Therefore, for each configuration $\omega$,  
$e_1,\cdots, e_J$ exist. So we have $\omega'$ and  $\bar{\bf \Gamma}({\bf k}, m)(\omega')$.  
If we change these edges in $\{e_1,\cdots, e_J\}$ from zero back to the original values, 
 $\bar{\bf \Gamma}({\bf k}, m)$, as a vertex set, exists  corresponding to $\omega$, 
but it will no longer  be a closed cutset.
We denote it by $\bar{\bf \Gamma}({\bf k}, m)(\omega)$, as this vertex set for the configuration $\omega$.
Note that $\bar{\bf \Gamma}({\bf k}, m)(\omega)$ is uniquely determined for each $\omega$.
We claim that for each configuration,
 $$\tau(\bar{\bf \Gamma}({\bf k}, m)(\omega))=\chi({\bf k}, m)(\omega).\eqno{(5.17)}$$
 To see this, note that $\bar{\bf \Gamma}({\bf k}, m)(\omega)$ is a cutset, so
$$\tau(\bar{\bf \Gamma}({\bf k}, m)(\omega))\geq \chi ({\bf k}, m)(\omega).\eqno{(5.18)}$$
On the other hand,  the other edges in $\bar{\bf \Gamma}({\bf k}, m)(\omega)$, except for $e_1,\cdots, e_J$, are all zero-edges,  
and ${\bf X}({\bf k}, m)$ uses all the edges $e_1,\cdots, e_J$, so
$$\tau(\bar{\bf \Gamma}({\bf k}, m)(\omega))\leq \tau({\bf X}({\bf k}, m))=\chi({\bf k}, m)(\omega).\eqno{(5.19)}$$
Therefore, (5.17) follows.
By the definition of $N({\bf k}, m)$ and (5.17), note that $\bar{\bf \Gamma}({\bf k}, m)$ as
a vertex set is the same for either $\omega$ or $\omega'$, so we have
$$|\bar{\bf \Gamma}({\bf k}, m)(\omega')|_e=|\bar{\bf \Gamma}({\bf k}, m)(\omega)|_e\geq N({\bf k}, m)(\omega) . \eqno{(5.20)}$$

If these edges $e_1$, $\cdots$, $e_J$ are zero-edges, as we mentioned above, $\bar{\bf \Gamma}({\bf k}, m)$ 
is a zero-cutset contained inside $\bar{\bf \Gamma}_t({\bf k}, m)$.
Since  $\bar{\bf \Gamma}_t({\bf k}, m)$ is
also a cutset, $\bar{\bf \Gamma}_t({\bf k}, m)$ and ${\bf L}$ must intersect. We denote by $B_t({\bf z})(\omega')$ the cube
in $\bar{\bf \Gamma}_t({\bf k}, m)(\omega')$ that intersects ${\bf L}$. If there are many such cubes, we simply select
${\bf z}$ with the largest $x_1$ value. By the same argument of (5.3), if $l$ is denoted the number of cubes 
with the lower corners at  ${\bf L}$ from the origin to ${\bf z}$, then
$$\mbox{the number of cubes in $\bar{\bf \Gamma}_t({\bf k}, m)(\omega')$ is larger than $l$}.\eqno{(5.21)}$$
For $\omega'$ and  each $t$-cube $B_t({\bf u})$ in $\bar{\bf \Gamma}_t({\bf k},m)(\omega')$, by Lemma 5, there is a $t$-cube in $\bar{B}_t({\bf u})$ that
 has the blocked  property. By our definition,  if  $B_t({\bf u})$ and
$B_t({\bf v})$ for fixed ${\bf u}$ and ${\bf v}$ satisfy that
$$\mbox{int}(\bar{B}_t({\bf u}))\cap \mbox{int}(\bar{B}_t({\bf v}))=\emptyset,$$
then
$$\{B_t({\bf u})\mbox{ has a blocked  property} \} \mbox{ and }\{B_t({\bf v})\mbox{ has a blocked  property}\}$$
are independent. 
Therefore, we need to decompose $\bar{\bf \Gamma}_t({\bf k}, m)(\omega')$ into $3t$-cubes such that their center cubes with blocked property
belong to $\bar{\bf \Gamma}_t({\bf k}, m)(\omega')$. By (2.3) and a standard  estimation (see Grimmett and Kesten, page 345 (1984) or
Zhang, pages 21 (2008), or  Steele and Zhang (2003), Lemma 6 by using Turan's theorem), if 
the number of cubes in $\bar{\bf \Gamma}_t({\bf k}, m)(\omega')$ is $s$, then
\begin{eqnarray*}
&&\exists \mbox{ at least $s/ 2^{4d}$  disjoint $3t$-cubes such that  their center cubes 
with block property}\\
&&\mbox{in  $\bar{\bf \Gamma}_t({\bf k}, m)(\omega')$}.\hskip 12cm {(5.22)}
\end{eqnarray*}
Let ${\bf M}_{3t}({\bf k},m)(\omega')$ be all the disjoint $3t$-cubes in (5.22) and 
${ M}_{3t}({\bf k},m)(\omega')$ be the number of the
$3t$-cubes in ${\bf M}_{3t}({\bf k},m)(\omega')$.
 Note that each $t$-cube has $2t^d$ edges, so by (5.20) and (5.22), if
$N({\bf k}, m)(\omega)=j$, then
$$\!\!\!\!\!\!{M}_{3t}({\bf k},m)(\omega')\geq {\#\mbox{ cubes in }\bar{\bf \Gamma}_t({\bf k}, m)(\omega')\over 2^{4d}}\geq {|\bar{\bf \Gamma}({\bf k}, m)(\omega')|_e\over 2^{4d+1} t^d}\geq 
{N({\bf k}, m)(\omega)\over 2(16t)^d}={j\over 2(16t)^d} .\eqno{(5.23)}$$

Furthermore, if
$$J(\omega) \leq j /(4(16t)^d)\mbox{ and }N({\bf k}, m)(\omega)=j \mbox{ and }\#\mbox{ cubes in }\bar{\bf \Gamma}_t({\bf k}, m)(\omega')=s,$$
then by (5.22) and (5.23), for each $\omega'$, there are at least
$${ M}_{3t}({\bf k},m)-{j\over (4(16t)^d)}= {{ M}_{3t}({\bf k},m)\over 2}+ {{ M}_{3t}({\bf k},m)\over 2}-{j \over (4 (16t)^d)}\geq { { M}_{3t}({\bf k},m)\over 2} \geq {s\over 2^{4d+1}} \eqno{(5.24)}$$
 disjoint  $3t$-cubes in 
${\bf M}_{3t}({\bf k}, m)(\omega')$ such that
their center $t$-cubes have  the blocked property,
and these interior $3t$-cubes do not contain  $e_1,\cdots, e_J$. 
Note that these disjoint $s/2^{4d+1}$ cubes always have the blocked  property whether $\tau(e_i)$ is positive or
zero for $i=1,2, \cdots, J$, since they do not contain these edges in  their interiors. 
We call them  {\em permanent} blocked cubes.
Now we change these edges in $\{e_1,\cdots, e_J\}$ from zero back to the original values. We still have 
$s/2^{4d+1}$ permanent block or disjoint $3t$-cubes.

In summary, for each $\omega$, if $N({\bf k}, m)(\omega)=j,$
and $J(\omega)\leq j /(4(16t)^d)$, note that by (5.20),
$$\#\mbox{ cubes in }\bar{\bf \Gamma}_t({\bf k}, m)(\omega')=s\geq |\bar{\bf \Gamma}({\bf k}, m)(\omega)|_e/(2t)^d\geq N({\bf k}, m)(\omega)/ (2t)^d =j /(2t)^d, \eqno{(5.25)}$$
so by (5.21)--(5.24),
there are  ${\bf L}^d$-connected $s\geq j/(2t)^d$  $t$-cubes containing $B_t({\bf z})$ 
for ${\bf z}\in {\bf L}$ with $\|{\bf z}\|\leq s$
such that\\
(a) there are at least $s/2^{4d+1}$ disjoint $3t$-cubes containing the above $t$-cubes as their center cubes and \\
(b)  each  center $t$-cube in these $3t$-cubes in (a) has the  permanent blocked property,\\
where the blocked property in (b) corresponds to the configuration $\omega$.
We denote the event of (a) and (b) by ${\cal E}(s, j, {\bf z})$.

Now we  try to estimate ${\cal E}(s, j, {\bf z})$ by fixing these $3t$-cubes in the following steps.  
We suppose that the connected $t$-cubes in event ${\cal E}(s, j, {\bf z})$ is ${\bf \Lambda}_t$ with $s$ cubes.
First, we fix $B_t({\bf z})$.  By (5.21), there are at most $s$ choices for this cube. With this cube, note that  ${\bf \Lambda}_t$ is ${\bf L}^d$-connected, so
by using  a standard computation technique (see (4.24) in Grimmett (1999)),  there are
  at most $s 7^{2ds}$ choices for this set ${\bf \Lambda}_t$. If
${\bf \Lambda}_t$ is fixed, we select these disjoint $3t$-cubes in ${\bf \Lambda}_t$
in (a) and (b).
There are at most
$$\sum_{i=1}^s {s\choose i}=2^s$$
choices for these $3t$-cubes. 
If a $t$-cube $B_t({\bf u})$  is not ${\bf L}^d$-connected to ${\bf B}({\bf k}, m)$, by Lemma 7,
there exist  $C_i=C_i(F(0), d)$ for $i=3,4$ such that
$${\bf P} [ B_t({\bf u}) \mbox{ has a blocked property}]\leq  C_3 \exp(-C_4 t).\eqno{(5.26)}$$
If  $B_t({\bf u})$  is  ${\bf L}^d$-connected to ${\bf B}({\bf k}, m)$, by the same proof of Lemma 7, (5.26) still holds.  Therefore, by the observations above,
there are $C_i=C_i(F, d, \epsilon, t)$ for $i=1,2$ and $C_i=C_i(F(0), d)$
for $i=3,4$ such that
\begin{eqnarray*}
&&{\bf P}\left[ N({\bf k}, m) \geq n\right]\\
&&\leq  \sum_{j\geq n}{\bf P}\left [ N({\bf k}, m) =j, J\leq j/ (2 (16t)^d)\right]+C_1\exp(-C_2n).\\
&&\leq \sum_{j\geq n} \sum_{\scriptstyle s\geq j/(2t)^d \atop \scriptstyle \|{\bf z}\|\leq s} {\bf P}\left[ {\cal E}(s,j,{\bf z})\right]\\
&&\leq \sum_{j\geq n}\sum_{s\geq {j/(2t)^d}} s^27^{2ds}2^s  \left[ C_3\exp(-C_4t)\right]^{s/2^{4d+1}}+C_1\exp(-C_2n).
\end{eqnarray*}
If we take $t$ large (but only dependent  on $d$), there exist $C_i=C_i( F, d, \beta, \epsilon,t)$ for $i=1,2$ such that
$${\bf P}[ N({\bf k}, m) \geq n]\leq C_1 \exp(-C_2 n).\eqno{(5.26)}$$
The approach in (5.26) is commonly called a {\em Peierls argument}.
Therefore, Theorem 1 follows.

\section{ Proof of Theorem 2.} 
Since  section 6  focuses on the edges inside ${\bf B}({\bf k}, m)$, we use
${\bf P}_{{\bf k}, m}(\cdot)$ to denote  the probability measure. In addition, we assume that $F(0) < 1-p_c$
in this section.
Let $\bar{\alpha}({\bf k}, m)$ be vertical edges between   $\prod_{i=1}^d[0, k_i]\times \{0\}$ and 
$\prod_{i=1}^d[0, k_i]\times \{1\}$ inside ${\bf B}({\bf k}, m)$. 
Note that $\bar{\alpha}({\bf k})$ is a cutset that cuts ${\bf F}_0$ from ${\bf F}_m$, so
$$\tau ({{\bf W}}({\bf k},  m))\leq \tau(\bar{\alpha}({\bf k}, m)).\eqno{(6.0)}$$
By a similar  large deviation
result for $\alpha({\bf k}, m)$ in the last section, 
note that there are $\|{\bf k}\|_v$ edges in $\bar{\alpha}({\bf k}, m)$.
Thus, if (1.1) holds and $u\geq 2 {\bf E}\tau(e) \|{\bf k}\|_v$, then there are $C_i=C_i(F, d)$ such that
$${\bf P}_{{\bf k}, m}\left[\tau ({{\bf W}}({\bf k}, m))\geq u\right]\leq {\bf P}_{{\bf k}, m}\left[\tau (\bar{\alpha} ({\bf k}, m))\geq u\right]\leq C_1\exp(-C_2 u).\eqno{(6.1)}$$

Recall that  
$e$ is  said to be an $\epsilon^+$-edge or $\epsilon^-$-edge if  $\tau(e) > \epsilon$ or $0< \tau(e) \leq \epsilon$.
Let
${N}^+({\bf k}, m)$ and ${N}^-({\bf k},  m)$ be  the numbers of $\epsilon^+$-edges  and $\epsilon^-$-edges
in ${{\bf W}}({\bf k}, m)$, respectively.
Note that
$$\epsilon {N}^+({\bf k}, m)\leq \tau({\bf W}({\bf k}, m)),$$
so if we take $\beta_1= 2{\bf E}\tau (e)$, by   (6.1) for $n \geq \epsilon^{-2}\beta_1 \|{\bf k}\|_v$
$${\bf P}_{{\bf k}, m}\left[ N^+({\bf k}, m) \geq  \epsilon n\right]\leq C_1\exp(-C_2 n).\eqno{(6.2)}$$ 

Now we take care of the $\epsilon^-$-edges in the cutset. As in (5.2), we assume that
$${\bf P}_{{\bf k}, m}[ e \mbox{ is an $\epsilon^-$ edge}] \leq F(\epsilon)-F(0)=\delta_1=\delta_1 (\epsilon),\eqno{(6.3)}$$
where $\delta_1\rightarrow 0$ as $\epsilon \rightarrow 0$.
With a small $\delta_1$, we estimate the following probability:
$${\bf P}_{{\bf k}, m}\left[ \bar{N}({\bf k}, m) \geq n, N^-({\bf k}, m)\geq -(D \log^{-1} (\delta_1))\bar{N}({\bf k}, m)\right].$$
Similar to  (5.4),  we need  to fix a vertex in ${\bf W}({\bf k}, m)\cap {\bf L}$.
Since ${\bf L}$ must intersect ${\bf W}({\bf k}, m)$, we select the intersection  
${\bf z}$ with the largest 
$x_{1}$-coordinate. There are at most $m$ choices for ${\bf z}$, since ${\bf W}({\bf k}, m)$ stays inside ${\bf B}({\bf k}, m)$.
By our assumption in Theorem 2,   note that   $n\geq \beta \|{\bf k}\|_v\geq \|{\bf k}\|_v$, 
$$m \leq \exp( \|{\bf k}\|_v)\leq \exp(n).\eqno{(6.4)}$$
 When ${\bf z}$ at ${\bf W}({\bf k}, m)$ is fixed, by Lemma 11 and the same
estimate in (5.4), we have
\begin{eqnarray*}
&&{\bf P}_{{\bf k}, m}\left[ \bar{N}({\bf k}, m) \geq n, N^-( {\bf k}, m)\geq -(D\log^{-1} (\delta_1))\bar{N}({\bf k}, m)\right ]\\
&&=\sum_{j=n}^{\infty}  {\bf P}_{{\bf k}, m}\left [ \bar{N}({\bf k}, m) =j, N^-({\bf k}, m)\geq - Dj
\log^{-1} (\delta_1)\right]\\
&&\leq \sum_{j=n}^{\infty} \exp( n)7^{d3^{d+1}j}2^{3^{d+1}j}(2d)^j\\
&&\hskip 2cm \times \max_{\Gamma} {\bf P}_{{\bf k}, m}\left[ |\Gamma|_e=j, 
\mbox{ $\Gamma$ contains more than $-Dj \log^{-1} (\delta_1)$ $\epsilon^-$-edges}\right]
\end{eqnarray*}
where $\Gamma$ is a fixed cutset that cuts ${\bf F}_0$ from ${\bf F}_m$ such that  the number of its edges
 $|\Gamma|_e=j$,  and
the maximum takes over all possible fixed vertex sets $\Gamma$. For a fixed set $\Gamma$, by the same estimate from (5.5)--(5.11), we have for a small $\delta_1>0$,
\begin{eqnarray*}
&&{\bf P}_{{\bf k}, m}\left[\Gamma, |\Gamma|_e=j, 
\mbox{ $\Gamma$ contains more than $- Dj\log^{-1} (\delta_1)$ $\epsilon^-$-edges}\right]\\
&&=\sum_{i\geq -Dj \log^{-1} (\delta_1)}^{j} {j \choose i} \delta_1^i \\
&&\leq \sum_{i\geq -Dj\log^{-1} (\delta_1)}^{j} \exp(Di\log (\delta_1)/2)\\
&&\leq 2 \exp(-Dj/2).
\end{eqnarray*}
With this observation,  by taking $D=D(d)$, 
there exists $\beta \geq 1$, and $\epsilon$, and $C_i=C_i(F, d, \epsilon)$ for $i=1,2$
such that
for all $n\geq \beta \|{\bf k}\|_v$,
\begin{eqnarray*}
&&{\bf P}_{{\bf k}, m}\left[ \bar{N}({\bf k}, m) \geq n, N^-( {\bf k}, m)\geq -D\log^{-1} (\delta_1)\bar{N}({\bf k}, m)\right ]\\
&&\leq \exp( n)\sum_{j=n}^{\infty}   7^{d3^{d+1}j}2^{3^{d+1}j}(2d)^j\exp(-Dj/2)\\
&&\leq C_1 \exp(-C_2 n).\hskip 4.5in (6.5)
\end{eqnarray*}

Therefore, 
for a small  $\epsilon$ and corresponding $\delta_1$, by (6.2) and (6.5), there exists $\beta=\beta(\epsilon)\geq \epsilon^{-2} \beta_1$
such that  for $n\geq \beta \|{\bf k}\|_v$, 
\begin{eqnarray*}
&&{\bf P}_{{\bf k}, m}\left[ \bar{N}({\bf k}, m) \geq n\right]\\
&&\leq  \sum_{j\geq n}{\bf P}_{{\bf k}, m}\left[ \bar{N}({\bf k}, m) =j, N^+({\bf k}, m) \leq  \epsilon j, N^-({\bf k}, m)\leq -Dj\log^{-1} (\delta_1) j\right]+  C_1 \exp(-C_2 n).
\end{eqnarray*}
Similarly, we denote by $J$ the number of all $\epsilon^\pm$-edges in ${\bf W}({\bf k}, m)$, and $\{e_1,\cdots, e_J\}$
are these $\epsilon^\pm $ edges.
On 
$$\left \{N^+({\bf k}, m) \leq  \epsilon j, N^-({\bf k}, m)\leq -Dj\log^{-1} (\delta_1) \right\},$$
 we have
$$J \leq \left(\epsilon -D\log^{-1} (\delta_1)\right) j.\eqno{(6.6)}$$
Therefore, for any large $t$, we take $\beta$ large such that for all $n\geq \beta \|{\bf k}\|_v$,
$$
{\bf P}_{{\bf k}, m}\left[ \bar{N}({\bf k}, m) \geq n\right]
\leq  \sum_{j\geq n}{\bf P}_{{\bf k}, m}\left[ \bar{N}({\bf k}, m) =j, J\leq j/ (2 (16t)^d)\right]+  C_1 \exp(-C_2 n).\eqno{(6.7)}
$$
For a configuration $\omega$,  since  $e_1,\cdots, e_J$ are the only non-zero edges in ${\bf W}({\bf k}, m)$, 
$$\tau_{\min}({\bf k}, m)=\sum_{i=1}^J \tau(e_i).\eqno{(6.8)}$$
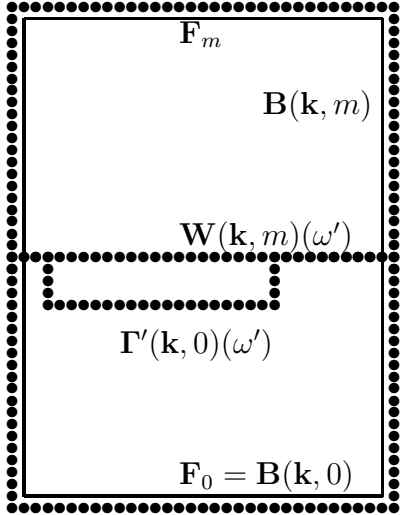
\begin{figure}\label{F:alphabeta}
\begin{center}
\setlength{\unitlength}{0.0125in}%
\begin{picture}(200,130)(67,800)
\thicklines
\put(100, 750){\line(0,1){200}}
\put(100, 745){\circle*{4}}
\put(105, 745){\circle*{4}}
\put(110, 745){\circle*{4}}
\put(115, 745){\circle*{4}}
\put(120, 745){\circle*{4}}
\put(125, 745){\circle*{4}}
\put(130, 745){\circle*{4}}
\put(135, 745){\circle*{4}}
\put(140, 745){\circle*{4}}
\put(145, 745){\circle*{4}}
\put(150, 745){\circle*{4}}
\put(155, 745){\circle*{4}}
\put(160, 745){\circle*{4}}
\put(165, 745){\circle*{4}}
\put(170, 745){\circle*{4}}
\put(175, 745){\circle*{4}}
\put(180, 745){\circle*{4}}
\put(185, 745){\circle*{4}}
\put(190, 745){\circle*{4}}
\put(195, 745){\circle*{4}}
\put(200, 745){\circle*{4}}
\put(205, 745){\circle*{4}}
\put(210, 745){\circle*{4}}
\put(215,745){\circle*{4}}
\put(220, 745){\circle*{4}}
\put(225, 745){\circle*{4}}
\put(230, 745){\circle*{4}}
\put(235, 745){\circle*{4}}
\put(240, 745){\circle*{4}}
\put(245, 745){\circle*{4}}
\put(250, 745){\circle*{4}}

\put(100, 955){\circle*{4}}
\put(105, 955){\circle*{4}}
\put(110, 955){\circle*{4}}
\put(115, 955){\circle*{4}}
\put(120, 955){\circle*{4}}
\put(125, 955){\circle*{4}}
\put(130, 955){\circle*{4}}
\put(135, 955){\circle*{4}}
\put(140, 955){\circle*{4}}
\put(145, 955){\circle*{4}}
\put(150, 955){\circle*{4}}
\put(155, 955){\circle*{4}}
\put(160, 955){\circle*{4}}
\put(165, 955){\circle*{4}}
\put(170, 955){\circle*{4}}
\put(175, 955){\circle*{4}}
\put(180, 955){\circle*{4}}
\put(185, 955){\circle*{4}}
\put(190, 955){\circle*{4}}
\put(195, 955){\circle*{4}}
\put(200, 955){\circle*{4}}
\put(205, 955){\circle*{4}}
\put(210, 955){\circle*{4}}
\put(215,955){\circle*{4}}
\put(220, 955){\circle*{4}}
\put(225, 955){\circle*{4}}
\put(230, 955){\circle*{4}}
\put(235, 955){\circle*{4}}
\put(240, 955){\circle*{4}}
\put(245, 955){\circle*{4}}
\put(250, 955){\circle*{4}}

\put(100, 750){\line(1,0){150}}
\put(100, 850){\circle*{4}}
\put(105, 850){\circle*{4}}
\put(110, 850){\circle*{4}}
\put(115, 850){\circle*{4}}
\put(120, 850){\circle*{4}}
\put(125, 850){\circle*{4}}
\put(130, 850){\circle*{4}}
\put(135, 850){\circle*{4}}
\put(140, 850){\circle*{4}}
\put(145, 850){\circle*{4}}
\put(150, 850){\circle*{4}}
\put(155, 850){\circle*{4}}
\put(160, 850){\circle*{4}}
\put(165, 850){\circle*{4}}
\put(170, 850){\circle*{4}}
\put(175, 850){\circle*{4}}
\put(180, 850){\circle*{4}}
\put(185, 850){\circle*{4}}
\put(190, 850){\circle*{4}}
\put(195, 850){\circle*{4}}
\put(200, 850){\circle*{4}}
\put(205, 850){\circle*{4}}
\put(210, 850){\circle*{4}}
\put(215, 850){\circle*{4}}
\put(220, 850){\circle*{4}}
\put(225, 850){\circle*{4}}
\put(230, 850){\circle*{4}}
\put(235, 850){\circle*{4}}
\put(240, 850){\circle*{4}}
\put(245, 850){\circle*{4}}
\put(250, 850){\circle*{4}}

\put(250, 750){\line(0,1){200}}
\put(95, 745){\circle*{4}}
\put(95, 750){\circle*{4}}
\put(95, 755){\circle*{4}}
\put(95, 760){\circle*{4}}
\put(95, 765){\circle*{4}}
\put(95, 770){\circle*{4}}
\put(95, 775){\circle*{4}}
\put(95, 780){\circle*{4}}
\put(95, 785){\circle*{4}}
\put(95, 790){\circle*{4}}
\put(95, 795){\circle*{4}}
\put(95, 800){\circle*{4}}
\put(95, 805){\circle*{4}}
\put(95, 810){\circle*{4}}
\put(95, 815){\circle*{4}}
\put(95, 820){\circle*{4}}
\put(95, 825){\circle*{4}}
\put(95, 830){\circle*{4}}
\put(95, 835){\circle*{4}}
\put(95, 840){\circle*{4}}
\put(95, 845){\circle*{4}}
\put(95, 850){\circle*{4}}
\put(95, 845){\circle*{4}}
\put(95, 850){\circle*{4}}
\put(95, 855){\circle*{4}}
\put(95, 860){\circle*{4}}
\put(95, 865){\circle*{4}}
\put(95, 870){\circle*{4}}
\put(95, 875){\circle*{4}}
\put(95, 880){\circle*{4}}
\put(95, 885){\circle*{4}}
\put(95, 890){\circle*{4}}
\put(95, 895){\circle*{4}}
\put(95, 900){\circle*{4}}
\put(95, 905){\circle*{4}}
\put(95, 910){\circle*{4}}
\put(95, 915){\circle*{4}}
\put(95, 920){\circle*{4}}
\put(95, 925){\circle*{4}}
\put(95, 930){\circle*{4}}
\put(95, 935){\circle*{4}}
\put(95, 940){\circle*{4}}
\put(95, 945){\circle*{4}}
\put(95, 950){\circle*{4}}
\put(95, 955){\circle*{4}}

\put(100, 950){\line(1,0){150}}

\put(255, 745){\circle*{4}}
\put(255, 750){\circle*{4}}
\put(255, 755){\circle*{4}}
\put(255, 760){\circle*{4}}
\put(255, 765){\circle*{4}}
\put(255, 770){\circle*{4}}
\put(255, 775){\circle*{4}}
\put(255, 780){\circle*{4}}
\put(255, 785){\circle*{4}}
\put(255, 790){\circle*{4}}
\put(255, 795){\circle*{4}}
\put(255, 800){\circle*{4}}
\put(255, 805){\circle*{4}}
\put(255, 810){\circle*{4}}
\put(255, 815){\circle*{4}}
\put(255, 820){\circle*{4}}
\put(255, 825){\circle*{4}}
\put(255, 830){\circle*{4}}
\put(255, 835){\circle*{4}}
\put(255, 840){\circle*{4}}
\put(255, 845){\circle*{4}}
\put(255, 850){\circle*{4}}
\put(255, 845){\circle*{4}}
\put(255, 850){\circle*{4}}
\put(255, 855){\circle*{4}}
\put(255, 860){\circle*{4}}
\put(255, 865){\circle*{4}}
\put(255, 870){\circle*{4}}
\put(255, 875){\circle*{4}}
\put(255, 880){\circle*{4}}
\put(255, 885){\circle*{4}}
\put(255, 890){\circle*{4}}
\put(255, 895){\circle*{4}}
\put(255, 900){\circle*{4}}
\put(255, 905){\circle*{4}}
\put(255, 910){\circle*{4}}
\put(255, 915){\circle*{4}}
\put(255, 920){\circle*{4}}
\put(255, 925){\circle*{4}}
\put(255, 930){\circle*{4}}
\put(255, 935){\circle*{4}}
\put(255, 940){\circle*{4}}
\put(255, 945){\circle*{4}}
\put(255, 950){\circle*{4}}
\put(255, 955){\circle*{4}}
\put(165,940){${\bf F}_m$}

\put(165,855){${\bf W}({\bf k}, m)(\omega')$}
\put(200,910){${\bf B}({\bf k}, m)$}
\put(165,755){${\bf F}_0={\bf B}({\bf k}, 0)$}

\put(100, 850){\circle*{4}}
\put(105, 850){\circle*{4}}
\put(110, 850){\circle*{4}}
\put(110, 845){\circle*{4}}
\put(110, 840){\circle*{4}}
\put(110, 835){\circle*{4}}
\put(110, 830){\circle*{4}}

\put(115, 830){\circle*{4}}
\put(120, 830){\circle*{4}}
\put(125, 830){\circle*{4}}
\put(130, 830){\circle*{4}}
\put(135, 830){\circle*{4}}
\put(140, 830){\circle*{4}}
\put(145, 830){\circle*{4}}
\put(150, 830){\circle*{4}}
\put(155, 830){\circle*{4}}
\put(160, 830){\circle*{4}}
\put(165, 830){\circle*{4}}
\put(170, 830){\circle*{4}}
\put(175, 830){\circle*{4}}
\put(180, 830){\circle*{4}}
\put(185, 830){\circle*{4}}
\put(190, 830){\circle*{4}}
\put(195, 830){\circle*{4}}
\put(200, 830){\circle*{4}}
\put(205, 830){\circle*{4}}
\put(205, 835){\circle*{4}}
\put(205, 840){\circle*{4}}
\put(205, 845){\circle*{4}}
\put(205, 850){\circle*{4}}

\put(210, 850){\circle*{4}}
\put(215, 850){\circle*{4}}
\put(220, 850){\circle*{4}}
\put(225, 850){\circle*{4}}
\put(230, 850){\circle*{4}}
\put(235, 850){\circle*{4}}
\put(240, 850){\circle*{4}}
\put(245, 850){\circle*{4}}
\put(250, 850){\circle*{4}}
\put(140,810){${\bf \Gamma}'({\bf k}, 0)(\omega')$}

\end{picture}
\end{center}
\vskip 1in
\caption{\em The surface edges of ${\bf B}({\bf k}, m)$, denoted by ${\bf \alpha }({\bf k}, m)$, are the dotted lines
 outside  ${\bf B}({\bf k}, m)$. We can use
the  surface together with ${\bf W}({\bf k}, m)(\omega')$ to construct a cutset that cuts ${\bf F}_0$ from $\infty$.
Thus, any open path from ${\bf F}_m$ to ${\bf F}_0$ must use an edge of ${\bf W}({\bf k}, m)$.
$\bar{\bf \Gamma}'({\bf k}, m)$ only uses the edges   inside  ${\bf B}({\bf k}, m)$.}
\end{figure}

To use the proof of Theorem 1, we need to construct a cutset that cuts  ${\bf B}({\bf k}, 0)={\bf F}_0$ from $\infty$.
We need to use  the surface edges $\alpha ({\bf k}, m)$ defined  in section 3 (see Fig. 3) and  ${\bf W}({\bf k}, m)$.
In particular, the surface edges  of $\alpha ({\bf k}, m)$ adjacent to ${\bf F}_m$ are called the {\em top} surface edges.
Moreover, let all the surface edges be closed. Note that the surface edges
are outside of ${\bf B}({\bf k}, m)$, so it will not affect our measure
${\bf P}_{{\bf k}, m}(\cdot)$. With the closed surface edges,
any path from ${\bf B}({\bf k}, 0)$ to $\infty$
must use at least  one surface edge. 
Thus, the closed surface  consists of a zero-cutset, so ${\cal G}({\bf k}, 0)$ occurs.
 Therefore, $\bar{\bf \Gamma}_t({\bf k}, 0)$ defined in section 2 exists and 
it contains a zero-cutset ${\bf \Gamma}({\bf k}, 0)$. 
Note that $\partial_e{\bf C}({\bf k}, 0)$ cannot be outside of the surface boundary,
so we may choose our ${\bf \Gamma}({\bf k}, 0)$ such that
$$\bar{\bf \Gamma}({\bf k}, 0) \mbox{ uses only edges of the surface and  the edges in ${\bf B}({\bf k}, m)$}.\eqno{(6.9)}$$
Let
$$\bar{\bf \Gamma}'({\bf k}, 0)=\bar{\bf \Gamma}({\bf k}, 0)\cap {\bf B}({\bf k}, m)\mbox{ and }
\bar{\bf \Gamma}'_t({\bf k}, 0)=\{ B_t({\bf u}): B_t({\bf u}) \cap \bar{\bf \Gamma}'({\bf k}, 0)\neq \emptyset\}.$$
For each configuration $\omega$, if we change all $e_1, e_2,\cdots, e_J$ from
$\epsilon^\pm$ to zero, we have another configuration $\omega'$.  With these changes, 
${\bf W}({\bf k}, m)(\omega')$ is a zero-cutset that cuts ${\bf F}_0$ from ${\bf F}_m$. Furthermore,  we will show that 
$$\bar{\bf \Gamma}'({\bf k}, 0)(\omega')\mbox{ is a closed set  inside }\bar{\bf \Gamma}'_t({\bf k}, 0)(\omega')\mbox{ that also 
cuts ${\bf F}_0$ from ${\bf F}_m$}.\eqno{(6.10)}$$ 
 Before showing (6.10), we first show that for $\omega'$, 
$\bar{\bf \Gamma}({\bf k}, 0)(\omega')$ is a zero-cutset that cuts ${\bf F}_0$ from ${\bf F}_m$. 
Intuitively, the surface edges in $\alpha({\bf k}, m)$ and the edge of ${\bf W}({\bf k}, m)(\omega')$ consist of a 
zero-cutset, so 
$\bar{\bf \Gamma}'({\bf k}, 0)(\omega')$ only uses the edges inside $\bar{\bf B}({\bf k}, m)(\omega')$ (see Fig. 3).
If $\bar{\bf \Gamma}({\bf k}, 0)(\omega')$ is not a cutset that cuts ${\bf F}_0$ from ${\bf F}_m$, 
then there exists a path (not necessarily open)
from ${\bf F}_0$ to ${\bf F}_m$ without using an edge of $\bar{\bf \Gamma}({\bf k}, 0)(\omega')$.
The path must reach  a vertex of  an edge, denoted by $e$, in the top surface. 
By (6.9),  $e\in \bar{\bf \Gamma}({\bf k}, 0)(\omega')$; otherwise, we can construct a path from ${\bf B}({\bf k}, 0)$
to $\infty$ without using  an edge of $\bar{\bf \Gamma}({\bf k}, 0)(\omega')$.
In other words, it reaches an edge in $\bar{\bf \Gamma}({\bf k}, 0)(\omega')$ and the edge is also adjacent to ${\bf F}_m$
from outside of ${\bf B}({\bf k}, m)$. Let  
$B_t({\bf u})\in \bar{\bf \Gamma}_t({\bf k}, m)(\omega')$ be the $t$-cube that contains the edge. Since
$k_i/t$ and $m/t$ are integers, $B_t({\bf u})$ and ${\bf B}({\bf k}, m)$ do not have  other  vertices in common, 
except  for vertices at ${\bf F}_m$.
 By Lemma 4, there exists an open path from 
$B_t({\bf u})$ to ${\bf B}({\bf k}, 0)={\bf F}_0$. Note that the surface is closed, so the open path must
 go from ${\bf F}_m$ to ${\bf F}_0$ inside ${\bf B}({\bf k}, m)$ (see Fig. 3). However, this situation contradicts
 the fact that
 ${\bf W}({\bf k}, m)(\omega')$ is a zero-cutset.
This contradiction shows that $\bar{\bf \Gamma}({\bf k}, 0)(\omega')$ is indeed  a cutset that cuts ${\bf F}_0$ from ${\bf F}_m$. Furthermore, note that $\omega'$ has more zeros than $\omega$'s and all edges in 
$\bar{\bf \Gamma}({\bf k}, 0)(\omega)$ are all zero-edges, so $\bar{\bf \Gamma}({\bf k}, 0)(\omega')$ is a zero-cutset that cuts ${\bf F}_0$ from ${\bf F}_m$. 
Note that the edges of $\bar{\bf \Gamma}({\bf k}, 0)(\omega')$ outside of ${\bf B}({\bf k}, m)$ will not affect
whether or not $\bar{\bf \Gamma}({\bf k}, 0)(\omega')$ cuts ${\bf F}_0$ from ${\bf F}_m$ inside ${\bf B}({\bf k}, m)$,
so (6.10) follows. In addition, $\bar{\bf \Gamma}'({\bf k}, 0)(\omega')$  is a self-avoiding cutset, since 
$\bar{\bf \Gamma}({\bf k}, 0)(\omega')$ is self-avoiding. 

With (6.10), the remaining proof is similar to he proof in Theorem 1.
If we change $\omega'$ back to $\omega$, $e_i$ changes from zero back to original values.
$\bar{\bf \Gamma}'({\bf k}, 0)$, as a vertex set, exists.  But 
$\bar{\bf \Gamma}'({\bf k}, 0)$ will no longer be a zero-cutset. We denote by $\bar{\bf \Gamma}'({\bf k}, 0)(\omega)$
as the set corresponding to configuration $\omega$. Note that the other edges except for
$e_i$ are all zero-edges in both $\bar{\bf \Gamma}'({\bf k}, 0)(\omega)$ and ${\bf W}({\bf k}, m)(\omega)$, so
$$\tau(\bar{\bf \Gamma}'({\bf k}, 0)(\omega))=\tau({\bf W}({\bf k}, m)(\omega))=\tau_{\min}({\bf k}, 0)(\omega).\eqno{(6.11)}$$
Therefore, for each $\omega$,
$$|\bar{\bf \Gamma}'({\bf k}, 0)(\omega')|_e=|\bar{\bf \Gamma}'({\bf k}, 0)(\omega)|_e\geq \bar{N}({\bf k}, m)(\omega).\eqno{(6.12)}$$

For each $\omega$, we focus on $\omega'$. As we mentioned above, $\bar{\bf \Gamma}'({\bf k}, 0)(\omega')$ 
is a self-avoiding zero-cutset contained inside $\bar{\bf \Gamma}_t'({\bf k}, 0)(\omega')$.
Note that ${\bf L}$, defined as the line below (5.2), must intersect $\bar{\bf \Gamma}'({\bf k}, 0)(\omega')$ 
inside $\bar{\bf B}({\bf k}, m)(\omega')$, otherwise
$\bar{\bf \Gamma}'({\bf k}, 0)(\omega')$  will not be a cutset. We denote by ${\bf z}$ the intersection vertex. If it is not unique, we select the 
one with the largest $x_1$-coordinate. 
Thus, there are at most $m$ choices in $\bar{\bf \Gamma}'_t({\bf k}, 0)(\omega')$ for the cube that contains ${\bf z}$,
since $\bar{\bf \Gamma}'_t({\bf k}, 0)(\omega')$ must stay inside ${\bf B}({\bf k}, m)$. 
As we discussed in the proof of Theorem 1, if the number of cubes of
$\bar{\bf \Gamma}_t'({\bf k}, 0)(\omega')$ is $s$, then
\begin{eqnarray*}
&&\exists \mbox{ at least $s/ 2^{4d}$ disjoint $3t$-cubes such that their center cubes with block proerty}\\
&&\mbox{ in  $\bar{\bf \Gamma}_t'({\bf k}, 0)(\omega')$}.\hskip 12cm {(6.13)}
\end{eqnarray*}

Let ${\bf M}_{3t}({\bf k},m)(\omega')$ be all $3t$-cubes and 
${ M}_{3t}({\bf k},m)(\omega')$ be the number of these
$3t$-cubes in ${\bf M}_{3t}({\bf k},m)(\omega')$.
 Note that each $t$-cube has $2t^d$ edges, so by (6.12), if
$\bar{N}({\bf k}, m)(\omega)=j$, by (6.12)
$${M}_{3t}({\bf k},m)(\omega')\geq {\# \mbox{ cubes in }\bar{\bf \Gamma}_t({\bf k}, 0)(\omega')\over 2^{4d}}\geq {|\bar{\bf \Gamma}({\bf k}, 0)(\omega')|_e\over 2^{4d+1} t^d }
\geq {\bar{N}({\bf k}, m)(\omega)\over 2(16t)^d}={j\over 2(16t)^d} .\eqno{(6.14)}$$

Furthermore, if
$$J(\omega) \leq j /(2(16t)^d)\mbox{ and }\bar{N}({\bf k}, m)(\omega)=j \mbox{ and the number of cubes in } \bar{\bf \Gamma}_t'({\bf k}, 0)(\omega')\mbox{ is }s,$$
by (6.14), for each $\omega'$, there are at least
$${ M}_{3t}({\bf k},m)-j/4(16t)^d= { M}_{3t}({\bf k},m)/2+ { M}_{3t}({\bf k},m)/2-j /(4 (16t)^d)\geq { M}_{3t}({\bf k},m)/2\geq  s/2^{4d+1} $$
  center cubes in $\bar{\bf \Gamma}_t'({\bf k}, 0)(\omega')$ with  the blocked  property
and they do not contain  $e_1,\cdots, e_J$ in their interiors.  
Recall that they are called  the permanent blocked  cubes.
 Now we change these edges in $\{e_1,\cdots, e_J\}$ from zero back to the original values. We still have 
$s/2^{2d+1}$ permanent blocked 
or disjoint $t$-cubes. 
Also, by (6.12),
$$\# \mbox { cubes in }\bar{\bf \Gamma}_t({\bf k}, 0)(\omega')=s\geq |\bar{\bf \Gamma}({\bf k}, 0)(\omega')|_e/(2t)^d\geq \bar{N}({\bf k}, m)(\omega)/ (2t)^d =j /(2t)^d. \eqno{(6.15)}$$
Finally, by Lemma 11, $\bar{\bf \Gamma}'({\bf k}, 0)$ is ${\bf Z}^d$-connected, so
 $$\bar{\bf \Gamma}'_t({\bf k}, 0)(\omega')\mbox{  is ${\bf L}^d$-connected.}\eqno{(6.16)}$$

In summary, for each $\omega$, if $N({\bf k}, m)(\omega)=j$ and $J(\omega)\leq j /(2(16t)^d)$,
then there are $s\geq j/ (2t)^d$ and $\|{\bf z}\|\leq m\leq \exp(\|{\bf k}\|_v)$ such that
${\cal E}(s,j,{\bf z})$ occurs, where   ${\cal E}(s,j,{\bf z})$ is the event defined in section 5
after (5.25). 
Therefore, by the same estimate as (5.26), there are $C_i=C_i(F, d, \epsilon, t)$ for $i=1,2$ and $C_i=C_i(F(0), d)$
for $i=3,4$ such that
\begin{eqnarray*}
&&{\bf P}_{{\bf k}, m}\left[ \bar{N}({\bf k}, m) \geq n\right]\\
&&\leq \sum_{j\geq n}\sum_{s\geq {j\over (2t)^d}} \exp(\|{\bf k}\|_v) 7^{2ds}2^s \left[ C_3\exp(-C_4t)\right]^{s/2^{4d+1}}+C_1\exp(-C_2n)\\
&&\leq \sum_{j\geq n}\sum_{s\geq {j\over (2t)^d}} \exp(n/\beta) 7^{2ds}2^s \left[ C_3\exp(-C_4t)\right]^{s/2^{4d+1}}+C_1\exp(-C_2n).
\end{eqnarray*}
If we take $t$ large  and  $\beta$ large, there exist $C_i=C_i(F, d, t, \epsilon,\beta)$  for $i=1,2$
such that for all $n\geq \beta \|{\bf k}\|_v$,
$${\bf P}_{{\bf k}, m}[ \bar{N}({\bf k}, m) \geq n]\leq C_1 \exp(-C_2 n).\eqno{(6.17)}$$
Therefore, Theorem 2 follows.

\section{Patching cutsets. }
Given a cutset ${\bf W}({\bf k}, m)$ as we defined in section 1, we shall now discuss a few basic properties of 
this cutset. Let ${\bf k}'=(k_1'\cdots, k_{d-1}')$ and ${\bf k}=(k_1,\cdots, k_{d-1})$ be two vectors. We say 
$${\bf k}'\leq {\bf k} \mbox{ if }0\leq  k_i'\leq k_i \mbox { for all }i=1,\cdots, d-1.$$
We also denote by ${\bf F}_0'$ and ${\bf F}_m'$  the bottom and the
top faces of the box ${\bf B}({\bf k}', m)$.
With these definitions we have the following lemma.\\

{\bf Lemma 12.} 
{\em If ${\bf k}'\leq {\bf k}$, then}
\begin{eqnarray*}
&&\mbox{ (a)} {\bf W}({\bf k}, m)\cap {\bf B}({\bf k}', m)\mbox{  is a cutset that cuts ${\bf F}_0'$ from ${\bf F}_m'$ in }
{\bf B}({\bf k}', m),\\
&&\mbox{ (b) }\tau({\bf W}({\bf k}', m) )\leq \tau ({\bf W}({\bf k}, m)) ,\\
&&\mbox{ (c) } \tau ({\bf W}({\bf k},  m))\leq \tau ({\bf W}({\bf k'}, m))+ \sum_{e\in {\bf B}({\bf k}, m)\setminus {\bf B}({\bf k}', m)}\tau (e).
\end{eqnarray*}

{\bf Proof.} 
To prove (a),
 we only need to show that any path in ${\bf B}({\bf k}', m)$ 
from ${\bf F}_0'$ to ${\bf F}_m'$ must use at least an edge of ${\bf W}({\bf k}, m)\cap {\bf B}({\bf k}', m)$.
Note that such a path is also a path from ${\bf F}_0$ to ${\bf F}_m$ in ${\bf B}({\bf k}, m)$ and note also
that ${\bf W}({\bf k}, m)$ is a cutset, so
any such path must use at least one edge of ${\bf W}({\bf k}, m)$. On the other hand, any such path
must stay in ${\bf B}({\bf k}', m)$, so it must use at least one edge of ${\bf W}({\bf k}, m)\cap {\bf B}({\bf k}', m)$. Therefore, (a) follows.
With (a), (b) follows from the definitions of ${\bf W}({\bf k}', m) )$ and ${\bf W}({\bf k}, m) $ directly.

Now we show (c). By the same argument as (a), we can show ${\bf W}({\bf k}', m)\cup [{\bf B}({\bf k}, m)\setminus {\bf B}({\bf k}', m)]$ is a cutset for ${\bf B}({\bf k}, m)$, so  (c) follows.
$\Box$\\

Now we want to patch two smaller  cutsets into a larger cutset. To do it, we need to study the
traces of the cutset in the boundary of the box ${\bf B}({\bf k}, m)$.
We denote the hyperplane  by
$${\bf L}_{n}=\{(x_1,\cdots, x_d): x_1=n\}.$$
For a cutset ${\bf W}({\bf k}, m)$, we define its {\em trace} in the hyperplane ${ \bf L}_{k_1}$ by
$${\bf I}({\bf k}, m)={\bf W}({\bf k}, m)\cap {\bf L}_{k_1}.$$
Let edges in ${\bf I}({\bf k}, m)$ be ${\bf I}_e({\bf k}, m)$.
If we  remove  all the edges of ${\bf I}_e({\bf k}, m)$ from ${\bf L}_{k_1}$, but leave the vertices of these
edges, the new graph, after removing
these edges,  consists
 of several clusters on ${ \bf L}_{k_1}$. Note that there might be a few clusters with only one isolated vertex.
We now analyze these clusters on the hyperplane. 

\begin{figure}\label{F:alphabeta}
\begin{center}
\setlength{\unitlength}{0.0125in}%
\begin{picture}(200,120)(67,820)
\thicklines
\put(0, 800){\line(0,1){200}}
\put(100, 750){\line(0,1){200}}
\put(0, 900){\line(2,-1){100}}
\put(0, 1000){\line(2,-1){100}}
\put(0, 800){\line(2,-1){100}}
\put(0, 900){\line(1,0){150}}
\put(150, 900){\line(2,-1){100}}

\put(100, 750){\line(1,0){150}}
\put(0, 1000){\line(1,0){150}}
\put(100, 850){\line(1,0){150}}

\put(250, 750){\line(0,1){200}}

\put(100, 950){\line(1,0){150}}
\put(150, 1000){\line(2,-1){100}}
\put(50, 920){\oval(20,30)}
\put(80, 880){\oval(20,15)}

\put(50,935){\circle*{2}}
\put(55,935){\circle*{2}}
\put(60,933){\circle*{2}}
\put(65,930){\circle*{2}}
\put(70,925){\circle*{2}}
\put(75,920){\circle*{2}}
\put(81,910){\circle*{2}}

\put(78,915){\circle*{2}}
\put(85,905){\circle*{2}}
\put(87,900){\circle*{2}}
\put(88,895){\circle*{2}}
\put(88,890){\circle*{2}}
\put(90,885){\circle*{2}}

\put(50,905){\circle*{2}}
\put(55,902){\circle*{2}}
\put(60,898){\circle*{2}}
\put(65,895){\circle*{2}}
\put(68,890){\circle*{2}}
\put(70,882){\circle*{2}}
\put(45,915){${\bf S}'$}
\put(15,915){${\bf T}$}

\put(165,855){${\bf W}({\bf k}, m)$}
\put(65,825){${\bf T}'$}
\put(15,825){${\bf S}$}

\put(65,745){${\bf L}_{k_1}$}

\put(70, 830){\oval(20,30)}
\put(130, 870){\oval(10,10)}
\put(70,845){\circle*{2}}
\put(75,845){\circle*{2}}
\put(80,845){\circle*{2}}
\put(85,845){\circle*{2}}
\put(70,845){\circle*{2}}
\put(75,845){\circle*{2}}
\put(80,845){\circle*{2}}
\put(85,845){\circle*{2}}
\put(90,847){\circle*{2}}
\put(95,849){\circle*{2}}
\put(100,851){\circle*{2}}
\put(110,853){\circle*{2}}
\put(115,857){\circle*{2}}
\put(120,861){\circle*{2}}
\put(125,865){\circle*{2}}

\put(70,815){\circle*{2}}
\put(75,815){\circle*{2}}
\put(80,815){\circle*{2}}
\put(85,815){\circle*{2}}
\put(70,815){\circle*{2}}
\put(75,815){\circle*{2}}
\put(80,815){\circle*{2}}
\put(85,815){\circle*{2}}
\put(90,817){\circle*{2}}
\put(95,819){\circle*{2}}
\put(100,821){\circle*{2}}
\put(110,823){\circle*{2}}
\put(115,827){\circle*{2}}
\put(120,831){\circle*{2}}
\put(125,835){\circle*{2}}
\put(128,840){\circle*{2}}
\put(131,845){\circle*{2}}
\put(135,850){\circle*{2}}
\put(137,855){\circle*{2}}
\put(137,860){\circle*{2}}
\put(137,865){\circle*{2}}

\end{picture}
\end{center}
\vskip 1in
\caption{\em This graph shows the exits of  upper tunnels  and lower tunnels on hyperplane ${\bf L}_{k_1}$. 
The middle plane is the cutset ${\bf W}({\bf k}, m)$. There are four exits of tunnels. 
${\bf S}$ below and ${\bf T}$ above
the cutset are two trivial exits. The circled ${\bf S}'$ and ${\bf T}'$, above and below the cutset, respectively,
are the other two exits.
One  can use the  exits of the circled  tunnels ${\bf S}'$ or ${\bf T}'$
 from ${\bf F}_m$ or ${\bf F}_0$ to ${\bf T}'$ or ${\bf S}'$ to cross the middle surface without using its edges. }

\end{figure}
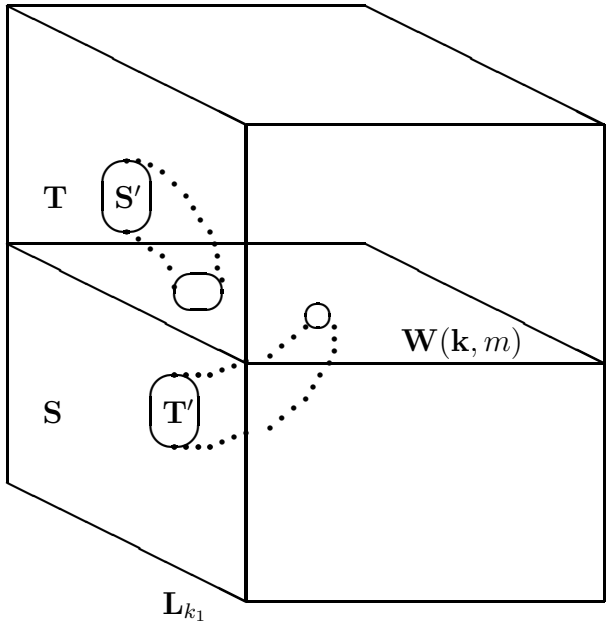
We denote by ${\bf T}_{k_1, 1}({\bf k}, m),\cdots, {\bf T}_{k_1,t}({\bf k}, m)\subset {\bf L}_{k_1}\setminus {\bf I}_e({\bf k}, m)$ (see Fig. 4) the first  kind of clusters such that 
each of their vertices is connected to ${\bf F}_m$ by a path lying in ${\bf B}({\bf k}, m)$
without using any edge of ${\bf W}({\bf k}, m)$ (see Fig. 4).
Note that $t\geq 1$ since 
$${\bf L}_{k_1}\cap {\bf F}_m\neq \emptyset. $$
Here we may view ${\bf T}$ as  both a vertex and an edge set.

We also denote by ${\bf S}_{k_1,1}({\bf k}, m),\cdots {\bf S}_{k_1,s}({\bf k}, m)\subset {\bf L}_{k_1}\setminus {\bf I}_e({\bf k}, m)$ (see Fig. 4)
as the second kind of clusters 
such that 
each of their vertices is connected to ${\bf F}_0$ by a path lying in ${\bf B}({\bf k}, m)$
without using any edge of ${\bf W}({\bf k}, m)$.
Similarly, we have $s\geq 1$. We write these ${\bf T}$ and ${\bf S}$ for the {\em exits of upper tunnels } and {\em exits of lower tunnels}, respectively.
If we do not work on a specific box, we may just write ${\bf T}_{k_1,j}$ and ${\bf S}_{k_1,i}$ rather than ${\bf T}_{k_1,j}({\bf k}, m)$
and ${\bf S}_{k_1,i}({\bf k}, m)$ as the exits of the upper and the lower tunnels. 
With these definitions, we have the following lemma.\\

{\bf Lemma 14.} {\em For all configurations, }
$$ {\bf S}_{k_1,i}({\bf k}, m)\cap {\bf T}_{k_1,j}({\bf k}, m)=\emptyset\mbox{ for  all }i=1,\cdots, s\mbox{ and }j=1,\cdots, t.$$

{\bf Proof.} If there exists a common
vertex belonging to ${\bf S}_{k_1,i}\cap {\bf T}_{k_1,j}$ for some $i$ and $j$, then there exist paths from ${\bf F}_0$ to ${\bf F}_m$
in ${\bf B}({\bf k}, m)$ 
without using an edge of the cutset ${\bf W}({\bf k}, m)$.  This  contradicts the definition of
${\bf W}({\bf k}, m)$. $\Box$\\

Also by our definitions of ${\bf S}$ and ${\bf T}$, we have the following lemma.\\

{\bf Lemma 15.}
 {\em (a) The boundary edges of ${\bf T}_{k_1,j}$ and ${\bf S}_{k_1,i}$, on 
${\bf L}_{k_1}\cap {\bf B}({\bf k}, m)$, belong to ${\bf W}({\bf k}, m)$.\\
(b) For any vertex $v\in {\bf L}_{k_1}\cap {\bf B}({\bf k}, m)$, if there exists a path from $v$ to  ${\bf S}_{k_1,i}$  (to ${\bf T}_{k_1,j}$)
in ${\bf B}({\bf k}, m)$ without using any edge of
${\bf W}({\bf k}, m)$, then $v\in {\bf S}_{k_1,i}$  ($v\in {\bf T}_{k_1,j}$).\\
(c) ${\bf W}({\bf k}, m)$, ${\bf T}_{k_1,j}$ and ${\bf S}_{k_1,i}$ only depend on the configurations of edges in 
${\bf B}({\bf k}, m)$.}\\
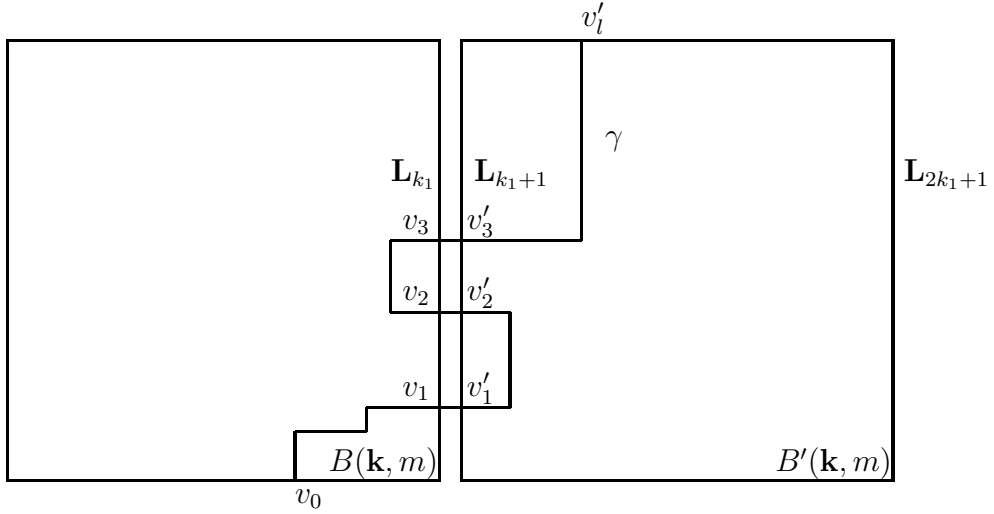
\begin{figure}\label{F:alphabeta}
\begin{center}
\setlength{\unitlength}{0.0125in}%
\begin{picture}(200,120)(67,800)
\thicklines
\put(100, 760){\line(0,1){20}}
\put(100, 780){\line (1,0){30}}
\put(130, 780){\line (0,1){10}}
\put(130, 790){\line (1,0){60}}
\put(100,750){\mbox{$v_0$}}
\put(145,795){\mbox{$v_1$}}
\put(172,795){\mbox{$v_1'$}}
\put(190, 790){\line (0,1){40}}
\put(190, 830){\line (-1,0){50}}
\put(140, 830){\line (0,1){30}}
\put(145,835){\mbox{$v_2$}}
\put(172,835){\mbox{$v_2'$}}
\put(140, 860){\line (1,0){80}}
\put(220, 860){\line (0,1){84}}
\put(145,865){\mbox{$v_3$}}
\put(172,865){\mbox{$v_3'$}}
\put(220,950){\mbox{$v_l'$}}
\put(230,900){\mbox{$\gamma$}}
\put(170,760){\framebox(180,180)[br]{\mbox{$B'({\bf k}, m)$}}}
\put(-20,760){\framebox(180,180)[br]{\mbox{$B({\bf k},m)$}}}
\put(140,885){\mbox{${\bf L}_{k_1}$}}
\put(175,885){\mbox{${\bf L}_{k_1+1}$}}
\put(355,885){\mbox{${\bf L}_{2k_1+1}$}}

\end{picture}
\end{center}
\vskip 1in
\caption{ \em This cross graph shows that if all  the upper and lower tunnels of  ${\bf B}({\bf k}, m)$ and 
 of ${\bf B}'({\bf k}, m)$, respectively,  in ${\bf L}_{k_1}$ are the same, and if there is a path from
${\bf F}_0\cup {\bf F}_0'$ to ${\bf F}_m\cup {\bf F}'_m$ without using ${\bf W}({\bf k},m)\cup {\bf W}'({\bf k}, m)$,
then the path, from ${v_1}_1$ to $v_l'$, must go around the tunnels. Hence, the result is that the exits of the lower and upper tunnels for
${\bf W}({\bf k}, m)$ have a common vertex, which is a contradiction.}
\end{figure}

With these observations, we are ready to patch two cutsets on two adjacent boxes. 
Before we state our result, we would like introduce more
definitions (see Fig. 5). 
Let ${\bf W}({\bf k}, m)$ and ${\bf W}'({\bf k},m)$ be two cutsets from ${\bf F}_0$ to ${\bf F}_m$
and from ${\bf F}_0'$ to ${\bf F}_m'$ of ${\bf B}({\bf k}, m)$
and ${\bf B}'({\bf k}, m)$, respectively, where
$${\bf B}'({\bf k}, m)=[k_1+1, 2k_1+1]\times [0,k_2]\times \cdots \times [0, k_{d-1}]\times [0,m]$$
and ${\bf F}_0'$ and ${\bf F}_m'$ are the bottom and top faces of ${\bf B}'({\bf k}, m)$.
Here we select ${\bf W}'({\bf k},m)$ to be a self-avoiding cutset with the minimal passage time in ${\bf B}'({\bf k}, m)$.
If there is more than one such cutset, we simply select one with the unique method, but it may not necessarily 
follow the same rule as the selection for ${\bf W}({\bf k},m)$.

Similarly, we define the hyperplane of ${\bf B}'({\bf k}, m)$
next to ${\bf L}_{k_1}$ by
$${\bf L}_{k_1+1}=\{(x_1,\cdots, x_d\}: x_1=k_1+1\}.$$
Let ${\bf I}'_e({\bf k}, m)$ be the edge set with vertices in
$${\bf W}'({\bf k}, m)\cap {\bf L}_{k_1+1}.$$
We denote 
$${\bf B}({\bf k}, m)\cup {\bf B}'({\bf k}, m)=[0,2k_1+1]\times [0,k_2]\times \cdots \times [0,k_{d-1}]\times [0,m]$$
and use ${\bf F}_0\cup {\bf F}_0'$ and ${\bf F}_m\cup {\bf F}_m'$ to
denote its bottom and top faces.

Similarly, let  ${\bf T}'_{k_1+1, j} $ and ${\bf S}_{k_1+1,i}'$ 
be the exits of  the upper or lower tunnels of ${\bf W}'({\bf k}, m)$ 
on the hyperplanes ${\bf L}_{k_1+1}$. 
We say ${\bf T}$ is shifted  $l$ units invariance as ${\bf T}'$ if 
$${\bf T}'=\{(u_1+l,u_2, \cdots, u_d): (u_1,\cdots, u_d)\in {\bf T}\}.$$
We write ${\bf T}\stackrel{l}{=}{\bf T}'$ for the above ${\bf T}$ and ${\bf T}'$.\\

{\bf Lemma 16.} {\em Let $\{{\bf T}_{k_1,1}, {\bf T}_{k_1, 2},\cdots, {\bf T}_{k_1, t}\}$ and 
$\{{\bf S}_{k_1,1}, {\bf S}_{k_1, 2},\cdots, {\bf S}_{k_1, s}\}$ be the exits of the upper and the lower tunnels in ${\bf L}_{k_1}$ for  ${\bf W}({\bf k}, m)$. Let  $\{{\bf T}_{k_1,1}', {\bf T}_{k_1, 2}',\cdots, {\bf T}_{k_1, t}'\}$ and 
$\{{\bf S}_{k_1,1}', {\bf S}_{k_1, 2}',\cdots, {\bf S}_{k_1, s}'\}$ be the exits of the upper and the lower tunnels in ${\bf L}_{k_1+1}$ for  ${\bf W}'({\bf k}, m)$. If
 ${\bf T}_{k_1,j}\stackrel{1}{=}{\bf T}_{k_1+1,j}'$ and ${\bf S}_{k_1,i}\stackrel{1}{=}{\bf S}_{k_1,i}'$ for all $i$ and $j$, then
${\bf W}({\bf k}, m)\cup {\bf W}({\bf k}'m)$ is a cutset that cuts from ${\bf F}_0\cup {\bf F}_0'$ to ${\bf F}_m\cup {\bf F}_m'$
in the box ${\bf B}({\bf k}, m)\cup {\bf B}'({\bf k}, m)$}.\\

{\bf Proof.} Under the hypotheses of Lemma 16,  we suppose that there is a path $\gamma$ in
${\bf B}({\bf k}, m)\cup {\bf B}'({\bf k},m)$ from ${\bf F}_0\cup {\bf F}_0'$ to ${\bf F}_m\cup {\bf F}_m'$ without using 
any edge of ${\bf W}({\bf k}, m)\cup{\bf W}({\bf k}'m)$. Since ${\bf W}({\bf k}, m)$ and ${\bf W}({\bf k}'m)$
are cutsets of ${\bf B}({\bf k}, m)$ and ${\bf B}'({\bf k}, m)$,   $\gamma$ cannot lie in 
${\bf B}({\bf k}, m)$ or in ${\bf B}'({\bf k}, m)$, respectively. 
The path $\gamma$ should  be  a snake-shaped between two boxes 
${\bf B}({\bf k}, m)$ and ${\bf B}'({\bf k}, m)$ (see Fig. 5). 
We then go along $\gamma$ from ${\bf F}_0\cap {\bf F}_0'$ to ${\bf F}_m\cup {\bf F}_m'$.
Without loss of generality, we assume that $\gamma$ starts at ${\bf F}_0$ and ends at ${\bf F}_m'$. With this definition, 
$\gamma$ must go out of the hyperplane ${\bf L}_{k_1}$. Let $v_1$ be the first vertex that $\gamma$ exits from
${\bf L}_{k_1}$. After that,  $\gamma$ must go through   ${\bf L}_{k_1+1}$ at $v'_1$. Let $e_{v_1, v'_1}$  be the edge
with vertices $v_1$ and $v'_1$ (see Fig. 5). Note that $e_{v_1, v'_1}$ is neither in ${\bf B}({\bf k}, m)$ nor 
in ${\bf B}({\bf k}', m)$, but  just between these two boxes. We then continue following $\gamma$ from $v_1'$. 
If it can reach ${\bf F}_m'$ inside ${\bf B}({\bf k}', m)$ directly, then we stop our trip. If  it cannot, 
let $v_2'$ be the vertex in ${\bf L}_{k_1+1}$ that $\gamma$ first goes out of ${\bf B}'({\bf k}, m)$.  Similarly,
we will have the vertex $v_2\in \gamma\cap {\bf L}_{k_1}$ such that $\gamma$ first goes back at $v_2$
from ${\bf B}'({\bf k}, m)$. Let $e_{v_2', v_2}$ be the edge with the vertices $v_2$ and $v_2'$ between these two boxes.
We continue this process until $\gamma$ reaches  ${\bf F}_m'$.
Let $v_0\in {\bf F}_0$ and $v_l'\in {\bf F}_m'$ be the starting vertex and the ending vertex of $\gamma$.
Our $\gamma$ contains the following  vertices and edges between ${\bf B}({\bf k}, m)$ and ${\bf B}'({\bf k}, m)$:
$$v_0, v_1, e_{v_1, v_1'}, v_1', v_2',e_{v_2', v_2},v_2, v_3, e_{v_3, v_3'}, v_3',\cdots, v_l'.$$
Note that $\gamma$ never uses an edge of ${\bf W}({\bf k}, m)$,
so $v_1\in \cup_j {\bf S}_{k_1,i}$. By the assumption of Lemma 16, 
$v_1'\in \cup_j {\bf S}_{k_1+1,i}'$. By Lemma 15 (b),
$v_{2}'\in \cup_j {\bf S}_{k_1+1,i}'$ 
so $v_2\in \cup_j {\bf S}_{k_1,i}$. If we iterate this way, we finally have  $v_l'\in \cup_i {\bf S}_{k_1+1,i}'$. However, 
 note that $\gamma$ never uses an edge of ${\bf W}({\bf k}, m)$, so
$v_l\in \cup_j {\bf T}_{k_1+1,j}'$. Therefore, this result shows that $\cup_i {\bf S}_{k_1,i}'$ and $\cup_j {\bf T}'_{k_1,j}$ have a common vertex, but it 
contradicts  Lemma 14, so Lemma 16 follows. $\Box$\\

\section{ Estimates for the  boundary size of a cutset.}

A cutset   ${\bf W}_r({\bf k}, m)$ in
${\bf B}({\bf k}, m)$ cutting ${\bf F}_0$ from ${\bf F}_m$  is said to be {\em regular} if
$$|{\bf W}_r({\bf k}, m)|\leq \bar{\beta} \|{\bf k}\|_v, \eqno{(8.0)}$$
where $\bar{\beta}=2d \beta$ for the $\beta $  defined in Theorem 2.
We select a regular cutset, still denoted  by ${\bf W}_r({\bf k}, m)$, with the minimum passage time.
We may also denote $\tau({\bf W}_r({\bf k}, m))=\tau_r({\bf k}, m)$. 
If ${\bf W}_r({\bf k}, m)$ is not unique, we select it with the minimum number of edges using the unique method of selection.
In particular, if  ${\bf W}({\bf k}, m)$, defined in section 1,
 satisfies (8.0), we only select 
$${\bf W}_r({\bf k}, m)={\bf W}({\bf k}, m).$$
Clearly,
$$\tau_{\min} ({\bf k}, m)\leq \tau_r ({\bf k},m).$$
If $|{\bf W}({\bf k}, m)|\geq \bar{\beta} \|{\bf k}\|_v$, then 
$$|{\bf W}({\bf k}, m)|_e\geq {\beta} \|{\bf k}\|_v.$$
By this observation and  Theorem 2, there exist $C_i=C_i(\beta, F)$ such that
$${\bf P}[ {\bf  W}({\bf k}, m)\neq {\bf W}_r({\bf k}, m)]\leq {\bf P}[N({\bf k}, m)\geq \beta \|{\bf k}\|_v]
\leq C_1 \exp(-C_2 \|{\bf k}\|_v).\eqno{(8.1)}$$
Now we only focus on regular cutsets. 
Under (8.0), note that there are $k_1$ disjoint  hyperplanes 
in ${\bf B}({\bf k}, m)$
perpendicular to the first coordinate,
so the number of vertices of ${\bf W}_r({\bf k}, m)$ on a few of these hyperplanes should be much less that
$\bar{\beta} \|{\bf k}\|_v$. Recall that ${\bf L}_i$ is defined in section 7 as  the hyperplane of $\{x_1=i\}$.
Now we try to find two such hyperplanes.
We account for the size of  $\{{\bf W}_r({\bf k}, m)\cap {\bf L}_{k_1}\}\cup \{{\bf W}_r({\bf k}, m)\cap {\bf L}_{0}\}$ 
to see whether
$$|\{{\bf W}_r({\bf k}, m)\cap {\bf L}_{k_1}\}\cup \{{\bf W}_r({\bf k}, m)\cap {\bf L}_{0}\}|\leq 
\bar{\beta }k_1^{\delta/2} k_2\cdots k_{d-1}, \eqno{(8.2)}$$
where $\delta$ is defined in (1.11).
If the cutset satisfies (8.2), we select ${\bf L}_0$ and ${\bf L}_{k_1}$. If it does not, we account for the size of 
${\bf W}_r({\bf k}, m)\cap {\bf L}_{1}$ and ${\bf W}_r({\bf k}, m)\cap {\bf L}_{k_1-1}$ to see whether 
$$|\{{\bf W}_r({\bf k}, m)\cap {\bf L}_{k_1-1}\}\cup \{{\bf W}_r({\bf k}, m)\cap {\bf L}_{1}\}|\leq \bar{\beta} k_1^{\delta/2}k_2\cdots k_{d-1}. \eqno{(8.3)}$$
If the cutset satisfies (8.3), we select ${\bf L}_1$ and ${\bf L}_{k_1-1}$. If it does not, we continue this process
until we find the first hyperplanes ${\bf L}_\tau$ and  ${\bf L}_{k_1-\tau}$ such that 
$$|\{{\bf W}_r({\bf k}, m)\cap {\bf L}_{k_1-\tau}\}\cup \{{\bf W}_r({\bf k}, m)\cap {\bf L}_{\tau}\}|\leq \bar{\beta }k_1^{\delta/2}k_2\cdots k_{d-1}. \eqno{(8.4)}$$
Note that the total number of vertices in  a regular cutset is less than $\bar{\beta} k_1k_2\cdots k_{d-1}$,
so we need to do this process at most $k_1^{1-\delta/2}$ times to find the hyperplanes. In other words,
$$\tau \leq k_1^{1-\delta/2}.\eqno{(8.5)}$$
 By (8.5), there exists $0< l< k^{1-\delta/2}$ such that
$${\bf P}[{\cal B}({\bf k}, m, l)]\geq {1\over 2k^{1-\delta/2}_1},\eqno{(8.6)}$$
where ${\cal B}({\bf k}, m, l)$ is the event  that $l$ is the first hyperplane
with the property (8.4).

For the fixed $l \leq  k^{1-\delta/2}$ defined in (8.6),
we collect all  cutsets $\{{\bf W}_r({\bf k}, m,l)\}$ in
$${\bf B}({\bf k}, m, l)=[l, k_1-l]\times [0, k_2]\times \cdots [0,k_{d-1}]\times [0, m]$$
cutting the bottom from the top of ${\bf B}({\bf k}, m, l)$ such that
$$ |{\bf W}_r({\bf k}, m, l)|\leq \bar{\beta} \|{\bf k}\|_v,
|\{{\bf W}_r({\bf k}, m, l)\cap {\bf L}_{l}\}\cup \{{\bf W}_r({\bf k}, m, l)\cap {\bf L}_{k_1-l}\}|\leq \bar{\beta} k_1^{\delta/2}k_2\cdots k_{d-1}.\eqno{(8.7)}$$
We select one from these cutsets, still denoted by ${\bf W}_r({\bf k}, m,l)$, with the minimum passage time:
$$\tau_{r} ({\bf k}, l,m)= \tau({\bf W}_r({\bf k}, l,m))$$
If ${\bf W}_r({\bf k}, m,l)$ is not unique, we select ${\bf W}_r({\bf k}, m,l)$
with the minimum number of edges in a unique method of selection. 
By our definition,
$$\mbox{ $\tau_r({\bf k}, m, l)$ only depends on the configurations of the edges in }{\bf B}({\bf k}, m, l).\eqno{(8.8)}$$

{\bf Lemma 17.} 
{\em On ${\cal B}({\bf k}, m, l)$,}
$$\tau_{r} ({\bf k}, m, l)\leq \tau_{r}({\bf k}, m).$$

{\bf Proof.} By Lemma 12 (a), 
$${\bf W}_r({\bf k}, m)\cap {\bf B}({\bf k}, m, l)\eqno{(8.9)}$$
is a cutset that cuts the bottom from the top of ${\bf B}({\bf k}, m, l)$.
On the other hand, on ${\cal B}({\bf k}, m, l)$, 
the cutset in (8.9) satisfies (8.7). Therefore, Lemma 17 follows. $\Box$\\

We use  $\{{\bf T}_{(l,j)}\},$ $\{{\bf T}_{(k_1-l,j)}\}$, $\{{\bf S}_{(l,i)}\}$, and $\{{\bf S}_{(k_1-l,i)}\}$ 
to denote all the exits of the upper and the lower tunnels on the hyperplanes
${\bf L}_l$ and ${\bf L}_{k_1-l}$ 
for the cutset ${\bf W}_r({\bf k}, m,l)$, respectively.
For given positive integers $t_1, t_2, s_1, s_2$, we now define the events
\begin{eqnarray*}
\{{\cal I}_{t_1,t_2}\}=&&\{ \exists \mbox{ $t_1$ exits of the upper tunnels ${\bf T}_{(l,1)},\cdots, {\bf T}_{(l,t_1)}$ on ${\bf L}_l$}\\
&&\mbox{ and $\exists$ $t_2$ exits of the upper tunnels ${\bf T}_{(k_1-l,1)},\cdots, {\bf T}_{(k_1-l,t_2)}$ on ${\bf L}_{k_1-l}$
}  \}\\
\{{\cal J}_{s_1,s_2}\}=&&\{ \exists \mbox{ $s_1$ exits of  the lower tunnels ${\bf S}_{(l,1)},\cdots, {\bf S}_{(l,s_1)}$ on ${\bf L}_l$}\\
&&\mbox{ and $\exists$ $s_2$ exits of the lower tunnels ${\bf S}_{(k_1-l,1)},\cdots, {\bf S}_{(k_1-l, s_2)}$ on ${\bf L}_{k_1-l}$
}  \}.
\end{eqnarray*}
On $\{{\cal I}_{t_1, t_2}\}\cap \{{\cal J}_{s_1, s_2}\}$,   note that ${\bf W}_r({\bf k}, m,l)$ is uniquely selected, so 
the exits of the lower and upper tunnels are also uniquely determined. Thus, 
we decompose  the exits of the tunnels to fixed sets:
\begin{eqnarray*}
1&=&{\bf P}[\exists \mbox{ cutset }{\bf W}_r({\bf k}, m,l)]\\
&=& \sum_{t_1, t_2}\sum_{s_1, s_2} \sum_{\scriptstyle {\Gamma_{(l,1)},\cdots, \Gamma_{(l,t_1)}} 
\atop{\beta_{(l,1)},\cdots, \beta_{(l,s_1)}}}
\sum_{\scriptstyle {\Gamma_{(k_1-l,1)},\cdots, \Gamma_{(k_1-l, t_2)}} \atop{\beta_{(k_1-l, 1)},\cdots, \beta_{(k_1-l, s_2)}}}
{\bf P}[\,\, \exists \,\,{\bf W}_r({\bf k}, m, l), {\cal I}_{t_1,t_2},{\cal J}_{s_1, s_2} ,\bigcap_{j=1}^{t_1}\{{\bf T}_{(l,j)}=\Gamma_{(l,j)}\}, 
 \\
&&\hskip 3cm \bigcap_{i=1}^{s_1} \{{\bf S}_{(l,i)}=\beta_{l,i}\},
\bigcap_{j=1}^{t_2} \{{\bf T}_{(k_1-l,j)}=\Gamma_{(k_1-l,j)}\},\bigcap_{i=1}^{s_2}\{{\bf S}_{(k_1-l,i)}=\beta_{(k_1-l,i)}\}],
\,\,\,\, (8.10)
\end{eqnarray*}
where the first two sums above take over all possible $t_{1}$, $t_2$ $s_1$, and $s_2$, and the last two sums take all possible groups of 
fixed clusters  such that each group of clusters
$$\Gamma_{(l,1)}, \Gamma_{(l,2)}, \cdots, \Gamma_{(l,t_1)}, \beta_{(l,1)},\beta_{(l,2)},\cdots, \beta_{(l,s_1)}\subset {\bf L}_{l}\cap {\bf B}({\bf k}, m, l)\eqno{}$$
and 
$$\Gamma_{(k_1-l,1)}, \Gamma_{(k_1-l, 2)}, \cdots, \Gamma_{(k_1-l,t_2)}, \beta_{(k_1-l,1)},\beta_{k_1-l, 2},\cdots, \beta_{(k_1-l, s_2)}\subset {\bf L}_{k_1-l}\cap {\bf B}({\bf k}, m, l).\eqno{(8.11)}$$

For simplicity, we denote each group of clusters by, 
\begin{eqnarray*}
&&{\Gamma}_{(l,I)} =\{{\Gamma}_{l,1,}, \cdots, {\Gamma}_{l,{t}_1}\},\,\,\,\,\, {\Gamma}_{(l,II)} =\{{\Gamma}_{k_1-l,1,}, \cdots, {\Gamma}_{l,{t}_2}\}, \\
&&{\beta}_{(l, I)}=\{{\beta}_{l,1,}, \cdots, {\beta}_{l,{s}_1}\},\,\,\,\,\,{\beta}_{(l, II)}=\{{\beta}_{k_1-l,1,}, \cdots, {\beta}_{k_1-l,{s}_2}\}.
\end{eqnarray*}
We also denote the event in the probability of the right side of (8.10) 
for the group of clusters in  (8.11) by
$${\cal D}_1(t_1, {t}_2, {s}_1,{s}_2, {\Gamma}_{(l, I)}, {\Gamma}_{(k_1-l, II)},
{\beta}_{(l, I)}, {\beta}_{(k_1-l, II)}).$$
Note that for some groups of clusters, we have
$${\bf P}\left[ {\cal D}_1({t}_1, {t}_2, {s}_1,{s}_2, {\Gamma}_{(l, I)}, {\Gamma}_{(k_1-l, II)},
{\beta}_{(l, I)}, {\beta}_{(k_1-l, II)})\right]=0.$$
In these cases, the groups of clusters are trivial  and we will not account for these terms in the four sums in (8.10).
Note also that 
the four sums only take finitely many terms, so there is a term with the largest probability among the others. 
We denote this largest term with the indexes $\bar{t}_1, \bar{t}_2, \bar{s}_1, \bar{s}_2$ and denote the group of
clusters by
\begin{eqnarray*}
&&\bar{\Gamma}_{(l,I)} =\{\bar{\Gamma}_{l,1,}, \cdots, \bar{\Gamma}_{l,\bar{t}_1}\},\,\,\,\,\, \bar{\Gamma}_{(l,II)} =\{\bar{\Gamma}_{k_1-l,1,}, \cdots, \bar{\Gamma}_{l,\bar{t}_2}\}, \\
&&\bar{\beta}_{(l, I)}=\{\bar{\beta}_{l,1,}, \cdots, \bar{\beta}_{l,\bar{s}_1}\},\,\,\,\,\,\bar{\beta}_{(l, II)}=\{\bar{\beta}_{k_1-l,1,}, \cdots, \bar{\beta}_{k_1-l,\bar{s}_2}\}.
\end{eqnarray*}
We also define
\begin{eqnarray*}
&&\max_{\scriptstyle{t_1, t_2, s_1, s_2, \Gamma_{l,1},\cdots \Gamma_{l,t_1},\beta_{l,1}\cdots \beta_{l,s_1},}\atop{ \Gamma_{k_1-l, 1},\cdots  \Gamma_{k_1-l, t_2},
\beta_{k_1-l,1},\cdots \beta_{k_1-l, s_2} }}
{\bf P}[ {\cal I}_{t_1,t_2}, {\cal J}_{s_1, s_2} ,\bigcap_{j=1}^{t_1}\{T_{l,j}=\Gamma_{l,j}\}, 
 \bigcap_{i=1}^{s_1} \{S_{l,i}=\beta_{l,i}\}\\
&& \hskip 6cm \bigcap_{j=1}^{t_2} \{T_{k_1-l,j}=\Gamma_{k_1-l,j}\},\bigcap_{i=1}^{s_2}\{S_{k_1-l,i}=\beta_{k_1-l,1})\}] \\
&&:={\bf P}\left[{\cal D}_1(\bar{t}_1, \bar{t}_2, \bar{s}_1,\bar{s}_2, \bar{\Gamma}_{(l, I)}, \bar{\Gamma}_{(k_1-l, II)},
\bar{\beta}_{(l, I)}, \bar{\beta}_{(k_1-l, II)})\right]. \hskip 2in 
\end{eqnarray*}

It is possible that there is another group of clusters with the same largest probability. If this occurs, we
select one group in a unique method.
We will account for the number of non-trivial groups of the clusters in the four sums in (8.10). 
In other words, we need to account 
for all possible groups of 
clusters on ${\bf L}_{k_1-l}\cup {\bf B}({\bf k}, m, l)$ and  ${\bf L}_{k_1-l}\cup {\bf B}({\bf k}, m, l)$ such that
they are the exits of upper or lower tunnels for ${\bf W}_r({\bf k}, m,l)$. Let $N_r({\bf k}, m, l)$ be the number of
 all the possible non-trivial groups of clusters above. We will then give an upper bound estimate.
For fixed $l\leq k^{1-\delta/2}$, let
$${\bf I}({\bf k}, m,l)= \{{\bf W}_r({\bf k},m,l)\cap {\bf L}_{l}\}\cup \{{\bf W}_r({\bf k},m,l)\cap {\bf L}_{k_1-l}\} .$$
By the definition,
$$|{\bf I}({\bf k}, m,l)|\leq \bar{\beta} k_1^{\delta/2}k_2\cdots k_{d-1}. \eqno{}$$
If we use ${\bf I}_e({\bf k}, m,l)$ to denote  the edge set in ${\bf I}({\bf k}, m, l)$, then
$$\mbox{ the number of edges in }{\bf I}_e({\bf k}, m, l)\leq 2d\bar{\beta} k_1^{\delta/2}k_2\cdots k_{d-1}. \eqno{(8.12)}$$
Note that the total number of vertices of ${\bf B}({\bf k}, m)$
on the two faces
$$|{\bf B}({\bf k}, m)\cap {\bf L}_{l}\cup {\bf B}({\bf k}, m)\cap {\bf L}_{k_1-l}|\leq  2k_1\cdots k_{d-1}m.\eqno{(8.13)}$$
Note also that for a cluster 
on ${\bf L}_l$ and on ${\bf L}_{k_1-l}$, 
if its boundary edges are fixed, then
 the precise location of the cluster  is uniquely fixed.
By Lemma 15 (a), the  boundary edges of the exits of the upper and the lower tunnels belong to 
${\bf I}_e({\bf k}, m,l)$.
Note that if we remove ${\bf I}_e({\bf k}, m,l)$ from both ${\bf L}_l$ and ${\bf L}_{k_1-l}$, we can view the remaining
edges as many clusters. These clusters are the exits of the 
upper and the lower tunnels. With these clusters, we need to identify the upper or the lower exists from them.
Given a fixed ${\bf I}_e({\bf k}, m,l)$, suppose that there are $q$ clusters, as
the exits of the upper and the lower tunnels on 
both ${\bf L}_{l}$ and  ${\bf L}_{k_1-l}$,
after removing the edges of ${\bf I}_e({\bf k}, m,l)$ from ${\bf L}_l$ 
and ${\bf L}_{k_1-l}$.
Note that if we remove one edge, it can separate  one cluster into at most two clusters.
Therefore,  by (8.12) after removing ${\bf I}_e({\bf k}, m,l)$, the total number of the clusters
of these exits of the upper and the lower tunnels is 
$$q\leq 2|{\bf I}_e({\bf k}, m,l)|\leq 4d \bar{\beta} k_1^{\delta/2}k_2\cdots k_{d-1}. \eqno{(8.14)}$$

Among these $q$ fixed clusters, we select some of them as the exits of the  upper  and lower tunnels
on ${\bf L}_l$ and on ${\bf L}_{k_1-l}$.
By (8.12) and (8.14), the number of selections is at most
$$\sum _{t_1=1}^q \sum_{t_2=1}^q {q\choose t_1 }{q\choose t_2 }  \leq 2^{2q}\leq 2^{8d \bar{\beta} k_1^{\delta/2}k_2\cdots k_{d-1}}.\eqno{(8.15)}$$

With these observations, 
 we first select  ${\bf I}({\bf k}, m,l)$, defined above,
on ${\bf L}_l$ and ${\bf L}_{k_1-l}$. With the first selection, ${\bf I}_e({\bf k},  m,l)$ is determined. After removing 
${\bf I}_e({\bf k},  m,l)$, the remaining clusters, the exits of the upper and the lower tunnels,  are determined. 
We then select  the exits for  the upper  and for the lower tunnels from
these clusters. With these selections and (8.12)--(8.15), 
 the total number, $N_r({\bf k}, m, l)$,  of all the possible exits of the upper and the 
lower tunnels 
is at most
$$N_r({\bf k}, m, l)\leq  \bar{\beta} k_1^{\delta/2} k_2\cdots k_{d-1} {2k_1 k_2\cdots k_{d-1}m  \choose \bar{\beta } k_1^{\delta/2} k_2\cdots k_{d-1}} 2^{8d \bar{\beta} k_1^{\delta/2}k_2\cdots k_{d-1}}.\eqno{(8.16)}
$$
By using Corollary 2.6.2 in Engle (1997), 
$$N_r({\bf k}, m, l)\leq C_1 \exp\left[C_2 k_1^{\delta/2}\cdots k_{d-1} \log (k_1\cdots k_{d-1} m)\right]\eqno{(8.17)}$$
for $C_i=C_i(F, d, \beta, \delta)$, $i=1,2$.
Furthermore, if we assume that $k_1\leq k_2\leq \cdots k_{d-1}$ with $ k_{d-1}\leq  2\exp(10k_1^{1-5\delta/6})$ and 
$\log m \leq k_{d-1}^{1-\delta}$ (the assumptions for $m$ in (1.17)), there exist $C_i=C_i(F, d, \delta)$ such that 
$$ N_r({\bf k}, m,l)\leq C_1 \exp\left[C_2 k_1^{\delta/2} k_2\cdots k_{d-1} \log (k_{d-1}m)\right]\leq
C_3 \exp\left[C_4 k_1^{1-\delta/3} k_2\cdots k_{d-1}\right].\eqno{(8.18)}$$
Therefore,  the number of all terms in these 
four sums in (8.10) is at most $N_r({\bf k}, m,l)$. 
With these observations, by (8.18),
\begin{eqnarray*}
1&=&{\bf P}[\exists \mbox{ the cutset }{\bf W}_r({\bf k}, m,l)]\\
&=& \sum_{t_1, t_2}\sum_{s_1, s_2} \sum_{\scriptstyle {\Gamma_{(l,1)},\cdots \Gamma_{(l,t_1)}} 
\atop{\beta_{(l,1)}\cdots \beta_{(l,s_1)}}}
\sum_{\scriptstyle {\Gamma_{(k_1-l,1)},\cdots \Gamma_{(k_1-l, t_2)}} \atop{\beta_{(k_1-l, 1)}\cdots \beta_{(k_1-l, s_2)}}}
{\bf P}\left[{\cal D}_1({t}_1, {t}_2, {s}_1,{s}_2, {\Gamma}_{(l, I)}, {\Gamma}_{(k_1-l, II)},
{\beta}_{(l, I)}, {\beta}_{(k_1-l, II)})\right]\\
&\leq &
{\bf P}\left[{\cal D}_1(\bar{t}_1, \bar{t}_2, \bar{s}_1,\bar{s}_2, \bar{\Gamma}_{(l, I)}, \bar{\Gamma}_{(k_1-l, II)},
\bar{\beta}_{(l, I)}, \bar{\beta}_{(k_1-l, II)})\right]N_r({\bf k}, m,l)\\
&\leq & {\bf P}\left[{\cal D}_1(\bar{t}_1, \bar{t}_2, \bar{s}_1,\bar{s}_2, \bar{\Gamma}_{(l, I)}, \bar{\Gamma}_{(k_1-l, II)},
\bar{\beta}_{(l, I)}, \bar{\beta}_{(k_1-l, II)})\right]
C_1 \exp\left(C_2 k_1^{1-\delta/3} k_2\cdots k_{d-1} \right).
\end{eqnarray*}
If we simply denote by 
$${\cal D}_1(\bar{t}_1, \bar{t}_2, \bar{s}_1,\bar{s}_2, \bar{\Gamma}_{(l, I)}, \bar{\Gamma}_{(k_1-l, II)},
\bar{\beta}_{(l, I)}, \bar{\beta}_{(k_1-l, II)})={\cal D}_1,$$
we summarize the above result as the following lemma.\\

{\bf Lemma 18.} {\em If $k_1\leq k_2\leq \cdots \leq k_{d-1}$ with $ k_{d-1}\leq  2\exp(10k_1^{1-5\delta/6})$, 
and $\log m\leq k_{d-1}^{1-\delta}$, then there are constants $C_i=C_i(F, d, \beta,\delta)$ for $i=1,2$ such that}
$$C_1 \exp\left(-C_2 k_1^{1-\delta/3} k_2\cdots k_{d-1} \right)\leq {\bf P}[{\cal D}_1].$$

If we work on $k_j$'s direction rather than $k_1$'s, similar to ${\cal D}_1$, let ${\cal D}_j$ be the event 
corresponding to the $j$-th coordinate. By the same estimate, we have the following
the  result in  Lemma 18 for ${\cal D}_j$ holds.\\

{\bf Lemma 19.} {\em If $k_1\leq k_2\leq \cdots \leq  k_{d-1}$ with $ k_{d-1}\leq  2\exp(10k_1^{1-5\delta/6})$, and $\log m\leq k_{d-1}^{1-\delta}$, then there are constants $C_i=C_i(F, d, \beta,\delta)$ for $i=1,2$ such that}
$$C_1 \exp\left(-C_2 k_1\cdots  k_j^{1-\delta/3} \cdots k_{d-1} \right)\leq {\bf P}[{\cal D}_j].$$
{\em In particular, if $k_1\leq k_2\leq \cdots \leq  k_{d-1}$ and $\log m\leq k_{d-1}^{1-\delta}$, then there are constants $C_i=C_i(F, d, \beta,\delta)$ for $i=1,2$ such that}
 $$C_1 \exp\left(-C_2 k_1\cdots  k_{d-2} \cdots k_{d-1}^{1-\delta/3} \right)\leq {\bf P}[{\cal D}_{d-1}].$$

\section{ Concentration of $\tau_r({\bf k}, m)$ from its mean.}
In general, there are two major methods to estimate the concentration inequalities.
Kesten (1993) has investigated the concentration for the first passage percolation by using a martingale argument. 
Later, Talagrand (1995)  obtained a better  result by using the {\em isoperimetric} inequality.
Both ways can be carried out to investigate the concentration for the passage time of a minimal cutset from its mean. 
We use the Talagrand method in this paper. 
Denote by ${\cal S}$ the sets of all regular cutsets $\{{{\bf Z}_r({\bf k}, m)}\}$,
defined in section 8, with the minimum passage time. Let 
$$\alpha =\sup_{{\bf Z}_r({\bf k}, m)\in {\cal S}} |{\bf Z}_r({\bf k}, m)|.$$
It follows from this definition
$$\alpha\leq \bar{\beta} \|{\bf k}\|_v.\eqno{(9.0)}$$
Denote by $M$ a median of $\tau_{r}({\bf k}, m)$. By Theorem (8.3.1) (see Talagrand (1995)) there exist constants
$C$ and $C_1$ 
such that
$${\bf P}\left[|\tau_{r}({\bf k}, m)-M|\geq u\right]\leq C\exp\left (-C_1\min\left \{{u^2\over \alpha}, u\right\}\right).\eqno{(9.1)}$$

By (9.0) and (9.1), for all $u >0 $, 
$${\bf P}[|\tau_{r}({\bf k}, m)-M|\geq u]\leq C\exp\left(-C_1\min\left\{{u^2\over \bar{\beta }\|{\bf k}\|_v}, u\right \}\right).\eqno{(9.2)}$$
If we  select  $u$  satisfying 
$$(\|{\bf k}\|_v)^{2/3}\leq u, \eqno{(9.3)}$$
then
$${\bf P}[ |\tau_{r}({\bf k}, m)-M|\geq u]\leq C_1\exp\left(-C_2 \|{\bf k}\|_v^{1/3}\right).\eqno{(9.4)}$$
By (9.4),
$$|{\bf E}\tau_{r}({\bf k}, m)-M|\leq {\bf E}|\tau_{r}({\bf k}, m)-M|
= \sum_{i=1} P(|\tau_{r}({\bf k}, m)-M|\geq i)\leq C (\|{\bf k}\|_v)^{2/3}.\eqno{(9.5)}$$
Therefore, for all large $k_1,\cdots, k_{d-1}$, and for $u$ with
$$ \max\{2C(\|{\bf k}\|_v)^{2/3},\bar{\beta}(\|{\bf k}\|_v)^{2/3}\}\leq u \mbox{ for the $C$ in (9.5)},\eqno{(9.6)}$$
then by (9.2) and (9.5),
\begin{eqnarray*}
&&{\bf P}[|\tau_{r}({\bf k}, m)-{\bf E}\tau_{r}({\bf k}, m)|\geq u]\\
&\leq &
{\bf P}[|\tau_{r}({\bf k}, m)-M| +|M-{\bf E}\tau_{r}({\bf k}, m)|\geq u]\leq
C\exp\left(-C_1{u^2\over \|{\bf k}\|_v}\right).\hskip 1.3in {(9.7)}
\end{eqnarray*}

If we focus on $\tau_r({\bf k}, m, l)$, the passage time of cutsets ${\bf W}_r({\bf k}, m, l)$,
then by the same estimates in (9.1)--(9.7), for all large $k_1,\cdots, k_{d-1}$, and for the $u$ satisfying (9.6), we have
$${\bf P}\left[|\tau_{r}({\bf k}, m,l)-E\tau_{r}({\bf k}, m,l)|\geq u\right]\leq C\exp\left(-C_1{u^2\over \|{\bf k}\|_v}\right).\eqno{(9.8)}$$

Now we will try to use the concentration property to estimate the means of $\tau_{r}({\bf k}, m)$ and
$\tau_{r}({\bf k}, m,l)$ on some event ${\cal E}$ that may depend on ${\bf k}$ and $m$.\\

{\bf Lemma 20.} {\em Under (1.1), there exist $C_i=C_i(F, d,  \beta, \delta)$ for $i=1,2$ such that
for each $k_j$, $j=1,2,\cdots, d-1$ and $0< \delta\leq 1$, 
\begin{eqnarray*}
&&|{\bf E} [\tau_{r}({\bf k}, m)]-{\bf E} [ \tau_{r}({\bf k}, m)\,\, \,|\,\,\, {\cal E}\,\,\,]|\\
&&\leq  
C_1 k_1  \cdots k_j^{(1-\delta/8)}\cdots k_{d-1}
+ C_1\left\{ {\bf P} ({\cal E})\right\}^{-1} \exp\left(-C_2 k_1 \cdots k_j^{(1-\delta/4)}\cdots k_{d-1}\right)\hskip 3cm {}
\end{eqnarray*}
and 
\begin{eqnarray*}
&&|{\bf E}[\tau_{r}({\bf k}, m,l)]-{\bf E} [ \tau_{r}({\bf k}, m,l)\,\, \,|\,\,\, {\cal E}\,\,\,]|\\
&&\leq  
C_1 k_1  \cdots k_j^{(1-\delta/8)}\cdots k_{d-1}
+ C_1\left\{ {\bf P} ({\cal E})\right\}^{-1} \exp\left(-C_2 k_1 \cdots k_j^{(1-\delta/4)}\cdots k_{d-1}\right).
\end{eqnarray*}}

{\bf Proof.} Without loss of generality, we show Lemma 20 for $j=1$.
We begin with an estimate of 
$${\bf E} [|{\bf E}\tau_{r}({\bf k}, m)-  \tau_{r}({\bf k}, m)|\,\, \,|\,\,\, {\cal E}\,\,\,]. $$
Denote the event ${\cal L}({\bf k}, m)$ by
$${\cal L}({\bf k}, m)=\left\{|{\bf E}\tau_{r}({\bf k}, m)-\tau_{r}({\bf k}, m)|>\bar{\beta} k_1^{(1-\delta/8)}k_2\cdots k_{d-1}\right\}.$$
We divide
\begin{eqnarray*}
&&{\bf E} \left[|{\bf E}\tau_{r}({\bf k}, m)-\tau_{r}({\bf k}, m)|
\,\, \,|\,\,\, {\cal E}\,\,\,\right]\\
&&={\bf E}\left[|{\bf E}\tau_{r}({\bf k}, m)-\tau_{r}({\bf k}, m)|I( {\cal L}({\bf k}, m))\,\, \,|\,\,\, {\cal E}\,\,\right]\\
&&+ {\bf E} \left[|{\bf E}\tau_{r}({\bf k}, m)-\tau_{r}({\bf k}, m)| (1-I({\cal L}({\bf k}, m))\,\, \,|\,\,\, 
{\cal E}\,\,\,\right]\\
&&=I+II,
\end{eqnarray*}
where $I({\cal A})$ is the indicator for the event ${\cal A}$.
By the definition of ${\cal L}({\bf k}, m)$,
$$II\leq \bar{\beta }k_1^{(1-\delta/8)}k_2\cdots k_{d-1}.\eqno{(9.9)}$$
We  estimate $I$.  By (9.7), there exist $C$ and $C_1$ such that
\begin{eqnarray*}
I &&\leq \sum_{i\geq \bar{\beta} k_1^{(1-\delta/8)}k_2\cdots k_{d-1}} {\bf P}\left[ |\tau_r({\bf k}, m)-{\bf E}(\tau_r({\bf k}, m)|\geq i\,\,\, |\,\,\, {\cal E}\,\,\,\right]\\
&& =\sum_{i\geq \bar{\beta }k_1^{(1-\delta/8)}k_2\cdots k_{d-1}} {{\bf P}\left[ |\tau_r({\bf k}, m)-{\bf E}(\tau_r({\bf k}, m)|\geq i\right]\over {\bf P}[ {\cal E}]}\\
&&\leq  C\{ {\bf P} [{\cal E}]\}^{-1} \exp\left(-C_1 k_1^{(1-\delta/4)}k_2\cdots k_{d-1}\right).\hskip 6cm (9.10)
\end{eqnarray*}
Combining (9.9) and (9.10), we have
\begin{eqnarray*}
&&{\bf E} \left[|{\bf E}\tau_{r}({\bf k}, m)-  \tau_{r}({\bf k}, m)|\,\, \,|\,\,\, {\cal E}\,\,\,\right]\\
&&\leq C k_1^{(1-\delta/8)} k_2 \cdots k_{d-1}
+ C\{ {\bf P} [{\cal E}]\} ^{-1}\exp\left(-C_1 k_1^{(1-\delta/4)} k_2\cdots k_{d-1}\right).\hskip 2.5cm {(9.11)}
\end{eqnarray*}
Using the same estimate of (9.11) together with (9.8),
\begin{eqnarray*}
&&{\bf E} \left[|{\bf E}\tau_{r}({\bf k}, m,l)-  \tau_{r}({\bf k}, m,l)|\,\, \,|\,\,\, {\cal E}\,\,\,\right]\\
&&\leq Ck_1^{(1-\delta/8)} k_2 \cdots k_{d-1}
+ C\{ {\bf P} [{\cal E}]\}^{-1} \exp\left(-C k_1^{(1-\delta/4)} k_2\cdots k_{d-1}\right).\hskip 3cm {(9.12)}
\end{eqnarray*}

With (9.11) and (9.12), let us  show Lemma 20.
We then have 
$${\bf E}\tau_{r}({\bf k}, m)={\bf E} \left[ \tau_{r}({\bf k}, m)\,\,\,|\,\,\, {\cal E}\,\,\right]
+{\bf E} \left[{\bf E}\tau_{r}({\bf k}, m)-\tau_{r}({\bf k}, m)]\,\,\,|\,\,\, {\cal E}\,\,\,\right].\eqno{(9.13)}$$
By (9.13), we have
$$|{\bf E} \left[\tau_{r}({\bf k}, m)\right]-{\bf E} \left[ \tau_{r}({\bf k}, m)\,\, \,|\,\,\, {\cal E}\,\,\right]|
\leq {\bf E}\left[|{\bf E}\tau_{r}({\bf k}, m)-  \tau_{r}({\bf k}, m)|\,\, \,|\,\,\, {\cal E}\,\,\right].\eqno{(9.14)}$$
Therefore, by (9.11) and (9.14), there exists $C_i=C_i(F, d,\beta,\delta)$ for $i=1,2$ such that
\begin{eqnarray*}
&&|{\bf E} \left[\tau_{r}({\bf k}, m)]-{\bf E} [ \tau_{r}({\bf k}, m)\,\, \,|\,\,\, {\cal E}\,\,\right]|\\
&&\leq 
C_1 k_1^{(1-\delta/8)} k_2 \cdots k_{d-1}
+ C_1\{ {\bf P}[{\cal E}]\}^{-1} \exp\left(-C_2 k_1^{1-\delta/4} k_2\cdots k_{d-1}\right).\hskip 3cm {(9.15)}
\end{eqnarray*}
The same estimate in (9.15) also shows that
\begin{eqnarray*}
&&|{\bf E} \left[\tau_{r}({\bf k}, m,l)]-{\bf E} [ \tau_{r}({\bf k}, m,l)\,\, \,|\,\,\, {\cal E}\right]|\\
&&\leq  
C_1k_1^{(1-\delta/8)} k_2 \cdots k_{d-1}
+ C_1\{ {\bf P}[{\cal E}]\}^{-1} \exp\left(-C_2 k_1^{(1-\delta/4)} k_2\cdots k_{d-1}\right).\hskip 3cm {(9.16)}
\end{eqnarray*}
Lemma 20, for $j=1$, follows from (9.15) and (9.16). $\Box$\\

\section{Proof of Theorem 3.}
As we pointed out in section 1, we only need to show Theorem 3 when $F(0) < 1-p_c$. Thus, we assume that 
$F(0) < 1-p_c$ in this section. Note that 
$$0\leq {\bf E}\tau_{\min}({\bf k}, m)\leq {\bf E}\bar{\alpha}({\bf k}, m)\leq C\|{\bf k}\|_v {\bf E}(\tau(e)),\eqno{(10.0)}$$
so we assume that  there exist $0\leq \nu_1\leq \nu_2< \infty $
such that
$$\nu_1=\liminf_{k_1,\cdots, k_{d-1},m\rightarrow \infty} \left( {{\bf E} \tau_{\min}({\bf k}, m)\over \|{\bf k}\|_v}\right)
\leq \nu_2=\limsup_{k_1,\cdots, k_{d-1},m\rightarrow \infty} \left( {{\bf E} \tau_{\min}({\bf k}, m)\over \|{\bf k}\|_v}\right).
\eqno{(10.1)}$$
We first show that $\nu_1=\nu_2$. 
The key proof of this argument is to 
show a multiple subadditive property  for $E\tau_{\min} ({\bf k}, m)$.

Now we assume that
$$k_1\leq k_2\leq \cdots\leq k_{d-1} \mbox{ with } k_{d-1}\leq 4k_1. \eqno{(10.2)}$$
Besides ${\bf B}({\bf k}, m,l)$ defined in section 8, we also denote by (see Fig. 5)
$${\bf B}'({\bf k}, m,l)=[k_1-l+1, 2k_1-3l+1]\times [0, k_2]\times \cdots \times [0,k_{d-1}]\times [0,m], $$
$${\bf B}''({\bf k}, m,l)=[2k_1-3l+2, 3k_1-5l+2]\times [0, k_2]\times \cdots \times[0,k_{d-1}]\times [0,m]. $$
We denote by $\omega({\bf B}({\bf k}, m,l))$, $\omega({\bf B}'({\bf k}, m,l))$, and $\omega({\bf B}''({\bf k}, m,l))$
the configurations on ${\bf B}({\bf k}, m,l)$, ${\bf B}'({\bf k}, m,l)$, and ${\bf B}''({\bf k}, m,l)$, respectively.
For each $\omega({\bf B}'({\bf k}, m,l))$ and $\omega({\bf B}''({\bf k}, m,l))$, we can select the unique cutsets
${\bf W}_r'({\bf k}, m,l)$ and ${\bf W}_r''({\bf k}, m,l)$ in ${\bf B}'({\bf k}, m,l)$ and ${\bf B}''({\bf k}, m,l)$,
respectively, using the same rule for  selecting  ${\bf W}_r({\bf k}, m,l)$ on ${\bf B}({\bf k}, m, l)$.

Recall that  ${\cal D}_1$, defined in section 8, is the event with the largest probability for the fixed  
 exits of the upper and the lower tunnels for ${\bf W}_r({\bf k}, m,l)$. 
Similarly, let ${\cal D}'_1$ and ${\cal D}''_1$ be the events 
with the largest probabilities for  the fixed   exits of the upper and the lower tunnels for 
${\bf W}_r'({\bf k}, m,l)$ and  ${\bf W}_r''({\bf k}, m,l)$, respectively,   the same as for
${\bf W}_r({\bf k}, m,l)$ in the sense of translation.

 For each set ${\bf A}\subset {\bf B}'({\bf k}, m,l)$, we define a {\em mirror reflection} 
about ${\bf L}_{k_1-l+0.5}$ as follows.
For ${\bf u}=(u_1, \cdots , u_d)\in {\bf A}$, let 
$$\sigma_1(u)=(2k_1-2l+1-u_1, u_2, \cdots, u_d).$$
After the mirror reflection, $\sigma_1({\bf A})\subset {\bf B}({\bf k}, m,l)$. We move $\sigma_1({\bf A})$ along 
the first coordinate $k_1-2l+1$ units back to ${\bf B}'({\bf k}, m,l)$.
 More precisely, for each $u\in \sigma_1({\bf A})$, let
$$\sigma_2({\bf u})=(u_1+k_1-2l+1, u_2,\cdots , u_d).$$
We denote by 
$$\pi({\bf A})=\sigma_2(\sigma_1({\bf A})).$$
 With these changes, we have another vertex
set, denoted by $\pi ({\bf A})$ in  ${\bf B}'({\bf k}, m,l)$.  For each edge $e=({\bf u}, {\bf v})$ with
a configuration $\omega(e)$, let
$\pi (\omega(e)) $ be the same value $\omega(e)$  on the edge $(\pi({\bf u}), \pi ({\bf v}))$.
Thus, for $\omega({\bf B}'({\bf k}, m,l))$, $\pi (\omega({\bf B}'({\bf k}, m,l)))$ will be the configuration by changing
 each configuration
$ \omega(e)$ at $e$, to $\pi(\omega(e))$.
 With configurations $\{\pi (\omega({\bf B}'({\bf k}, m,l)))\}$, we consider  $\pi({\bf W}_r'({\bf k}, m,l))$ (see Fig. 6). By our definition, $\pi({\bf W}_r'({\bf k}, m,l))$ is still a self-avoiding regular
cutset
that cuts the bottom face from the top face of ${\bf B}'({\bf k}, m, l)$. Also, it has the minimum passage time
among all the other regular  cutsets. Since the selection of ${\bf W}_r'({\bf k}, m,l)$ is unique, the selection of $\pi({\bf W}_r'({\bf k}, m))$ is also unique. 
In addition, let (see Fig. 6)
$$\pi ({\cal D}'_1)=\{\pi (\omega({\bf B}'({\bf k}, m,l))): \omega({\bf B}'({\bf k}, m,l))\in {\cal D}'_1\}.$$
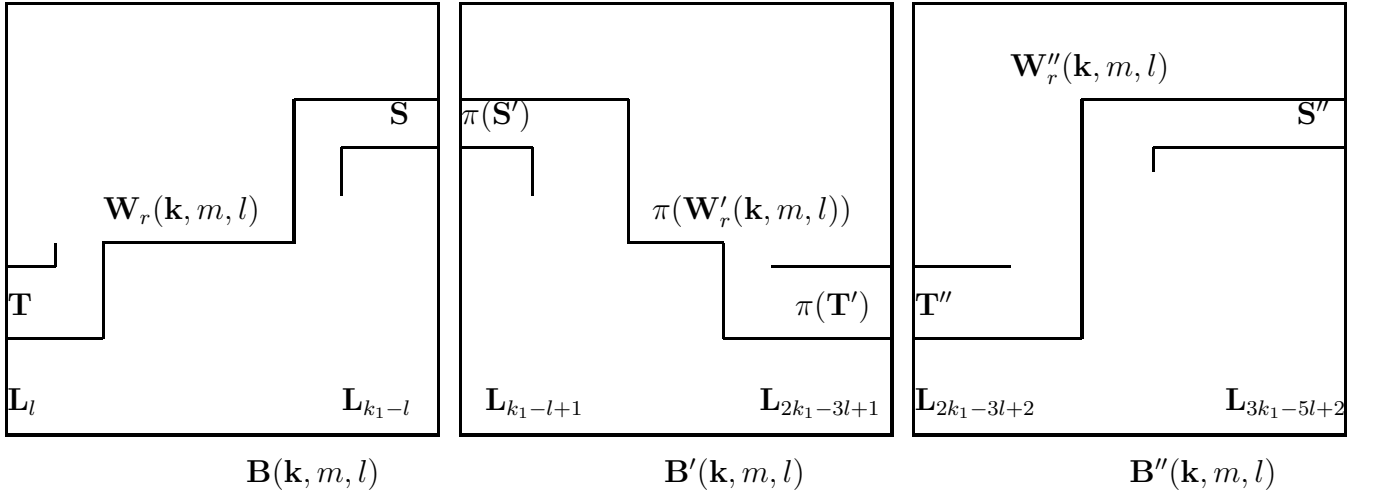
\begin{figure}\label{F:alphabeta}
\begin{center}
\setlength{\unitlength}{0.0125in}%
\begin{picture}(200,150)(67,800)
\thicklines
\put(-80, 840){\line (0,-1){10}}
\put(-80, 830){\line (-1,0){20}}
\put(-100,800){\line (1,0){40}}
\put(-60, 800){\line(0,1){40}}
\put(-60, 840){\line(1,0){80}}
\put(20, 840){\line(0,1){60}}
\put(20, 900){\line(1,0){60}}
\put(40, 860){\line(0,1){20}}
\put(40, 880){\line(1,0){40}}
\put(-100,770){\mbox{${\bf L}_l$}}
\put(-100,810){\mbox{${\bf T}$}}
\put(60,890){\mbox{${\bf S}$}}
\put(40,770){\mbox{${\bf L}_{k_1-l}$}}
\put(0,740){\mbox{${\bf B}({\bf k}, m, l)$}}

\put(100,770){\mbox{${\bf L}_{k_1-l+1}$}}
\put(90, 880){\line(1,0){30}}
\put(120, 880){\line(0,-1){20}}
\put(90, 900){\line(1,0){70}}
\put(160, 900){\line(0,-1){50}}
\put(90,890){\mbox{$\pi({\bf S}')$}}
\put(160, 850){\line(0,-1){10}}
\put(160, 840){\line(1,0){40}}
\put(200, 840){\line(0,-1){40}}
\put(200,800){\line(1,0){70}}
\put(270,830){\line(-1,0){50}}
\put(230,810){\mbox{$\pi({\bf T}')$}}
\put(215,770){\mbox{${\bf L}_{2k_1-3l+1}$}}
\put(175,740){\mbox{${\bf B}'({\bf k}, m,l)$}}

\put(280, 830){\line(1,0){40}}
\put(280, 800){\line(1,0){70}}
\put(350, 800){\line(0,1){100}}
\put(350, 900){\line(1,0){110}}
\put(460, 880){\line(-1,0){80}}
\put(380,880){\line(0,-1){10}}
\put(90,760){\framebox(180,180)[br]{\mbox{$$}}}
\put(-100,760){\framebox(180,180)[br]{\mbox{$$}}}
\put(280,760){\framebox(180,180)[br]{\mbox{$$}}}
\put(440,890){\mbox{${\bf S}''$}}
\put(280,810){\mbox{${\bf T}''$}}
\put(280,770){\mbox{${\bf L}_{2k_1-3l+2}$}}
\put(410,770){\mbox{${\bf L}_{3k_1-5l+2}$}}
\put(370,740){\mbox{${\bf B}''({\bf k}, m,l)$}}

\put(-60, 850){\mbox{${\bf W}_r({\bf k}, m,l)$}}
\put(170, 850){\mbox{$\pi({\bf W}_r'({\bf k}, m,l))$}}
\put(320, 910){\mbox{${\bf W}_r''({\bf k}, m,l)$}}

\end{picture}
\end{center}
\vskip 1in
\caption{ \em This cross-section graph shows that the three cutsets ${\bf W}_r({\bf k}, m,l)$, $\pi({\bf W}_r'({\bf k}, m, l))$, and
${\bf W}_r''({\bf k}, m, l)$ consist of  a three-times cutset when their corresponding exits are matched.}
\end{figure}
By  symmetry, we have 
\begin{eqnarray*}
&&{\bf E}[\tau(\pi({\bf W}_r'({\bf k}, m,l)))]={\bf E}[\tau_r({\bf k}, m,l)]={\bf E}[\tau({\bf W}_r''({\bf k}, m,l))]\\ 
&&{\bf P}[\pi({\cal D}'_1)]={\bf P}[{\cal D}_1]={\bf P}[{\cal D}''_1].\hskip 3.8in {(10.3)}
\end{eqnarray*}
By the definition of the mirror reflection $\sigma_1$ and the horizontal move $\sigma_2$, 
on ${\cal D}_1\cap \pi({\cal D}'_1)$,
the upper and the lower tunnels for ${\bf W}_r({\bf k}, m,l)$ and  $\pi({\bf W}_r({\bf k}, m, l))$
on ${\bf L}_{k_1-l}$ and  on ${\bf L}_{k_1-l+1}$ are matched, so by Lemma 16,
these two cutsets consist of a larger cutset (see Fig. 6):
\begin{eqnarray*}
&&{\bf W}_r({\bf k}, m, l)\cup \pi({\bf W}_r '({\bf k}, m, l))\mbox{ is a cutset that cuts the bottom from the top of }\\
&&[l,2 k_1-3l+1]\times [0,k_2]\times \cdots \times [0,k_{d-1}]\times[0,  m].\hskip 2.3in (10.4)
\end{eqnarray*}
Note also that the new cutset consists of ${\bf W}_r({\bf k}, m, l)$ and $\pi({\bf W}_r '({\bf k}, m, l))$.
Therefore, it is still regular.
With this observation, 
\begin{eqnarray*}
&&{\bf E}\left[\tau_{r} \left((2k_1-2l+1, k_2,\cdots, k_{d-1}), m,l\right) \,\,\, |\,\,\, {\cal D}_1\cap \pi({\cal D}')_1\right]\\
&&\leq {\bf E}[\tau({\bf W}_r({\bf k}, m, l))+
\tau(\pi({\bf W}_r'({\bf k}, m, l)))
\,\,\, |\,\,\, {\cal D}_1\cap \pi({\cal D}'_1)].\hskip 4.5cm (10.5)
\end{eqnarray*}

Note that ${\bf W}_r$ and ${\cal D}_1$, and $\pi({\bf W}_r')$ and $\pi({\cal D}'_1)$, 
 only depend on the configurations
of edges in different boxes, so
$$({\bf W}_r, {\cal D}_1), (\pi({\bf W}_r'), \pi({\cal D}'_1))\mbox{ are independent}.\eqno{(10.6)}$$
By  (10.5)--(10,6) and  symmetry,
$${\bf E}\left[\tau_{r} \left((2k_1-2l+1, k_2,\cdots, k_{d-1}), m,l\right)
 \,\,\, |\,\,\, {\cal D}_1\cap \pi({\cal D}'_1)\right]
\leq 2{\bf E}[\tau_r({\bf k}, m,l)\,\,\, | \,\,\, {\cal D}_1].\eqno{(10.7)}$$
We  need to use Lemma 20 to change conditional expectations in (10.7)
to unconditional  expectations.
By  (10.2)--(10.3) and Lemma 18, we  have
$$C_1 \exp[-C_2 k_1^{1-\delta/3} k_2\cdots k_{d-1} ]\leq {\bf P}[{\cal D}_1]={\bf P}[\pi({\cal D}'_1)].\eqno{(10.8)}$$
By (10.8) and Lemma 20, there exist $C_i=C_i(F, d, \beta,\delta)$ for $i=1,2,3$ such that for all  $k_1$,
$$
{\bf E}[\tau_r({\bf k}, m, l)\,\,\, |\,\,\, {\cal D}_1]
\leq {\bf E}\tau_{r}({\bf k}, m,l)+ Ck_1^{(1-\delta/8)} k_2\cdots k_{d-1}+C_1\exp\left(-C_2k_1^{1-\delta/4}k_2\cdots k_{d-1}\right).$$
Now we find $C=C(F, d, \beta,\delta)$ such that
$$
{\bf E}[\tau_r({\bf k}, m, l)\,\,\, |\,\,\, {\cal D}_1]
\leq {\bf E}\tau_{r}({\bf k}, m,l)+ Ck_1^{(1-\delta/8)} k_2\cdots k_{d-1}.\eqno{(10.9)}$$
As we mentioned in Lemma 17, 
$${\bf E}[\tau_{r}({\bf k},m,l)\,\,\, |\,\,\, {\cal B}({\bf k}, m,l)]
\leq {\bf E}[\tau_{r}({\bf k}, m)\,\,\, |\,\,\, {\cal B}({\bf k}, m,l)].\eqno{(10.10)}$$
By using Lemma 20 twice, (10.2), and (8.6), 
there exist constants $C_i=C_i(F, d, \beta, \delta)$ for $i=1,2,3$ such that for all large $k_1,\cdots, k_{d-1}$,
\begin{eqnarray*}
&&{\bf E}[\tau_{r}({\bf k}, m,l)]\\
&\leq &{\bf E}[\tau_{r}({\bf k}, m,l)\,\,\, |\,\,\, {\cal B}({\bf k}, m,l)]\\
&&+C_1k_1^{(1-\delta/8)} k_2\cdots k_{d-1}
+ C_1 k_1^{(1-\delta/2)}\exp\left(-C_2 k_1^{(1-\delta/4)} k_2\cdots k_{d-1}\right)\\
&\leq &{\bf E}[\tau_{r}({\bf k}, m)\,\,\, |\,\,\, {\cal B}({\bf k}, m,l)]\\
&&+C_1k_1^{(1-\delta/8)} k_2\cdots k_{d-1}
+ C_1 k_1^{(1-\delta/2)}\exp\left(-C_2 k_1^{(1-\delta/4)} k_2\cdots k_{d-1}\right)\\
&\leq&  {\bf E}\tau_{r}({\bf k}, m)+2C_1k_1^{(1-\delta/8)} k_2\cdots k_{d-1}
+ 2C_1 k_1^{(1-\delta/2)}\exp\left(-C_2 k_1^{(1-\delta/4)} k_2\cdots k_{d-1}\right)\\
&\leq&  {\bf E}\tau_{r}({\bf k}, m)+C_3k_1^{(1-\delta/8)} k_2\cdots k_{d-1}.\hskip 3in (10.11)
\end{eqnarray*} 
Let ${\cal H}({\bf k}, m)$ be the event that 
$$\tau_r({\bf k}, m)\leq 2{\bf E}\tau(e) \|{\bf k}\|_v.$$
Note that  
$$\tau_r({\bf k}, m) \leq \bar{\alpha}({\bf k}, m)\mbox{ and } {\bf E}\bar{\alpha}({\bf k}, m)={\bf E}\tau(e) \|{\bf k}\|_v,$$
where $\bar{\alpha}({\bf k}, m)$ is defined in (6.0).
Thus, by (1.1) and  a standard large deviation estimate,
$${\bf E}\tau_{r}({\bf k}, m)\leq 
{\bf E}[\tau_{r}({\bf k}, m)I({\cal H}({\bf k}, m))]+C_1\exp(-C_2\|{\bf k}\|_v). \eqno{(10.12)}$$
By Theorem 2,
\begin{eqnarray*}
&&{\bf E}\left[\tau_{r}({\bf k}, m)\right]\\
&\leq &{\bf E}\left[\tau_{\min}({\bf k}, m)\right]\\
&&+{\bf E}\left[\tau_{r}({\bf k}, m)I({\cal H}({\bf k}, m))\,\,\, | \,\,\, \bar{N}({\bf k}, m) \geq 
\beta \|{\bf k}\|_v\right]{\bf P}[ \bar{N}({\bf k}, m) \geq \beta \|{\bf k}\|_v]+ 
C_1\exp(-C_2\|{\bf k}\|_v) \\
&\leq & {\bf E}\left[\tau_{\min}({\bf k}, m)\right]+C_3 {\bf E}[\tau(e)] \|{\bf k}\|_v \exp(-C_4 \|{\bf k}\|_v)+C_1\exp(-C_2\|{\bf k}\|_v),
\end{eqnarray*}
where $\bar{N}({\bf k}, m)$ and $\beta$ are defined in Theorem 2.
Therefore,
$${\bf E}[\tau_{r}({\bf k}, m)]\leq {\bf E}[\tau_{\min}({\bf k}, m)]+C_1\exp(-C_2 \|{\bf k}\|_v).\eqno{(10.13)}$$
Together with (10.7)--(10.13), there is $C_1=C_1(F, d, \beta, \delta)$ such that
$${\bf E}(\tau_{r}({\bf k}, m,l)\,\,|\,\,\, {\cal D}_1)\leq {\bf E}\tau_{\min}({\bf k}, m)+C_1k_1^{1-\delta/8} k_2\cdots k_{d-1}.\eqno{(10.14)}$$

Now we  work on the lower bound of  (10.5). 
By the independent discussion and (10.8), 
$${\bf P} \left[{\cal D}_1\cap \pi({\cal D}'_1)\right]=\left({\bf P}[{\cal D}_1]\right)^2\geq C_1^2\exp\left(-2C_2 k_1^{1-\delta/3}k_2\cdots k_{d-1}\right).\eqno{(10.15)}$$
By Lemma 18 and (10.15), and by translation invariance, we may use  the same $C_1$ in (10.14) to have
\begin{eqnarray*}
&&{\bf E}\left[\tau_{r} \left((2k_1-2l+1, k_2,\cdots, k_{d-1}), m,l\right)\,\,\, |\,\,\, {\cal D}_1 \cap \pi({\cal D}'_1) \right]\\
&&\geq  {\bf E}\left[\tau_{r} \left((2k_1-2l+1, k_2,\cdots, k_{d-1}), m,l\right)\right]-
C_1 (2k_1-4l+1)^{(1-\delta/8)} k_2\cdots k_{d-1}\\
&&\geq  {\bf E}\left[\tau_{\min} \left((2k_1-4l+1, k_2,\cdots, k_{d-1}), m\right)\right]-
2C_1 k_1^{(1-\delta/8)} k_2\cdots k_{d-1}.\hskip 1.1in (10.16)
\end{eqnarray*}
Therefore,  by (10.5),  (10.13), and (10.16), 
$$ {\bf E}\left[\tau_{\min}\left((2k_1-4l+1, k_2,\cdots, k_{d-1}), m\right)\right]\leq 2 {\bf E}\tau_{\min}({\bf k}, m)+ 2Ck_1^{(1-\delta/8)} k_2\cdots k_{d-1}. \eqno{(10.17)}$$

We then use the same proof  for ${\bf W}_r''({\bf k}, m,l)$ on ${\bf B}''({\bf k}, m)$. On the 
event ${\cal D}_1\cap \pi({\cal D}'_1)\cap {\cal D}''_1$, we know that (see Fig. 6)
\begin{eqnarray*}
&&{\bf W}_r({\bf k}, m, l)\cup \pi({\bf W}_r '({\bf k}, m, l))\cup {\bf W}_r''({\bf k}, m, l) \mbox{ is a cutset that cuts the bottom  }\\
&&\mbox{ from the top of }[l,3k_1-5l+2]\times [0, k_2]\times\cdots \times [0, k_{d-1}]\times [0, m].
\end{eqnarray*}
By the same discussion from (10.5)--(10.17), there is $C=C(F, d, \beta, \delta)$ such that
$$ {\bf E}\left[\tau_{\min}\left((3k_1-6l+2, k_2,\cdots, k_{d-1}), m\right)\right]\leq 3 {\bf E}\tau_{\min}({\bf k}, m)+3 Ck_1^{(1-\delta/8)} k_2\cdots k_{d-1}.\eqno{(10.18)}$$
With  the same method, (10.2), and Lemmas 18 and 20 by replacing  2 with $w_1$ in (10.5)--(10.17),
we  patch $w_1$  cutsets on adjacent boxes together along the first coordinate to show
$${\bf E}\left[\tau_{\min} \left((w_1(k_1-2l)+w_1, k_2,\cdots, k_{d-1}),m  \right)\right]
\leq w_1{\bf E}\tau_{\min}({\bf k}, m)+Cw_1 k_1^{(1-\delta/8)}k_2\cdots k_{d-1},$$
where $C=C(F, d, \beta, \delta)$ is a constant.
 Note that by (8.5),  $l\leq k_1^{(1-\delta/2})$, so by Lemma 12 (b), for all  large $k_i$ satisfying (10.2)
and $w_1 k_1\leq 2\exp(k_{d-1}^{1-5\delta/6})$,
\begin{eqnarray*}
&&{\bf E}\left[\tau_{\min} \left((w_1( k_1-\lfloor k_1^{(1-\delta/3)}\rfloor), k_2,\cdots, k_{d-1}),m  \right)\right]\\
&&\leq  w_1{\bf E}\tau_{\min}({\bf k}, m)+Cw_1 k_1^{(1-\delta/8)}k_2\cdots k_{d-1}.\hskip 2.4in (10.19)
\end{eqnarray*}
We want to remark that $k_{d-1}$ cannot be arbitrarily larger than  $k_i$ for $i\leq d-2$ in (10.17), 
since we need Lemma 18. So (10.2) is good enough for (10.17).
However, if we work on the $d-1$-th direction, by using Lemma 19, we do not need a restriction for $k_{d-1}$.
More precisely, for all $k_1\leq \cdots \leq k_{d-1}$, 
\begin{eqnarray*}
&&{\bf E}\left[\tau_{\min}\left((k_1, k_2,\cdots k_{d-2}, 2(k_{d-1} -\lfloor k_1^{(1-\delta/3)}\rfloor) ), m\right)\right]\\
&&\leq 2 {\bf E}\tau_{\min}({\bf k}, m)+ 2Ck_1 k_2\cdots k_{d-1}^{(1-\delta/8)}. \hskip 2.6in {(10.20)}
\end{eqnarray*}

We next work on $w_1$ cutsets along the first coordinate and $w_2$ cutsets along the second coordinate.
Along the second coordinate, we have $w_2$ strips with a width $k_2$ for each strip.
We first use (10.19) to patch   $w_1$ cutsets in each strip.
After the first patching, we use the same method of (10.19) to patch $w_2$ patched cutsets in each strip to a cutset.
Note that the size of each cutset, after the first patching, is $w_1( k_1-\lfloor k_1^{(1-\delta/3)}\rfloor)$ along
the first coordinate.
Thus,  Lemma 19 may not be applied for large $w_1 k_1$. So we need to make an extra assumption:
$$w_1k_1\leq 2\exp\left(k_{d-1}^{1-5\delta/6}\right)\leq 2\exp\left(4 k_{2}^{1-5\delta/6}\right).$$
With this assumption, Lemma 19, and the same method of (10.19), we have
\begin{eqnarray*}
&&{\bf E}\left[\tau_{\min}\left( (w_1( k_1- \lfloor k_1^{(1-\delta/3)}\rfloor), w_2( k_2- \lfloor k_2^{(1-\delta/3)} \rfloor),k_3,\cdots, 
k_{d-1}),m  \right)\right]\\
&&\leq  w_1w_2{\bf E}\tau_{\min}({\bf k}, m)+
C\left[w_1 w_2 k_1^{(1-\delta/8)}k_2\cdots k_{d-1}+w_1w_2 k_1k_2^{(1-\delta/8)} \cdots k_{d-1}\right]. \hskip 1cm (10.21)
\end{eqnarray*}
If we continue to iterate this way for the third, ..., the $d-1$-th coordinates, we can show that for integers
$w_1,w_2,\cdots, w_{d-1}$ with $w_j k_j \leq 2\exp(k_{d-1}^{1-5\delta/6})$ for $j=1,\cdots, d-1$, there exists $C=C(F, d, \beta, \delta)$ such that for all large $k_1\leq k_2\leq \cdots \leq k_{d-1}$ with $k_{d-1}\leq 4k_1$, and $m$ that
satisfies (1.17),
\begin{eqnarray*}
&&{\bf E}\left[\tau_{\min} \left((w_1( k_1-\lfloor k_1^{(1-\delta/3)}\rfloor) , w_2( k_2-\lfloor k_2^{(1-\delta/3)}\rfloor),\cdots, w_{d-1}
( k_{d-1}- \lfloor k_{d-1}^{(1-\delta/3)}\rfloor )),m  \right)\right]\\
&&\leq  w_1w_2\cdots w_{d-1}{\bf E}\tau_{\min}({\bf k}, m)\\
&&+  Cw_1w_2\cdots w_{d-1}\left[ k_1^{(1-\delta/8)}k_2\cdots k_{d-1}+ k_1 k_2^{(1-\delta/8)} \cdots k_{d-1}+ \cdots
+k_1k_2\cdots k_{d-2}  k_{d-1}^{(1-\delta/8)} \right].\hskip 0.1cm (10.22)
\end{eqnarray*}

By (10.1), we pick large numbers $k_1',\cdots, k_{d-1}',$ and $m$, given their precise values later, such that for $\epsilon >0$,
$$\left( {{\bf E} \tau_{\min}({\bf k}', m)\over \|{\bf k}'\|_v}\right)\leq \nu_1+\epsilon.\eqno{(10.23)}$$
Now we need to justify the values of these $k_j'$'s such that they satisfy (10.2).
If $ k_{d-1}' \geq 4k_1'$, we may choose 
$0\leq s$ and $0\leq t\leq k_1'$ such that
$$ 2^{s}(k_1'+t)\leq k_{d-1} \leq 2^{s} (k_1'+t+1).\eqno{(10.24)}$$
We divide  $[0, 2^s (k_1'+t)]$, in
the $d-1$-th coordinate, to $2^s$ equal subsegments:
$${\bf D}_1, \cdots, {\bf D}_{2^s}.$$
We consider 
$${\bf T}(j)=[0, k_1']\times [0, k_2']\times \cdots [0, k_{d-2}']\times {\bf D}_j\times [0, m]\mbox{ for } j=1,2,\cdots, 2^s.$$
By using Lemma 12 (a), we know that ${\bf W}({\bf k}, m)\cap {\bf T}(j)$ is a cutset that cuts the bottom from the top of
${\bf T}(j)$. By  translation invariance and Lemma 12 (b), we have
$$2^s {\bf E}\tau_{\min} ((k_1', \cdots, k_{d-2}', (k_1'+t)), m)\leq  {\bf E}\tau_{\min} ({\bf k}', m).\eqno{(10.25)}$$
If we divide (10.25) by  $\|{\bf k}'\|_v$ and use (10.23), for all $k_1'\geq \kappa_1$, then
$$ {\bf E}{\tau_{\min} ((k_1', \cdots, k_{d-2}', (k_1'+t)), m)\over k_1'\cdots k_{d-2}' (k_1'+t)}
\leq {\bf E}{\tau_{\min} ({\bf k}', m)\over \|{\bf k}'\|_v }(1+\epsilon)\leq (\nu_1+\epsilon)(1+\epsilon).\eqno{(10.26)}$$
We use the same argument of (10.26) for the second, ...,  the $d-2$-th coordinates. Thus, by  symmetry,
there are $k_1\leq k_2 \leq \cdots\leq k_{d-1}$ with $k_{d-1} \leq 4 k_{1}$ such that for all $k_1 \geq \kappa_1$,
$$\left( {{\bf E} \tau_{\min}({\bf k}, m)\over \|{\bf k}\|_v}\right)\leq (\nu_1+\epsilon)(1+\epsilon)^d.\eqno{(10.27)}$$
By the assumption in Theorem 3, we can take
$$m \leq \exp(k_{d-1}^{1-\delta}).\eqno{(10.28)}$$

Now we assume that 
$$\lim_{{\bf n}, m'}{{\bf E}\tau_{\min}({\bf n}, m')\over \|{\bf n}\|_v}=\nu_2$$
for a subsequence in $({\bf n}, m')$. 
We select   ${\bf n}=(n_1, \cdots, n_{d-1})$ and $m'$ such that, for
$$ m\leq m'\mbox{ and }2\exp\left (k_{d-1}^{1-5\delta/6}\right) \leq n_j \mbox{ for }j=1,2,\cdots, d-1,$$
$$\nu_2-\epsilon \leq {{\bf E}\tau_{\min}({\bf n}, m')\over \|{\bf n}\|_v}.$$
Also, by symmetry, we take $n_1\leq n_2\leq \cdots \leq n_{d-1}.$
Note that
$${\bf E}\tau_{\min}({\bf n}, m')\leq {\bf E}\tau_{\min}({\bf n}, m) \mbox{ for the $m$ in (10.28)}.\eqno{(10.29)}$$
We assume that 
$$2^{s_j}[\exp(k_{d-1}^{1-5\delta/6})+t_j-1]\leq n_j\leq 2^{s_j}[\exp(k_{d-1}^{1-5\delta/6})+t_j]\eqno{(10.30)}$$
for $1\leq s_j $ and $0\leq t_j \leq \exp(k_{d-1}^{1-5\delta/6})$.
Let 
$$l_j= \exp(k_{d-1}^{1-5\delta/6})+t_j \mbox{ and } {\bf L}=(l_1, \cdots, l_{d-1}).$$
Here we assume that $l_j$ is an integer; otherwise,we just use $\lfloor l_j\rfloor$ to replace $l_j$.
By Lemma 12 (b), for all $k_1\geq \kappa_1$,
$${{\bf E}\tau_{\min}({\bf n}, m)\over \|{\bf n}\|_v}\leq { {\bf E}\tau_{\min}\left((2^{s_1}l_1,
\cdots, 2^{s_{d-1}}l_{d-1}), m\right)\over 2^{s_1+\cdots s_{d-1}} \|{\bf L}\|_v}(1+\epsilon).\eqno{(10.31)}$$
Let 
$$q-\lfloor q^{1-\delta/3}\rfloor= 2^{s_{d-1}-1}l_{d-1}.$$
We may take $k_1 \geq \kappa_2$ such that
$$q\leq 2\times 2^{s_{d-1}-1}l_{d-1}=2^{s_{d-1}}l_{d-1}.$$
Thus 
$$q=2^{s_{d-1}-1}l_{d-1}+ \lfloor q^{1-\delta/3}\rfloor\leq 2^{s_{d-1}-1}l_{d-1}+( 2^{s_{d-1}}l_{d-1})^{1-\delta/3}.$$
Under these observations, by (10.20) and Lemma 12 (b),
\begin{eqnarray*}
&&{\bf E}\tau_{\min}\left((2^{s_1}l_1, \cdots, 2^{s_{d-1}}l_{d-1}), m\right)\\
&&={\bf E}\tau_{\min}\left((2^{s_1}l_1, \cdots, 2(q-\lfloor q^{1-\delta/3}\rfloor)), m\right)\\
&&\leq 2{\bf E}\tau_{\min}\left((2^{s_1}l_1, \cdots, q), m\right)\\
&&+2C \left[2^{s_1+\cdots s_{d-2}} l_1 \cdots l_{d-2} (2^{s_{d-1}}l_{d-1})^{1-\delta/8}\right]\\
&&\leq 
2{\bf E}\tau_{\min}\left((2^{s_1}l_1, \cdots, 2^{s_{d-2}} l_{d-2},  2^{s_{d-1}-1}l_{d-1}+[ 2^{s_{d-1}}l_{d-1}]^{1-\delta/3}), m\right )\\
&&+2C \left[2^{s_1+\cdots s_{d-2}} l_1 \cdots l_{d-2} (2^{s_{d-1}}l_{d-1})^{1-\delta/8}\right].\hskip 2.7in {(10.32)}
\end{eqnarray*}
Note that $m\leq \exp(k_{d-1}^{1-\delta})$ and $l_{d-1} \geq \exp(k_{d-1}^{1-5\delta/6})$, so we may take $k_1\geq \kappa_3$
 such that
$$m\leq \exp(k_{d-1}^{1-\delta})\leq \exp\left({\delta\over 100}k_{d-1}^{1-5\delta/6}\right)\leq \left[l_{d-1}^{1-\delta/3}\right]^{\delta/100}.\eqno{(10.33)}$$
By Lemma 12 (c) and (10.33),
\begin{eqnarray*}
&&2{\bf E}\tau_{\min}\left((2^{s_1}l_1, \cdots, 2^{s_{d-2}} l_{d-2},  2^{s_{d-1}-1}l_{d-1}+(2^{s_{d-1}}l_{d-1})^{1-\delta/3}), m\right )\\
&&\leq 2{\bf E} \tau_{\min} \left((2^{s_1}l_1, \cdots, 2^{s_{d-1}-1}l_{d-1}), m\right)+ 2d ({\bf E}\tau(e)) 2^{s_1+\cdots +s_{d-2}} l_1 \cdots l_{d-2}
(2^{s_{d-1}} l_{d-1} )^{1-\delta/3} m\\
&&\leq 2{\bf E} \tau_{\min} ((2^{s_1}l_1, \cdots, 2^{s_{d-1}-1}l_{d-1}), m)+ 2d ({\bf E}\tau(e))2^{s_1+\cdots +s_{d-2}} l_1 \cdots l_{d-2}
(2^{s_{d-1}} l_{d-1} )^{1-\delta/4}.\hskip 0.1cm (10.34)
\end{eqnarray*}
Together with (10.32) and (10.34),  there exists $C=C(F, d, \beta, \delta)$ such that
\begin{eqnarray*}
&&{\bf E}\tau_{\min}\left((2^{s_1}l_1, \cdots, 2^{s_{d-1}}l_{d-1}), m\right)\\
&&\leq 2{\bf E}\tau_{\min}\left((2^{s_1}l_1, \cdots, 2^{s_{d-1}-1}l_{d-1}), m\right)
+C 2^{-\delta s_{d-1}/8} 2^{s_1+\cdots + s_{d-1}} (l_1 \cdots l_{d-2})l_{d-1}^{1-\delta/8}.\hskip 1cm (10.35)
\end{eqnarray*}

If
 $$2^{s_j}l_j =\max\{ 2^{s_1}l_1, \cdots, 2^{s_{d-1}-1}l_{d-1}\}\mbox{ for } j=1,2,\cdots, d-2,$$ 
then we  continue the process of (10.35) in $j$-th coordinate; otherwise, we still work on the $d-1$-th coordinate.
With this iteration, we can show that for all $k_1\geq \max\{\kappa_2,\kappa_3\}$,
\begin{eqnarray*}
&&{\bf E}\tau_{\min}\left( (2^{s_1}l_1, \cdots, 2^{s_{d-1}}l_{d-1}), m\right)\\
&&\leq 2^{s_1+\cdots + s_{d-1}} {\bf E}\tau_{\min}\left({\bf L}, m\right)
+C \sum_{j=1}^{d-1} \left [2^{s_1+\cdots + s_{d-1}} (l_1 \cdots l_{j-1})l_j^{1-\delta/8} l_{j+1}\cdots l_{d-1}\right]
\left [\sum_{i=1}^{s_j} 2^{-i\delta/8}\right].
\end{eqnarray*}
With these observations and (10.31), for all $k_1\geq \kappa_4$, we have
$${{\bf E}\tau_{\min}({\bf n}, m)\over \|{\bf n}\|_v}\leq { {\bf E}\tau_{\min}({\bf L}, m)\over  \|{\bf L}\|_v}(1+\epsilon).\eqno{(10.36)}$$

Now we need to investigate the relationship between ${\bf E}\tau_{\min}({\bf L}, m)$ and ${\bf E}\tau_{\min}({\bf k}, m)$.
We select  $w_i$ and $r_i$ for $i=1,2,...,d-1$ such that
$$l_i=w_i  (k_i-\lfloor k_i^{(1-\delta/3)}\rfloor) +r_i \mbox{ for } r_i\leq  k_i- \lfloor k_i^{(1-\delta/3)}\rfloor.$$
As we defined,
$${\bf W}\left( w_1(k_1-\lfloor k_1^{(1-\delta/3)}\rfloor), w_2(k_2-\lfloor k_2^{(1-\delta/3)}\rfloor),\cdots, w_{d-1}(k_{d-1}-\lfloor k_{d-1}^{(1-\delta/3)}\rfloor), m \right)$$
is a cutset that cuts the bottom from the top of
 $$[0, w_{1}( k_1- \lfloor k_{1}^{(1-\delta/3)}\rfloor )]\times \cdots \times  [0, w_{d-1}( k_{d-1}-\lfloor k_{d-1}^{(1-\delta/3)}\rfloor )]\times [0, m].$$
By  Lemma 12 (c),  
\begin{eqnarray*}
&& {\bf E}\tau_{\min} ({\bf L}, m)\\
&&\leq {\bf E}\left[\tau_{\min}\left((w_1\lfloor k_1-k_1^{(1-\delta/3)}\rfloor, w_2\lfloor k_2-k_2^{(1-\delta/3)}\rfloor,\cdots, w_{d-1}\lfloor k_{d-1}-k_{d-1}^{(1-\delta/3)}\rfloor) , m\right)\right]\\
&&\hskip 0.3cm + 2d ({\bf E}\tau(e))\left[r_1l_2\cdots l_{d-1} m+ l_1 r_2 l_3 \cdots l_{d-1} m + \cdots + l_1\cdots l_{d-2} r_{d-1} m\right] .\hskip 2cm (10.37)
\end{eqnarray*}

Note that $l_i \leq 2\exp(k_{d-1}^{1-5\delta/6})$ for $i=1,\cdots, d-1$, so by  (10.22) and  (10.37),
\begin{eqnarray*}
&&{\bf E}\tau_{\min} ({\bf L}, m)\\
&\leq & {\bf E}\left[\tau_{\min}\left((w_1\lfloor k_1-k_1^{(1-\delta/3)}\rfloor, w_2\lfloor k_2-k_2^{(1-\delta/3)}\rfloor,\cdots, w_{d-1}\lfloor k_{d-1}-k_{d-1}^{(1-\delta/3)}\rfloor),m \right)\right]\\
&&+2d {\bf E}(\tau(e))[r_1l_2\cdots l_{d-1} m+ l_1 r_2 l_3 \cdots l_{d-1} m + \cdots + l_1\cdots l_{d-2} r_{d-1} m]\\
&\leq & w_1w_2\cdots w_{d-1}{\bf E}\tau_{\min}({\bf k}, m)\\
&&+Cw_1w_2\cdots w_{d-1}\left[ k_1^{(1-\delta/8)}k_2\cdots k_{d-1}+ k_1 k_2^{(1-\delta/8)} \cdots k_{d-1}+ \cdots
+k_1k_2\cdots k_{d-2}  k_{d-1}^{(1-\delta/8)} \right]\\
&&+ 2d{\bf E}(\tau(e))\left[r_1l_2\cdots l_{d-1} m+ l_1 r_2 l_3 \cdots l_{d-1} m + \cdots + l_1\cdots l_{d-2} r_{d-1} m\right]. \hskip 1.8cm (10.38)
\end{eqnarray*}
Therefore, we  divide $\|{\bf L}\|_v$ on  both  sides of (10.38).
Now we work on the left side of (10.38). Note that $l_i\geq \exp(k_{d-1}^{1-5\delta/6})$, so for all $k_1\geq \kappa_5$,
$$ {w_ik_i\over l_i}\leq (1+\epsilon).\eqno{(10.39)}$$
Thus,  the first term in the left side of (10.38),  divided by $\|{\bf L}\|_v$,  is 
$${w_1w_2\cdots w_{d-1}{\bf E}\tau_{\min}({\bf k}, m)\over \|{\bf L}\|_v}\leq {{\bf E}\tau_{\min}({\bf k}, m)\over \|{\bf k}\|_v}(1+\epsilon)^{d-1}.
\eqno{(10.40)}$$
By (10.39), the second sum in the left side of (10.38), divided  by $\|{\bf L}\|_v$, is 
\begin{eqnarray*}
&&{Cw_1w_2\cdots w_{d-1}\left[ k_1^{(1-\delta/8)}k_2\cdots k_{d-1}+ k_1 k_2^{(1-\delta/8)} \cdots k_{d-1}+ \cdots
+k_1k_2\cdots k_{d-2}  k_{d-1}^{(1-\delta/8)} \right]\over \|{\bf L}\|_v} \\
&&\leq  (1+\epsilon)^{d-1}\sum_{i=1}^{d-1} {C\over k^{(1-\delta/8)}}\leq (d-1) \epsilon (1+\epsilon)^{d-1}\hskip 6cm (10.41)
\end{eqnarray*}
for all $k_1\geq \kappa_6$.
Finally, note that $r_i\leq k_i$ and
$$m\leq \exp(k_{d-1}^{1-\delta}) \mbox{ and } \exp(k_{d-1}^{1-5\delta/6})\leq l_j$$
for all $j=1,\cdots, d-1$.
Thus,  for all $k_1\geq \kappa_7$,
$$m r_j \leq r_j\exp(k_{d-1}^{1-\delta}) \leq k_j \exp(k_{d-1}^{1-\delta})\leq \exp(2k_{d-1}^{1-\delta})\leq
\exp(k_{d-1}^{1-5\delta/6}/2) \leq l_j^{1/2}.$$
With this observation, $$l_1 \cdots l_{j-1} r_j l_{j+1} \cdots l_{d-1} m\leq l_1 \cdots l_{j-1} l_j^{1/2} l_{j+1} \cdots l_{d-1}.\eqno{(10.42)}$$
Therefore, the third sum in the left side of
(10.38), divided by $\|{\bf L}\|_v$, is
$${2d{\bf E}(\tau(e))\left[r_1l_2\cdots l_{d-1} m+ l_1 r_2 l_3 \cdots l_{d-1} m + \cdots + l_1\cdots l_{d-2} r_{d-1} m\right]\over \|{\bf L}\|_v}\leq {2d^2\over k_i^{1/2}}\leq C\epsilon\eqno{(10.43)}$$
for all $k_1\geq \kappa_8$.

We now select $k_1 \geq \max\{\kappa_1,\kappa_2,\kappa_3,\kappa_4,\kappa_5,\kappa_6,\kappa_7,\kappa_8\}$ such that (10.27) holds.
Finally, if we put (10.29)--(10.43) together, we show that 
$$\nu_2-\epsilon\leq  \left( {{\bf E} \tau_{\min}({\bf n}, m)\over \|{\bf n}\|_v}\right)\leq \left( {{\bf E} \tau_{\min}({\bf k}, m)\over \|{\bf k}\|_v}\right)(1+\epsilon)^{d-1} +C \epsilon \leq \nu_1(1+\epsilon)^{d-1}+C_1\epsilon.
\eqno{(10.44)}$$
This shows that $\nu_1=\nu_2=\nu$.

Next we need to show the pointwise and $L_1$ convergence. By simply using  a Borel-Cantelli lemma together with the mean
convergence,  the concentration property in (9.7), and (8.1), we have 
$$\lim_{k_1,\cdots, k_{d-1},m\rightarrow \infty} \left( { \tau_{\min}({\bf k}, m)\over \|{\bf k}\|_v}\right)=\nu
\mbox{ a.s. and in }L_1.
\eqno{(10.45)}$$
Therefore, Theorem 3 follows.

\begin{center}
{\bf References}
\end{center}
Aizenman, M., Chayes, J. T., Chayes, L., Frohlich, J., and  Russo, L. (1983). On a sharp transition from area law to perimeter law in a system of random surfaces. {\em Comm. Math. Phys.} {\bf 92} 19--69.\\
Cerf, R. and Theret, M.
(2011). Lower large deviations for the maximal flow
through a domain of $R^d$ in first passage percolation. {\em Probab. Theory Related
Fields} {\bf 150} 635--661.\\
Chayes, L., and Chayes, J. (1986). Bulk transport properties and exponent inequalities for random resistor and flow networks. {\em Comm. Math. Phys.} {\bf 105} 133--152. \\
Engle, E. (1997). {\em  Sperner Theory}. Cambridge University Press, Cambridge, UK.\\
Grimmett, G. (1999). {\em Percolation}. Springer, Berlin.\\
 Grimmett, G.,  and  Kesten, H. (1984).
First-passage percolation, network flows, and electrical
resistances. \emph{Z. Wahrsch. verw. Gebiete}. {\bf 66} 335--366.\\
 Hammersley, J. M., and  Welsh, D. J. A. (1965).
First-passage percolation, subadditive processes,
stochastic networks and generalized renewal theory.
In {\em Bernoulli, Bayse, Laplace Anniversary Volume}
(J. Neyman and L. LeCam, eds.) 61--110. Springer, Berlin.\\
Kesten, H. (1982). {\em Percolation theory for mathematicians}. Birkhauser, Berlin.\\
Kesten, H. (1986). Aspects of first-passage percolation. {\em Lecture Notes in
Math.} {\bf 1180} 126--264. Springer, Berlin.\\
Kesten, H. (1987). Surfaces with minimal random weights and maximal flows: A higher-dimensional version of first-passage percolation. {\em Illinois J. Math.} {\bf 31} 99--166. \\
Kesten, H. (1993). On the speed of convergence in first passage percolation. {\em Ann. Appl. Probab.} {\bf 3} 296--338.\\
Kesten, H., and Zhang, Y. (1990). The probability of a large finite cluster in supercritical Bernoulli percolation. {\em Ann. Probab.} {\bf 18} 537--555.\\
Rossignol, R. and Theret, M.
(2010). Lower large deviations and laws of large
numbers for maximal flows through a box in first passage percolation.
{\em Ann. Inst.
Henri Poincare Probab. Stat.} {\bf 46} 1093--1131.\\
Steele, M., and Zhang, Y. (2003). Nondifferentiability of the time constants of first-passage percolation. 
{\em Ann. Probab.} {\bf 31}  1028--1051.\\
Talagrand, M. (1995). Concentration of measure and isoperimetric inequalities in product spaces. {\em Inst. Hautes Publ. Math. Etudes Sci.}  {\bf 81} 73--205.\\
Zhang, Y. (2000). Critical behavior for maximal flows on the cubic lattice. {\em J. Stat. Phys.} {\bf 98} 799--811. \\
Zhang, Y. (2008). Shape fluctuations are different in  different directions. {\em Ann. Probab.}\\
{\bf 36} 331--362.\\
Yu Zhang\\
Department of Mathematics\\
University of Colorado\\
Colorado Springs, CO 80933\\
yzhang3@uccs.edu\\

\end{document}